\documentclass[11pt]{amsart}
\usepackage[top=1.5in, bottom=1.5in, left=1.4in, right=1.4in]{geometry} 
\usepackage{amssymb}
\usepackage{epstopdf}
\usepackage{amsthm}
\usepackage{hyperref}
\usepackage{array}
\usepackage{mathtools}
\usepackage{enumerate}
\usepackage{tabularx}
\usepackage{ltablex} 
\usepackage{booktabs}
\usepackage{multicol}

\usepackage{verbatim}

\allowdisplaybreaks

\usepackage[all]{xy}
\usepackage{graphicx}  
\usepackage{color, colortbl}
\usepackage[table]{xcolor}
 
\usepackage{tikz}
\usetikzlibrary{decorations.markings}
\usetikzlibrary{arrows}
\usetikzlibrary{arrows.meta}
\usetikzlibrary{bending}
\usetikzlibrary{decorations.pathreplacing,calc}
\usetikzlibrary{decorations.pathmorphing,backgrounds,positioning,fit,matrix}
\usepackage{tkz-graph}

\numberwithin{equation}{section}

\usepackage{afterpage}

\newcommand{\arxiv}[1]{\href{https://arxiv.org/abs/#1}{\texttt{arXiv:#1}}}

\usepackage{etoolbox}
\apptocmd{\sloppy}{\hbadness 10000\relax}{}{}

\newtheorem{thm}{Theorem}[section]
\newtheorem{lemma}[thm]{Lemma}

\newtheorem{cor}[thm]{Corollary}
\newtheorem{prop}[thm]{Proposition}
\newtheorem{conj}[thm]{Conjecture}

\newtheorem{Definition}[thm]{Definition}
\newenvironment{definition}
  {\begin{Definition}\rm}{\end{Definition}}

\newtheorem{Example}[thm]{Example}
\newenvironment{example}
  {\begin{Example}\rm}{\end{Example}}
  
\newtheorem{Remark}[thm]{Remark}
\newenvironment{remark}
  {\begin{Remark}\rm}{\end{Remark}}

\newtheorem{Question}[thm]{Question}

\def\<{\langle}
\def\>{\rangle}

\newcommand\raU[1]{{\xrightarrow[#1]{}}}
\newcommand\raAstU[1]{{\xrightarrow[#1]{\ast}}}
\newcommand\baU[1]{{\xleftrightarrow[#1]{}}}
\newcommand\baAstU[1]{{\xleftrightarrow[#1]{\ast}}}

\def\ra{{\raU{}}}
\def\raAst{{\raAstU{}}}
\def\ba{{\baU{}}}
\def\baAst{{\baAstU{}}}

\def\Span{\mathrm{Span}}

\def\D{D}

\def\scl{0.3}
\def\chipscl{0.6}
\def\chiptikzscl{0.7}
\def\chipcoef{1.2}

\newcommand{\chip}[1]{\begin{tikzpicture}[baseline={([yshift=4pt]current bounding box.south)},block/.style={draw,circle, minimum width={width("11")+12pt},
font=\normalsize,scale=\chipscl}]
\node[block] (A) at (0,0) {$#1$};
\end{tikzpicture}}

\newcommand{\emailhref}[1]{\email{\href{#1}{#1}}}
  
\title{Root system chip-firing {I}: Interval-firing}
\author{Pavel Galashin}
\emailhref{galashin@mit.edu}
\author{Sam Hopkins} 
\emailhref{samuelfhopkins@gmail.com}
\author{Thomas McConville}
\emailhref{thomasmcconvillea@gmail.com }
\author{Alexander Postnikov}
\emailhref{apost@math.mit.edu}
\address{Department of Mathematics, Massachusetts Institute of Technology, Cambridge, MA 02139, USA}

\newboolean{bigimages} 
\setboolean{bigimages}{true}

\begin{document}

\begin{abstract}
Jim Propp recently introduced a variant of chip-firing on a line where the chips are given distinct integer labels. Hopkins, McConville, and Propp showed that this process is confluent from some (but not all) initial configurations of chips. We recast their set-up in terms of root systems: labeled chip-firing can be seen as a \emph{root-firing} process which allows the moves $\lambda \ra \lambda + \alpha$ for $\alpha\in \Phi^{+}$ whenever $\<\lambda,\alpha^\vee\>= 0$, where~$\Phi^{+}$ is the set of positive roots of a root system of Type~A and $\lambda$ is a weight of this root system. We are thus motivated to study the exact same root-firing process for an arbitrary root system. Actually, this \emph{central root-firing} process is the subject of a sequel to this paper. In the present paper, we instead study the \emph{interval root-firing} processes determined by  $\lambda \ra \lambda + \alpha$ for $\alpha\in \Phi^{+}$ whenever~$\<\lambda,\alpha^\vee\>\in [-k-1,k-1]$ or~$\<\lambda,\alpha^\vee\>\in [-k,k-1]$, for any $k \geq 0$. We prove that these interval-firing processes are always confluent, from any initial weight. We also show that there is a natural way to consistently label the stable points of these interval-firing processes across all values of $k$ so that the number of weights with given stabilization is a polynomial in~$k$. We conjecture that these \emph{Ehrhart-like polynomials} have nonnegative integer coefficients.
\end{abstract}

\date{\today}

\maketitle
\tableofcontents
\section{Introduction}
The Abelian Sandpile Model (ASM) is a discrete dynamical system that takes place on a graph. The states of this system are configurations of grains of sand on the vertices of the graph. A vertex with at least as many grains of sand as its neighbors is said to be \emph{unstable}. Any unstable vertex may \emph{topple}, sending one grain of sand to each of its neighbors. The sequence of topplings may continue forever, or it may terminate at a \emph{stable} configuration, where every vertex is stable. The ASM was introduced (in the special case of the two-dimensional square lattice) by the physicists Bak, Tang, and Wiesenfeld~\cite{bak1987self} as a simple model of self-organized criticality; much of the general, graphical theory was subsequently developed by Dhar~\cite{dhar1990self,dhar1999abelian}. The ASM is by now studied in many parts of both physics and pure mathematics: for instance, following the seminal work of Baker and Norine~\cite{baker2007riemann}, it is known that this model is intimately related to tropical algebraic geometry (specifically, divisor theory for tropical curves~\cite{gathmann2008riemann,mikhalkin2008tropical}); meanwhile, the ASM is studied by probabilists because of its remarkable scaling-limit behavior~\cite{pegden2013convergence,levine2016apollonian}; and there are also interesting complexity-theoretic questions related to the ASM, such as, what is the complexity of determining whether a given configuration stabilizes~\cite{kiss2015chip,farrell2016coeulerian}. For more on sandpiles, consult the short survey article~\cite{levine2010sandpile} or the recent textbook~\cite{corry2017divisors}.

Independently of its introduction in the statistical mechanics community, the same model was defined and studied from a combinatorial perspective by Bj\"{o}rner, Lov\'{a}sz, and Shor~\cite{bjorner1991chip} under the name of \emph{chip-firing}.\footnote{It is also worth mentioning that essentially the same model was studied even earlier, in the context of math pedagogy, by Engel~\cite{engel1975probabilistic,engel1976why} under the name of the \emph{probabilistic abacus}.} Instead of grains of sand, we imagine that chips are placed on the vertices of a graph; the operation of an unstable vertex sending one chip to each of its neighbors is now called \emph{firing} that vertex. One fundamental result of  Bj\"{o}rner-Lov\'{a}sz-Shor is that, from any initial chip configuration, either the chip-firing process always goes on forever, or it terminates at a stable configuration that does not depend on the choice of which vertices were fired. This is a \emph{confluence} result: it says that (in the case of termination) the divergent paths in the chip-firing process must come together eventually. This confluence property is the essential property which serves as the basis of all further study of the chip-firing process; it explains the adjective ``Abelian'' in ``Abelian Sandpile Model.''

A closely related chip-firing process to the one studied by Bj\"{o}rner-Lov\'{a}sz-Shor is where a distinguished vertex is chosen to be the \emph{sink}. The sink will never become unstable and is allowed to accumulate any number of chips; hence, any initial chip configuration will eventually stabilize to a unique stable configuration. This model was studied for instance by Biggs~\cite{biggs1999chip} and by Dhar~\cite{dhar1990self,dhar1999abelian}. Chip-firing with a sink has been generalized to several other contexts beyond graphs. One of the most straightforward but also nicest such generalizations is what is called \emph{M-matrix chip-firing} (see e.g.~\cite{gabrielov1993asymmetric,guzman2015chip}, or~\cite[\S13]{postnikov2004trees}). Rather than a graph, we take as input an integer matrix $\mathbf{C} = (\mathbf{C}_{ij})\in \mathbb{Z}^{n\times n}$. The states are vectors $c=(c_1,c_2,\ldots,c_n) \in \mathbb{Z}^n$, and a firing move replaces a state $c$ with $c-\mathbf{C}^te_i$ whenever~$c_i \geq \mathbf{C}_{ii}$ for~$i=1,\ldots,n$. (Here~$e_1,\dots,e_n$ are the standard basis vectors of $\mathbb{Z}^n$; i.e.,~$c-\mathbf{C}^te_i$ is $c$ minus the $i$th row of $\mathbf{C}$.) This firing move is denoted $c \to c-\mathbf{C}^te_i$. Setting $\mathbf{C}$ to be the reduced Laplacian of a graph (including possibly a directed graph, as in~\cite{bjorner1992chip}) recovers chip-firing with a sink. But in fact $\mathbf{C}$ does not need to be a reduced Laplacian of any graph for confluence to hold in this setting: the condition required to guarantee confluence (and termination), as first established by Gabrielov~\cite{gabrielov1993asymmetric}, is that $\mathbf{C}$ be an M-matrix. 

We will discuss M-matrix chip-firing, and its relation to our present research, in more detail later (see~\S\ref{sec:cartanmatrix}). But now let us explain the direct motivation for our work, namely, \emph{labeled} chip-firing.

Bj\"{o}rner, Lov\'{a}sz, and Shor were motivated to introduce the chip-firing process for arbitrary graphs by papers of Spencer~\cite{spencer1986balancing} and Anderson et al.~\cite{anderson1989disks} which studied the special case of chip-firing on a line. Jim Propp recently introduced a version of labeled chip-firing on a line that generalizes this original case. In ordinary chip-firing, the chips are all indistinguishable. But the states of the labeled chip-firing process are configurations of $N$ \emph{distinguishable} chips with integer labels $1,2,\ldots,N$ on the infinite path graph $\mathbb{Z}$. The firing moves consist of choosing two chips that occupy the same vertex and moving the chip with the lesser label one vertex to the right and the chip with the greater label one vertex to the left. Propp conjectured that if one starts with an even number of chips at the origin, this labeled chip-firing process is confluent and in particular the chips always end up in sorted order. Propp's conjecture was recently proved by Hopkins, McConville, and Propp~\cite{hopkins2017sorting}. Note crucially that confluence does not hold for labeled chip-firing if the initial number of chips at the origin is odd (e.g., three). Hence, compared to all the other models of chip-firing discussed above (for which confluence holds locally and follows from Newman's \emph{diamond lemma}~\cite{newman1942theories}), confluence is a much subtler property for labeled chip-firing.  

The crucial observation that motivated our present research is that we can generalize Propp's labeled chip-firing to ``other types,'' as follows. For any configuration of~$N$ labeled chips on the line, if we define the vector $c\coloneqq (c_1,c_2,\ldots,c_N) \in \mathbb{Z}^{N}$ by setting $c_i$ to be the position of the chip with label $i$, then for $i < j$ we are allowed to fire chips with label $i$ and $j$ in this configuration as long as $c$ is orthogonal to $e_i-e_j$; and doing so replaces the vector $c$ by $c+(e_i-e_j)$. Note that the vectors $e_i-e_j$ for~$1 \leq i < j \leq N$ are exactly the positive roots $\Phi^+$ of the root system~$\Phi$ of Type~$A_{N-1}$. 

So there is a natural candidate for a generalization of Propp's labeled chip-firing to arbitrary (crystallographic) root systems: let $\Phi$ be any root system living in some Euclidean vector space~$V$; then for a vector $v\in V$ and a positive root $\alpha \in \Phi^{+}$, we allow the firing move~$v \ra v+\alpha$ whenever $v$ is orthogonal to $\alpha$. We call this process \emph{central root-firing} (or just \emph{central-firing} for short) because we allow a firing move whenever our point $v$ lies on a certain central hyperplane arrangement (namely, the Coxeter arrangement of $\Phi$). 

Central-firing is actually the subject of our sequel paper~\cite{galashin2017rootfiring2}. 

\begin{figure}
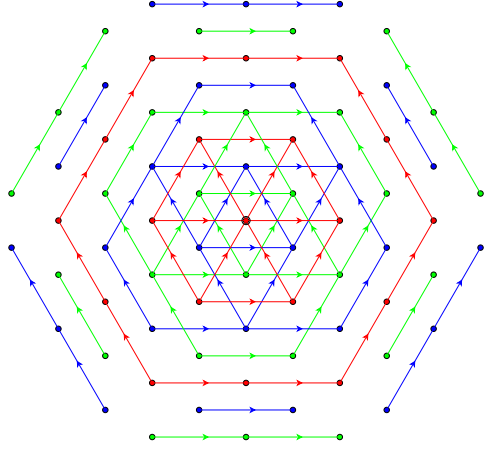

\scalebox{0.12}{
}
\caption{The $k=1$ symmetric interval-firing process for $\Phi=A_2$.} \label{fig:syma2k1}
\end{figure}

In the present paper we instead study two ``affine'' deformations of central-firing. Let us explain what these deformations look like. First of all, it turns out to be best to interpret the condition ``whenever $v$ is orthogonal to~$\alpha$'' as ``whenever~$\<v,\alpha^\vee\>=0$,'' where $\<\cdot,\cdot\>$ is the standard inner product on~$V$ and $\alpha^\vee$ is the \emph{coroot} associated to $\alpha$. Also, rather than consider all vectors $v \in V$ to be the states of our system, it is better to restrict to a discrete setting where the states are \emph{weights} $\lambda \in P$, where~$P$ is the \emph{weight lattice} of~$\Phi$ (this is akin to only allowing vectors $c\in \mathbb{Z}^N$ above). The central-firing moves thus become
\[\lambda \to \lambda+\alpha \textrm{ whenever $\<\lambda,\alpha^\vee\>= 0$ for $\lambda \in P$, $\alpha \in \Phi^{+}$}.\]
The deformations of central-firing we consider involve changing the values of~$\<\lambda,\alpha^\vee\>$ at which we allow the firing move $\lambda \ra \lambda +\alpha$ to be some wider interval. In fact, we study two very particular families of intervals. For $k\in\mathbb{Z}_{\geq 0}$, the \emph{symmetric interval root-firing process} is the binary relation $\raU{\mathrm{sym},k}$ on $P$ defined by
\[\lambda \raU{\mathrm{sym},k} \lambda + \alpha, \; \textrm{ for $\lambda \in P$ and $\alpha\in \Phi^+$ with $\<\lambda,\alpha^\vee\>+ 1 \in \{-k,-k+1,\ldots,k\}$}\]
and the \emph{truncated interval root-firing process} is the relation $\raU{\mathrm{tr},k}$ on $P$ defined by
\[\lambda \raU{\mathrm{tr},k} \lambda + \alpha, \; \textrm{ for $\lambda \in P$ and $\alpha\in \Phi^+$ with $\<\lambda,\alpha^\vee\>+ 1 \in \{-k+1,-k+2,\ldots,k\}$}.\]
We refer to these as \emph{interval-firing} processes for short. 

As mentioned, the central-firing process may or may not be confluent, depending on the initial weight we start at (e.g., our comment about three labeled chips above says that the central-firing process is not confluent from the origin for the root system of Type $A_2$). The first main result of the present paper is the following, which we prove in \S\ref{sec:sym_conf} and~\S\ref{sec:tr_conf}.

\begin{thm} \label{thm:confluence_intro}
For any $k\geq 0$, both the symmetric and truncated interval-firing processes are confluent from all initial weights.
\end{thm}

For example, Figure~\ref{fig:syma2k1} depicts the $k=1$ symmetric interval-firing process for~$\Phi=A_2$: the edges of this graph correspond to firing moves; that this process is confluent means that all paths starting from a given vertex must terminate at the same final vertex. For more such pictures, see Example~\ref{ex:rank2graphs}. 

We call these processes \emph{interval-firing} processes because they allow firing a root from a weight when the inner product of that weight with the corresponding coroot is in some fixed interval. Alternately, we could say that the firing moves are allowed when our weight belongs to a certain affine hyperplane arrangement whose hyperplanes are orthogonal translates of the Coxeter arrangement hyperplanes; this is precisely the sense in which these processes are ``affine.''  The \emph{symmetric} process is so called because the symmetric closure of the relation~$ \raU{\mathrm{sym},k}$ is invariant under the action of the Weyl group. The \emph{truncated} process is so-called because the interval defining it is truncated by one element on the left compared to the symmetric process.

Note that these processes are not truly ``deformations'' of central-firing in the sense that we cannot recover central-firing by specializing~$k$. But observe that the~$k=0$ case of symmetric interval-firing has the firing moves
\[\lambda \raU{\mathrm{sym},0} \lambda+\alpha \textrm{ whenever $\<\lambda,\alpha^\vee\>=-1$ for $\lambda \in P$, $\alpha \in \Phi^{+}$}\]
and the $k=1$ case of truncated interval-firing has the firing moves
\[\lambda \raU{\mathrm{tr},1} \lambda+\alpha \textrm{ whenever $\<\lambda,\alpha^\vee\>\in \{-1,0\}$ for $\lambda \in P$, $\alpha \in \Phi^{+}$}.\]
So these two interval-firing processes are actually very ``close'' to central-firing, and suggest that central-firing (in particular, labeled chip-firing) is somehow right on the ``cusp'' of confluence. Hence, it is not surprising that some of the tools we develop in the present paper are applied to the study of central-firing in the sequel paper~\cite{galashin2017rootfiring2}. We also note that these interval-firing processes themselves have a direct chip-firing interpretation in Type A; see Remark~\ref{rem:chipfiringinterpretation} for more details.

Moreover, we contend that these interval-firing processes are interesting not just because of their connection to central-firing (and hence labeled chip-firing), but also because of their remarkable geometric structure. To get a sense of this geometric structure, the reader is encouraged to look at the depictions of these interval-firing processes for the irreducible rank~$2$ root systems in Example~\ref{ex:rank2graphs}. As we will show, the symmetric and truncated interval-firing processes are closely related to \emph{permutohedra}, and indeed we will mostly investigate these processes from the perspective of convex, polytopal geometry. For example, a key ingredient in our proof of confluence is an exact formula for \emph{traverse lengths of root strings} in permutohedra.

The most striking geometric objects that come out of our investigation of interval-firing are certain ``Ehrhart-like'' polynomials that count the number of weights with given stabilization as we vary our parameter $k$. To make sense of ``with given stabilization,'' first we show that there is a consistent way to label the stable points of the symmetric and truncated interval-firing processes across all values of $k$: these stable points are (a subset of) $\eta^k(\lambda)$ for $\lambda \in P$, where $\eta\colon P\to P$ is a certain piecewise-linear ``dilation'' map depicted in Figure~\ref{fig:eta}. Then we ask: for $\lambda \in P$, how many weights stabilize to $\eta^k(\lambda)$, as a function of~$k$? Let us denote by $L^{\mathrm{sym}}_{\lambda}(k)$ (resp., $L^{\mathrm{tr}}_{\lambda}(k)$) the number of weights $\mu\in P$ that $\raU{\mathrm{sym},k}$-stabilize (resp., $\raU{\mathrm{tr},k}$-stabilize) to $\eta^k(\lambda)$. The following is our second main result, which we prove in \S\ref{sec:sym_Ehrhart} and~\S\ref{sec:tr_Ehrhart}.

\begin{thm}\label{thm:Ehrhart_intro}
  \leavevmode
  \begin{itemize}
  \item For any root system $\Phi$ and any $\lambda \in P$, $L^{\mathrm{sym}}_{\lambda}(k)$ is a polynomial in $k$ with integer coefficients.
  \item For any simply laced root system $\Phi$ and any $\lambda \in P$, $L^{\mathrm{tr}}_{\lambda}(k)$ is a polynomial in $k$ with integer coefficients.
  \end{itemize}
\end{thm}

We conjecture for all root systems~$\Phi$ that these functions are polynomials in~$k$ \emph{with nonnegative integer coefficients}. We call these polynomials \emph{Ehrhart-like} because they count the number of points in some discrete region as it is dilated, but we note that in general the set of weights with given stabilization is not the set of lattice points of any convex polytope, or indeed any convex set (although these Ehrhart-like polynomials do include the usual Ehrhart polynomials of regular permutohedra).

That these Ehrhart-like polynomials apparently have nonnegative integer coefficients suggests that our interval-firing processes may have a deeper connection to the representation theory or algebraic geometry associated to the root system~$\Phi$, although we have no precise idea of what such a connection would be. There is some similarity between our interval-firing processes and the space of \emph{quasi-invariants} of the Weyl group (see~\cite{etingof2003lectures}). We thank Pavel Etingof for pointing this out to us.

As for possible connections to algebraic geometry: one can see in the above definitions of the interval-firing processes that rather than record the intervals corresponding to the values of~$\<\lambda,\alpha^\vee\>$ at which we allow firing, we recorded the intervals corresponding to the values of~$\<\lambda,\alpha^\vee\>+ 1 =\<\lambda+\frac{\alpha}{2},\alpha^\vee\>$ at which we allow firing. This turns out to be more natural in many respects. And with this convention, the intervals defining the symmetric and truncated interval-firing processes are exactly the same as the intervals defining the \emph{(extended) $\Phi^\vee$-Catalan} and  \emph{(extended) $\Phi^\vee$-Shi} hyperplane arrangements~\cite{postnikov2000deformations, athanasiadis2000deformations}. The Catalan and Shi arrangements are known to have many remarkable combinatorial and algebraic properties, such as \emph{freeness}~\cite{edelman1996free, terao2002multiderivations, yoshinaga2004characterization}. Although we have no precise statement to this effect, empirically it seems that many of the remarkable properties of these families of hyperplane arrangements are reflected in the interval-firing processes. See Remark~\ref{rem:hyperplanes} for more discussion of connections with hyperplane arrangements.

Finally, we remark that a kind of ``chip-firing for root systems'' was recently studied by Benkart, Klivans, and Reiner~\cite{benkart2016chip}. However, what Benkart-Klivans-Reiner studied was in fact M-matrix chip-firing with respect to the Cartan matrix $\mathbf{C}$ of the root system~$\Phi$. As we discuss later (see~\S\ref{sec:cartanmatrix}), this Cartan matrix chip-firing is analogous to root-firing \emph{where we only allow firing of the simple roots of $\Phi$}. The root-firing processes we study in this paper allow firing of all the positive roots of~$\Phi$. Hence, our set-up is quite different than the set-up of Benkart-Klivans-Reiner: for instance, the simple roots are always linearly independent, but there are many linear dependencies among the positive roots. Establishing confluence for Cartan matrix chip-firing is easy since the fact that the simple roots are pairwise non-acute implies confluence holds locally; whereas two positive roots may form an acute angle and hence confluence for interval-firing processes is a much more delicate question. Nevertheless, we do explain in Remark~\ref{rem:bkr} how Cartan matrix chip-firing can be obtained from our interval-firing processes by taking a $k\to\infty$ limit.

Now let us outline the rest of the paper. In Part~\ref{part:symtrconfluence} we prove that the symmetric and truncated interval-firing processes are confluent. To do this, we first identify some Weyl group symmetries for both of the interval-firing processes (Theorem~\ref{thm:symmetry}); in particular, we demonstrate that symmetric interval-firing is invariant under the action of the whole Weyl group (explaining its name). We then introduce the map $\eta$ and explain how it labels the stable points for symmetric interval-firing~(Lemma~\ref{lem:symsinks}). We proceed to prove some polytopal results: we establish the aforementioned formula for traverse lengths of permutohedra (Theorem~\ref{thm:traverseformula}); this traverse length formula leads directly to a ``permutohedron non-escaping lemma'' (Lemma~\ref{lem:permtrap}) which says that interval-firing processes get ``trapped'' inside of certain permutohedra. The confluence of symmetric interval-firing (Corollary~\ref{cor:symconfluence}) follows easily from the permutohedron non-escaping lemma. Finally, we establish the confluence of truncated interval-firing (Corollary~\ref{cor:trconfluence}) by first explaining how the map $\eta$ also labels the stable points in the truncated case (Lemma~\ref{lem:trsinks}), and then combining the permutohedron non-escaping lemma with a careful analysis of truncated interval-firing in rank~$2$.

In Part~\ref{part:ehrhart} we study the Ehrhart-like polynomials. We establish the existence of the symmetric Ehrhart-like polynomials (Theorem~\ref{thm:symehrhart}) via some basic Ehrhart theory for zonotopes (see, e.g., Theorem~\ref{thm:polypluszoneehrhart}). Then, to establish the existence of the truncated Ehrhart-like polynomials in the simply laced case (Theorem~\ref{thm:trehrhart}), we study in detail the relationship between symmetric and truncated interval-firing and in particular how the connected components of the graphs of these processes ``decompose'' into smaller connected components in a way consistent with the labeling map  $\eta$ (see~\S\ref{sec:decompose}). In the final section, \S\ref{sec:iterate}, we explain how these Ehrhart-like polynomials also count the sizes of fibers of iterates of a certain operator on the weight lattice, another surprising property of these polynomials that would be worth investigating further.

\medskip

\noindent {\bf Acknowledgements:} We thank Jim Propp, both for several useful conversations and because his introduction of labeled chip-firing and his infectious enthusiasm for exploring its properties launched this project. We also thank the anonymous referee for paying close attention to our article and providing several useful comments. The second author was supported by NSF grant~\#1122374.

\part{Confluence of symmetric and truncated interval-firing} \label{part:symtrconfluence}

\section{Background on root systems} \label{sec:rootsystemdefs}

Here we review the basic facts about root systems we will need in the study of certain vector-firing processes we define in terms of a fixed root system $\Phi$. For details, consult~\cite{humphreys1972lie},~\cite{bourbaki2002lie}, or~\cite{bjorner2005coxeter}.

Fix $V$, an $n$-dimensional real vector space with inner product $\<\cdot,\cdot\>$. For a nonzero vector $\alpha \in V\setminus \{0\}$ we define its \emph{covector} to be $\alpha^\vee \coloneqq  \frac{2\alpha}{\<\alpha,\alpha\>}$. Then we define the \emph{reflection} across the hyperplane orthogonal to $\alpha$ to be the linear map $s_{\alpha}\colon V \to V$ given by $s_{\alpha}(v) \coloneqq  v - \<v,\alpha^\vee\>\alpha$.

\begin{definition}
A \emph{root system} is a finite collection~$\Phi \subseteq V \setminus \{0\}$ of nonzero vectors such that:
\begin{enumerate}
\item $\Span_{\mathbb{R}}(\Phi) = V$;
\item $s_{\alpha}(\Phi) = \Phi$ for all $\alpha \in \Phi$;
\item $\Span_{\mathbb{R}}(\{\alpha\}) \cap \Phi = \{\pm \alpha\}$ for all $\alpha \in \Phi$;
\item $\<\beta,\alpha^\vee\>\in \mathbb{Z}$ for all $\alpha,\beta \in \Phi$.
\end{enumerate} 
We remark that sometimes the third condition is omitted and those root systems satisfying the third condition are called \emph{reduced}. On the other hand, sometimes the fourth condition is omitted and those root systems satisfying the fourth condition are called \emph{crystallographic}. We will assume that all root systems under consideration are reduced and crystallographic and from now on will drop these adjectives.
\end{definition}

From now on in the paper we will fix a root system $\Phi$ in $V$. The vectors $\alpha\in \Phi$ are called \emph{roots}. The dimension of $V$ (which is $n$) is called the \emph{rank} of the root system. The vectors $\alpha^\vee$ for $\alpha\in \Phi$ are called \emph{coroots} and the set of coroots forms another root system, denoted $\Phi^\vee$, in $V$.

We use $W$ to denote the \emph{Weyl group} of $\Phi$, which is the subgroup of $GL(V)$ generated by the reflections $s_{\alpha}$ for $\alpha \in \Phi$. By the first and second conditions of the definition of a root system, $W$ is isomorphic as an abstract group to a subgroup of the symmetric group on $\Phi$, and hence is finite. Observe that the Weyl group of $\Phi^\vee$ is equal to the Weyl group of $\Phi$. Also note that all transformations in $W$ are orthogonal.

It is well-known that we can choose a set $\Delta \subseteq \Phi$ of \emph{simple roots} which form a basis of~$V$, and which divide the root system $\Phi = \Phi^{+} \cup \Phi^{-}$ into \emph{positive} roots $\Phi^{+}$ and \emph{negative} roots $\Phi^{-} \coloneqq  -\Phi^{+}$ so that any positive root $\alpha \in \Phi^{+}$ is a nonnegative integer combination of simple roots. The choice of $\Delta$ is equivalent to the choice of $\Phi^{+}$; one way to choose $\Phi^{+}$ is to choose a generic linear form and let $\Phi^{+}$ be the set of roots which are positive according to this form. There are many choices for $\Delta$ but they are all conjugate under $W$. From now on we will fix a set of simple roots $\Delta$, and thus also a set of positive roots $\Phi^{+}$.  It is known that any $\alpha\in \Phi$ appears in some choice of simple roots (in fact, every $\alpha \in \Phi$ is $W$-conjugate to a simple root appearing with nonzero coefficient in its expansion in terms of simple roots) and hence $W(\Delta) = \Phi$. We use~$\Delta =\{\alpha_1,\ldots,\alpha_n\}$ to denote the simple roots with an arbitrary but fixed order. The coroots $\alpha^\vee_i$ for $i=1,\ldots,n$ are called the \emph{simple coroots} and they of course form a set of simple roots for $\Phi^\vee$. We will always make this choice of simple roots for the dual root system, unless stated otherwise. With this choice of simple roots for the dual root system, we have $(\Phi^\vee)^{+}= (\Phi^{+})^\vee$. 

We use $\mathbf{C} \coloneqq  (\<\alpha_i,\alpha_j^\vee\>) \in \mathbb{Z}^{n\times n}$ to denote the \emph{Cartan matrix} of $\Phi$. Clearly one can recover the root system $\Phi$ from the Cartan matrix $\mathbf{C}$, which is encoded by its Dynkin diagram. The \emph{Dynkin diagram} of $\Phi$ is the graph with vertex set $[n]\coloneqq \{1,2,\dots,n\}$ obtained as follows: first for all $1 \leq i < j \leq n$ we draw $\<\alpha_i,\alpha_j^\vee\>\<\alpha_j,\alpha_i^\vee\>$ edges between~$i$ and $j$; then, if $\<\alpha_i,\alpha_j^\vee\>\<\alpha_j,\alpha_i^\vee\>\notin \{0,1\}$ for some $i$ and $j$, we draw an arrow on top of the edges between them, from $i$ to $j$ if $|\alpha_i| > |\alpha_j|$. If there are no arrows in the Dynkin diagram of $\Phi$ then we say that $\Phi$ is \emph{simply laced}.

There are two important lattices related to $\Phi$, the \emph{root lattice} $Q \coloneqq  \Span_{\mathbb{Z}}(\Phi)$ and the \emph{weight lattice} $P \coloneqq  \{v\in V\colon \<v,\alpha^\vee\>\in \mathbb{Z} \textrm{ for all $\alpha \in \Phi$}\}$. The elements of $P$ are called the \emph{weights} of $\Phi$. By the assumption that~$\Phi$ is crystallographic, we have $Q \subseteq P$. We use $\Omega \coloneqq  \{\omega_1,\ldots,\omega_n\}$ to denote the dual basis to the basis of simple coroots $\{\alpha_1^\vee,\ldots,\alpha_n^\vee\}$ (in other words, the $\omega_i$ are defined by $\<\omega_i,\alpha^\vee_j\>= \delta_{i,j}$); the elements of $\Omega$ are called \emph{fundamental weights}. Observe that $Q = \Span_{\mathbb{Z}}(\Delta)$ and $P = \Span_{\mathbb{Z}}(\Omega)$.

We use $P^{\mathbb{R}}_{\geq 0} \coloneqq  \Span_{\mathbb{R}_{\geq 0}}(\Omega)$, $P_{\geq 0} \coloneqq  \Span_{\mathbb{Z}_{\geq 0}}(\Omega)$ and similarly $Q^{\mathbb{R}}_{\geq 0} \coloneqq  \Span_{\mathbb{R}_{\geq 0}}(\Delta)$, $Q_{\geq 0} \coloneqq  \Span_{\mathbb{Z}_{\geq 0}}(\Delta)$. Note that $P^{\mathbb{R}}_{\geq 0}$ and $Q^{\mathbb{R}}_{\geq 0}$ are dual cones; moreover, because the simple roots are pairwise non-acute, we have $P^{\mathbb{R}}_{\geq 0} \subseteq Q^{\mathbb{R}}_{\geq 0}$. The elements of~$P_{\geq 0}$ are called \emph{dominant weights}. For every $\lambda \in P$ there exists a unique element in $W(\lambda)\cap P_{\geq 0}$ and we use $\lambda_{\mathrm{dom}}$ to denote this element. A dominant weight of great importance is the \emph{Weyl vector} $\rho \coloneqq  \sum_{i=1}^{n}\omega_i$. It is well-known (and easy to check) that $\rho = \frac{1}{2}\sum_{\alpha\in \Phi^{+}}\alpha$.

The connected components of $\{v \in V\colon \<v,\alpha^\vee\>\neq 0 \textrm{ for all $\alpha \in \Phi$}\}$ are called the \emph{chambers} of $\Phi$. The \emph{fundamental chamber} is $C_0 \coloneqq  \{v \in V\colon \<v,\alpha^\vee\>> 0 \textrm{ for all $\alpha \in \Phi^{+}$}\}$. The Weyl group acts freely and transitively on the chambers and hence every chamber is equal to $wC_0$ for some unique $w\in W$. Observe that $P^{\mathbb{R}}_{\geq 0}$ is the closure of $C_0$. 

If $U \subseteq V$ is any subspace spanned by roots, then $\Phi \cap U$ is a root system in $U$, which we call a \emph{sub-root system} of $\Phi$. The root lattice of $\Phi\cap U$ is of course $\Span_{\mathbb{Z}}(\Phi\cap U)$ while the weight lattice is the orthogonal (with respect to~$\< \cdot,\cdot \>$) projection of $P$ onto~$U$. Moreover, $\Phi^{+}\cap U$ is a set of positive roots for $\Phi\cap U$, although $\Delta\cap U$ may not be a set of simple roots for $\Phi\cap U$. We will always consider the positive roots of $\Phi\cap U$ to be $\Phi^{+}\cap U$ unless explicitly stated otherwise. The case of \emph{parabolic sub-root systems} (where in fact $\Delta \cap U$ is a set of simple roots for $\Phi\cap U$) is of special significance: for $I\subseteq [n]$ we set $\Phi_I \coloneqq  \Phi\cap \Span_{\mathbb{R}}(\{\alpha_i\colon i \in I\})$.

\begin{figure}
\begin{tikzpicture}
\node (A) at (0,0) {\begin{tikzpicture}
\def\scl{0.3}
\node[scale=\scl,draw,circle,fill=black] (1) at (-1,0) {};
\node[scale=\scl,draw,circle,fill=black] (2) at (0,0) {};
\node[scale=\scl,draw,circle,fill=black] (3) at (1,0) {};
\node[scale=\scl,draw,circle,fill=black] (4) at (2,0) {};
\node[anchor=north] at (1.south) {$1$};
\node[anchor=north] at (2.south) {$2$};
\node[anchor=north] at (3.south) {${n-1}$};
\node[anchor=north] at (4.south) {$n$};
\draw[thick] (1)--(2);
\draw[thick,dashed] (2)--(3);
\draw[thick] (3)--(4);
\end{tikzpicture}};
\node at (A.south) {$A_n$};
\node (B) at (0,-1.5) {\begin{tikzpicture}[decoration={markings,mark=at position 0.7 with {\arrow{>}}}]
\def\scl{0.3}
\node[scale=\scl,draw,circle] (1) at (-1,0) {};
\node[scale=\scl,draw,circle] (2) at (0,0) {};
\node[scale=\scl,draw,circle] (3) at (1,0) {};
\node[scale=\scl,draw,circle,fill=black] (4) at (2,0) {};
\node[anchor=north] at (1.south) {$1$};
\node[anchor=north] at (2.south) {$2$};
\node[anchor=north] at (3.south) {${n-1}$};
\node[anchor=north] at (4.south) {$n$};
\draw[thick] (1)--(2);
\draw[thick,dashed] (2)--(3);
\draw[thick,double,postaction={decorate}] (3)--(4);
\end{tikzpicture}};
\node at (B.south) {$B_n$};
\node (C) at (0,-3) {\begin{tikzpicture}[decoration={markings,mark=at position 0.7 with {\arrow{>}}}]
\def\scl{0.3}
\node[scale=\scl,draw,circle,fill=black] (1) at (-1,0) {};
\node[scale=\scl,draw,circle] (2) at (0,0) {};
\node[scale=\scl,draw,circle] (3) at (1,0) {};
\node[scale=\scl,draw,circle] (4) at (2,0) {};
\node[anchor=north] at (1.south) {$1$};
\node[anchor=north] at (2.south) {$2$};
\node[anchor=north] at (3.south) {${n-1}$};
\node[anchor=north] at (4.south) {$n$};
\draw[thick] (1)--(2);
\draw[thick,dashed] (2)--(3);
\draw[thick,double,postaction={decorate}] (4)--(3);
\end{tikzpicture}};
\node at (C.south) {$C_n$};
\node (D) at (0,-4.5) {\begin{tikzpicture}
\def\scl{0.3}
\node[scale=\scl,draw,circle,fill=black] (1) at (-1,0) {};
\node[scale=\scl,draw,circle] (2) at (0,0) {};
\node[scale=\scl,draw,circle] (3) at (1,0) {};
\node[scale=\scl,draw,circle,fill=black] (4) at (2,-0.3) {};
\node[scale=\scl,draw,circle,fill=black] (5) at (2,0.3) {};
\node[anchor=north] at (1.south) {$1$};
\node[anchor=north] at (2.south) {$2$};
\node[anchor=north] at (3.south) {${n-2}$};
\node[anchor=north] at (4.south) {${n-1}$};
\node[anchor=south] at (5.north) {$n$};
\draw[thick] (1)--(2);
\draw[thick,dashed] (2)--(3);
\draw[thick] (3)--(4);
\draw[thick] (3)--(5);
\end{tikzpicture}};
\node at (D.south) {$D_n$};
\node (G) at (6,1) {\begin{tikzpicture}[decoration={markings,mark=at position 0.7 with {\arrow{>}}}]
\def\scl{0.3}
\node[scale=\scl,draw,circle] (1) at (-1,0) {};
\node[scale=\scl,draw,circle] (2) at (0,0) {};
\node[anchor=north] at (1.south) {$1$};
\node[anchor=north] at (2.south) {$2$};
\draw[thick,double distance=2pt,postaction={decorate}] (2)--(1);
\draw[thick] (1)--(2);
\end{tikzpicture}};
\node at (G.south) {$G_2$};
\node (F) at (6,-0.5) {\begin{tikzpicture}[decoration={markings,mark=at position 0.7 with {\arrow{>}}}]
\def\scl{0.3}
\node[scale=\scl,draw,circle] (1) at (-1,0) {};
\node[scale=\scl,draw,circle] (2) at (0,0) {};
\node[scale=\scl,draw,circle] (3) at (1,0) {};
\node[scale=\scl,draw,circle] (4) at (2,0) {};
\node[anchor=north] at (1.south) {$1$};
\node[anchor=north] at (2.south) {$2$};
\node[anchor=north] at (3.south) {$3$};
\node[anchor=north] at (4.south) {$4$};
\draw[thick] (1)--(2);
\draw[thick,double,postaction={decorate}] (2)--(3);
\draw[thick] (3)--(4);
\end{tikzpicture}};
\node at (F.south) {$F_4$};
\node (E6) at (6,-2) {\begin{tikzpicture}
\def\scl{0.3}
\node[scale=\scl,draw,circle,fill=black] (1) at (-1,0) {};
\node[scale=\scl,draw,circle] (2) at (1,0.5) {};
\node[scale=\scl,draw,circle] (3) at (0,0) {};
\node[scale=\scl,draw,circle] (4) at (1,0) {};
\node[scale=\scl,draw,circle] (5) at (2,0) {};
\node[scale=\scl,draw,circle,fill=black] (6) at (3,0) {};
\node[anchor=north] at (1.south) {$1$};
\node[anchor=east] at (2.west) {$2$};
\node[anchor=north] at (3.south) {$3$};
\node[anchor=north] at (4.south) {$4$};
\node[anchor=north] at (5.south) {$5$};
\node[anchor=north] at (6.south) {$6$};
\draw[thick] (1)--(3)--(4)--(5)--(6);
\draw[thick] (2)--(4);
\end{tikzpicture}};
\node at (E6.south) {$E_6$};
\node (E7) at (6,-3.5) {\begin{tikzpicture}
\def\scl{0.3}
\node[scale=\scl,draw,circle] (1) at (-1,0) {};
\node[scale=\scl,draw,circle] (2) at (1,0.5) {};
\node[scale=\scl,draw,circle] (3) at (0,0) {};
\node[scale=\scl,draw,circle] (4) at (1,0) {};
\node[scale=\scl,draw,circle] (5) at (2,0) {};
\node[scale=\scl,draw,circle] (6) at (3,0) {};
\node[scale=\scl,draw,circle,fill=black] (7) at (4,0) {};
\node[anchor=north] at (1.south) {$1$};
\node[anchor=east] at (2.west) {$2$};
\node[anchor=north] at (3.south) {$3$};
\node[anchor=north] at (4.south) {$4$};
\node[anchor=north] at (5.south) {$5$};
\node[anchor=north] at (6.south) {$6$};
\node[anchor=north] at (7.south) {$7$};
\draw[thick] (1)--(3)--(4)--(5)--(6)--(7);
\draw[thick] (2)--(4);
\end{tikzpicture}};
\node at (E7.south) {$E_7$};
\node (E8) at (6,-5) {\begin{tikzpicture}
\def\scl{0.3}
\node[scale=\scl,draw,circle] (1) at (-1,0) {};
\node[scale=\scl,draw,circle] (2) at (1,0.5) {};
\node[scale=\scl,draw,circle] (3) at (0,0) {};
\node[scale=\scl,draw,circle] (4) at (1,0) {};
\node[scale=\scl,draw,circle] (5) at (2,0) {};
\node[scale=\scl,draw,circle] (6) at (3,0) {};
\node[scale=\scl,draw,circle] (7) at (4,0) {};
\node[scale=\scl,draw,circle] (8) at (5,0) {};
\node[anchor=north] at (1.south) {$1$};
\node[anchor=east] at (2.west) {$2$};
\node[anchor=north] at (3.south) {$3$};
\node[anchor=north] at (4.south) {$4$};
\node[anchor=north] at (5.south) {$5$};
\node[anchor=north] at (6.south) {$6$};
\node[anchor=north] at (7.south) {$7$};
\node[anchor=north] at (8.south) {$8$};
\draw[thick] (1)--(3)--(4)--(5)--(6)--(7)--(8);
\draw[thick] (2)--(4);
\end{tikzpicture}};
\node at (E8.south) {$E_8$};
\end{tikzpicture}
\caption{Dynkin diagrams of all irreducible root systems. The nodes corresponding to minuscule weights are filled in.} \label{fig:dynkinclassification}
\end{figure}
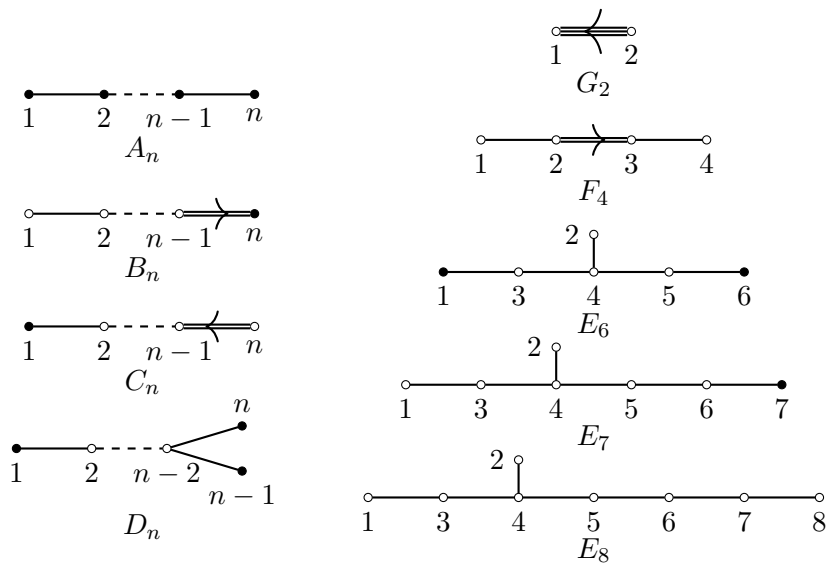

If there exists an orthogonal decomposition $V = V_1 \oplus V_2$ with $0 \subsetneq V_1,V_2 \subsetneq V$ such that $\Phi = \Phi_1 \cup \Phi_2$ with $\Phi_i \subseteq V_i$ for $i=1,2$, then we write $\Phi = \Phi_1\oplus\Phi_2$ and we say the root system $\Phi$ is \emph{reducible}. Otherwise we say that it is \emph{irreducible}. (Let us also declare by fiat that the empty set, although it is a root system, is not irreducible.) In other words, a root system is irreducible if and only if its Dynkin diagram is connected. The famous \emph{Cartan-Killing classification} classifies all irreducible root systems up to isomorphism, where an \emph{isomorphism} of root systems is a bijection between roots induced from an invertible orthogonal map, potentially composed with a global rescaling of the inner product. Figure~\ref{fig:dynkinclassification} shows the Dynkin diagrams of all the irreducible root systems: these are the classical infinite series $A_n$ for $n\geq 1$, $B_n$ for $n \geq 2$, $C_n$ for $n \geq 3$, $D_n$ for $n \geq 4$, together with the exceptional root systems $G_2$, $F_4$, $E_6$, $E_7$, and $E_8$. Our numbering of the simple roots is consistent with Bourbaki~\cite{bourbaki2002lie}. In every case the subscript in the name of the root system denotes the number of nodes of the Dynkin diagram, which is also the number of simple roots, that is, the rank of~$\Phi$. These labels $A_n$, $B_n$, etc. are the \emph{type} of the root system; we may also talk about, e.g., ``Type A'' root systems.

All constructions that depend on the root system $\Phi$ decompose in a simple way as a direct product of irreducible factors. Hence without loss of generality we will from now on {\bf assume that $\Phi$ is irreducible.}

In an irreducible root system, there are at most two values of lengths $|\alpha|$ among the roots $\alpha \in \Phi$. Those roots whose lengths achieve the maximum value are called \emph{long}, and those which do not are called \emph{short}. The Weyl group $W$ acts transitively on the long roots, and it also acts transitively on the short roots.

There is a natural partial order on $P$ called the \emph{root order} whereby $\mu\leq \lambda$ for $\mu,\lambda \in P$ if $\lambda-\mu \in Q_{\geq 0}$. When restricted to $\Phi^{+}$, this partial order is graded by \emph{height}; the height of $\alpha = \sum_{i=1}^{n}c_i\alpha_i \in \Phi$ is $\sum_{i=1}^{n}c_i$. Because we have assumed that~$\Phi$ is irreducible, there is a unique maximal element of $\Phi^{+}$ according to root order, denoted $\theta$ and called the \emph{highest root}. The highest root is always long. We use $\widehat{\theta}$ to denote the unique (positive) root such that $\widehat{\theta}^\vee$ is the highest root of the dual root system $\Phi^\vee$ (with respect to the choice of $\{\alpha_1^\vee,\ldots,\alpha_n^\vee\}$ as simple roots). If $\Phi$ is simply laced then $\theta=\widehat{\theta}$ and $\theta$ is the unique root which is a dominant weight; if~$\Phi$ is not simply laced then $\theta$ and $\widehat{\theta}$ are the two roots which are dominant weights. In the non-simply laced case we call $\widehat{\theta}$ the \emph{highest short root}: it is the maximal short root with respect to the root ordering.

The root lattice $Q$ is a full rank sublattice of $P$; hence the quotient $P/Q$ is some finite abelian group. Note that $P/Q \simeq \mathrm{coker}(\mathbf{C}^t)$ where we view the transposed matrix as a map $\mathbf{C}^t\colon \mathbb{Z}^n\to \mathbb{Z}^n$. The order of this group is called the \emph{index of connection} of~$\Phi$ and is denoted $f\coloneqq |P/Q|$. There is a nice choice of coset representatives of $P/Q$, which we now describe. A dominant, nonzero weight $\lambda \in P_{\geq 0}\setminus \{0\}$ is called \emph{minuscule} if $\<\lambda, \alpha^\vee\>\in \{-1,0,1\}$ for all $\alpha \in \Phi$. Let us use $\Omega_m$ to denote the set of minuscule weights. Note that $\Omega_m\subseteq \Omega$, i.e., a minuscule weight must be a fundamental weight. In Figure~\ref{fig:dynkinclassification}, the vertices corresponding to minuscule weights are filled in. In fact, there are $f-1$ minuscule weights and the minuscule weights together with zero form a collection of coset representatives of $P/Q$. We use $\Omega^{0}_{m} \coloneqq  \Omega_m \cup\{0\}$ to denote the set of these representatives. 

There is another characterization of minuscule weights that we will find useful. Namely, for a dominant weight $\lambda \in P_{\geq 0}$ we have that $\lambda \in \Omega^{0}_m$ if and only if $\lambda$ is the minimal element according to root order in $(Q+\lambda) \cap P_{\geq 0}$. 

This last characterization of minuscule weight can also be described in terms of certain polytopes called \emph{($W$)-permutohedra}. Permutohedra will play a key role for us in our understanding of interval-firing processes, so let us review these now. For~$v \in V$, we define the \emph{permutohedron} associated to~$v$ to be~$\Pi(v) \coloneqq  \mathrm{ConvexHull}\, W(v)$, a convex polytope in~$V$.  And for a weight~$\lambda \in P$, we define $\Pi^Q(\lambda) \coloneqq  \Pi(\lambda)\cap(Q+\lambda)$, which we call the \emph{discrete permutohedron} associated to~$\lambda$. 

The following simple proposition describes the containment of permutohedra (see also~\cite[1.2]{stembridge1998partial}):

\begin{prop} \label{prop:perm_containment}
For $u,v\in P_{\geq 0}^{\mathbb{R}}$ we have $\Pi(u)\subseteq \Pi(v)$ if and only if $v-u\in Q_{\geq 0}^{\mathbb{R}}$. Hence for $\mu,\lambda \in P_{\geq 0}$ we have $\Pi^Q(\mu)\subseteq \Pi^Q(\lambda)$ if and only if $\mu \leq \lambda$ (in root order).
\end{prop}
\begin{proof}
First suppose that $u$ and $v$ are strictly inside the fundamental chamber~$C_0$, i.e., that we have $\<u,\alpha_i^\vee\>> 0$ and $\<v,\alpha_i^\vee\>> 0$ for all $i\in[n]$. By the \emph{inner cone} of polytope at a vertex, we mean the affine convex cone spanned by the edges of the polytope incident to that vertex in the direction ``outward'' from that vertex. Note that a point belongs to a polytope if and only if it belongs to the inner cone of that polytope at every vertex. Since the walls of the fundamental chamber are orthogonal to the simple roots, it is easy to see that if $u$ and $v$ are strictly inside the fundamental chamber then the inner cone of $\Pi(u)$ at $u$ is spanned by the negatives of the simple roots, and ditto for the inner cone of $\Pi(v)$ and $v$. So if we do not have $v-u\in Q_{\geq 0}^{\mathbb{R}}$, then clearly $u$ does not belong to $\Pi(v)$. Hence suppose that $v-u\in Q_{\geq 0}^{\mathbb{R}}$. Every vertex of $\Pi(u)$ belongs to the inner cone of $\Pi(u)$ at $u$; i.e., $u-u' \in Q_{\geq 0}^{\mathbb{R}}$  for all $u'\in W(u)$. Thus for all $u'\in W(u)$ we have $v-u' \in Q_{\geq 0}^{\mathbb{R}}$; i.e., every point in $\Pi(u)$ is in the inner cone of $\Pi(v)$ at $v$. But then by the $W$-invariance of permutohedra, we conclude that every point in $\Pi(u)$ is in the inner cone of $\Pi(v)$ at every vertex of $\Pi(v)$, and hence that~$\Pi(u)\subseteq \Pi(v)$, as claimed.

For arbitrary $u,v\in P_{\geq 0}^{\mathbb{R}}$, note $\Pi(u) = \bigcap_{\varepsilon > 0} \Pi(u+\varepsilon \rho)$ and $\Pi(v) = \bigcap_{\varepsilon > 0} \Pi(v+\varepsilon \rho)$, and $u+\varepsilon \rho$ and $v+\varepsilon \rho$ will be strictly inside the fundamental chamber for all $\varepsilon > 0$. Thus the result for arbitrary $u,v\in P_{\geq 0}^{\mathbb{R}}$ follows from the preceding paragraph.
\end{proof}

So in light of Proposition~\ref{prop:perm_containment}, we see that minuscule weights can also be characterized as follows: for~$\lambda \in P_{\geq 0}$ we have~$\lambda \in \Omega^{0}_m$ if and only if $\Pi^Q(\lambda)=W(\lambda)$. For references for all these various characterizations of and facts about minuscule weights, see~\cite[Proposition 3.10]{benkart2016chip} (who in particular credit Stembridge~\cite{stembridge1998partial} for some of these facts).

\section{Background on binary relations and confluence} \label{sec:relations}

Interval-firing will formally be defined to be a binary relation on the weight lattice of~$\Phi$. Before giving the precise definition, we review some general notation and results concerning binary relations. Let $X$ be a set and $\ra$ a binary relation on $X$. We use $\Gamma_{\ra}$ to denote the directed graph (from now on, ``digraph'') with vertex set $X$ and with a directed edge~$(x,y)$ whenever $x \ra y$. Clearly~$\Gamma_{\ra}$ contains exactly the same information as~$\ra$ and we will often implicitly identify binary relations and digraphs (specifically, digraphs without multiple edges in the same direction) in this way. We use $\raAst$ to denote the reflexive, transitive closure of~$\ra$: that is, we write~$x \raAst y$ to mean that~$x = x_0 \ra x_1 \ra \cdots \ra x_k = y$ for some~$k \in \mathbb{Z}_{\geq 0}$. In other words, $x \raAst y$ means there is a path from $x$ to $y$ in~$\Gamma_{\ra}$. We use~$\ba$ to denote the symmetric closure of~$\ra$: $x \ba y$ means that~$x \ra y$ or~$y \ra x$. For any digraph $\Gamma$, we use $\Gamma^{\mathrm{un}}$ to denote the underlying undirected graph of $\Gamma$; in fact, we view~$\Gamma^{\mathrm{un}}$ as a digraph: it has edges $(x,y)$ and $(y,x)$ whenever $(x,y)$ is an edge of~$\Gamma$. Hence $\Gamma_{\ba} = \Gamma^{\mathrm{un}}_{\ra}$. Finally, we use $\baAst$ to denote the reflexive, transitive, symmetric closure of~$\ra$: $x \baAst y$ means that~$x = x_0 \ba x_1 \ba \cdots \ba x_k = y$ for some~$k \in \mathbb{Z}_{\geq 0}$. In other words, $x \baAst y$ means there is a path from $x$ to $y$ in~$\Gamma^{\mathrm{un}}_{\ra}$.

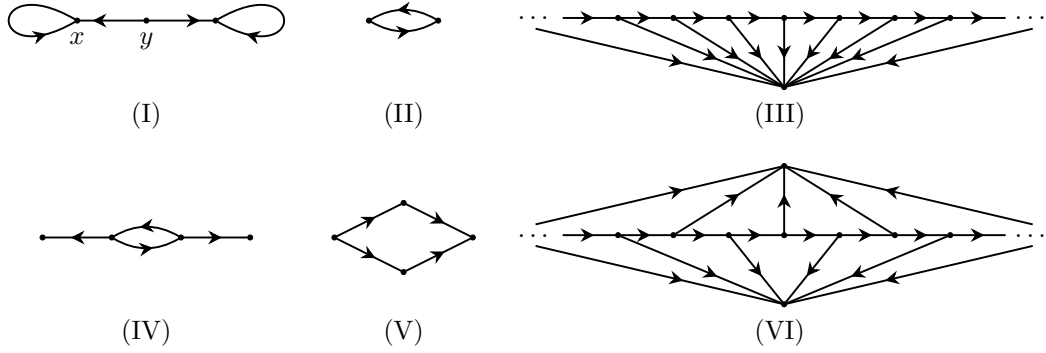
\begin{figure}
\scalebox{0.92}{
\begin{tabular}{ccc}

\begin{tikzpicture}[>={Stealth[width=2mm,length=2mm,bend]},decoration={markings,mark=at position 0.8 with {\arrow{>}}}]
\def\scl{0.2}
  \draw[white] (0,-1)--(0,0);
\node[scale=\scl,draw,circle,fill=black] (X) at (-1,0) {};
\node[anchor=north] at (X.south) {$x$};
\node[scale=\scl,draw,circle,fill=black] (Y) at (0,0) {};
\node[anchor=north] at (Y.south) {$y$};
\node[scale=\scl,draw,circle,fill=black] (Z) at (1,0) {};
\draw[postaction={decorate},thick] (Y)--(X);
\draw[postaction={decorate},thick] (X) to [out=150,in=210,looseness=100] (X);
\draw[postaction={decorate},thick] (Y)--(Z);
\draw[postaction={decorate},thick] (Z) to [out=30,in=-30,looseness=100] (Z);
\end{tikzpicture}
&
\begin{tikzpicture}[>={Stealth[width=2mm,length=2mm,bend]},decoration={markings,mark=at position 0.6 with {\arrow{>}}}]
  \def\scl{0.2}
  \draw[white] (0,-1)--(0,0);
\node[scale=\scl,draw,circle,fill=black] (Y) at (-1,0) {};
\node[scale=\scl,draw,circle,fill=black] (Z) at (0,0) {};
\draw[postaction={decorate},thick] (Y) to[bend right] (Z);
\draw[postaction={decorate},thick] (Z) to[bend right] (Y);
\end{tikzpicture}
  &
    \def\hordist{0.8}
\begin{tikzpicture}[>={Stealth[width=2mm,length=2mm,bend]},decoration={markings,mark=at position 0.6 with {\arrow{>}}}]
  \def\scl{0.2}
  \draw[white] (0,-1)--(0,0);
  \node[scale=1,anchor=east] (A) at ({-1*\hordist},0) {$\dots$};
  \node[scale=1,anchor=west] (B) at ({7*\hordist},0) {$\dots$};
  \node[scale=\scl,draw,circle,fill=black] (U) at ({3*\hordist},-1) {};
\foreach \i in {0,1,...,6}{
  \node[scale=\scl,draw,circle,fill=black] (X\i) at ({\i*\hordist},0) {};
  \draw[postaction={decorate},thick] (X\i) -- (U);
}
  \draw[postaction={decorate},thick] (A) -- (X0);
  \draw[postaction={decorate},thick] (X0) -- (X1);
  \draw[postaction={decorate},thick] (X1) -- (X2);
  \draw[postaction={decorate},thick] (X2) -- (X3);
  \draw[postaction={decorate},thick] (X3) -- (X4);
  \draw[postaction={decorate},thick] (X4) -- (X5);
  \draw[postaction={decorate},thick] (X5) -- (X6);
  \draw[postaction={decorate},thick] (X6) -- (B);
  \draw[postaction={decorate},thick] (A.south) -- (U);
  \draw[postaction={decorate},thick] (B.south) -- (U);
\end{tikzpicture}
  \\
  (I) & (II) & (III)\\ &&\\
\begin{tikzpicture}[>={Stealth[width=2mm,length=2mm,bend]},decoration={markings,mark=at position 0.6 with {\arrow{>}}}]
\def\scl{0.2}
  \draw[white] (0,-1)--(0,0);
\node[scale=\scl,draw,circle,fill=black] (X) at (-2,0) {};
\node[scale=\scl,draw,circle,fill=black] (Y) at (-1,0) {};
\node[scale=\scl,draw,circle,fill=black] (Z) at (0,0) {};
\node[scale=\scl,draw,circle,fill=black] (W) at (1,0) {};
\draw[postaction={decorate},thick] (Y) -- (X);
\draw[postaction={decorate},thick] (Y) to[bend right] (Z);
\draw[postaction={decorate},thick] (Z) to[bend right] (Y);
\draw[postaction={decorate},thick] (Z) -- (W);
\end{tikzpicture}
&
\begin{tikzpicture}[>={Stealth[width=2mm,length=2mm,bend]},decoration={markings,mark=at position 0.6 with {\arrow{>}}}]
\def\scl{0.2}
  \draw[white] (0,-1)--(0,0);
\node[scale=\scl,draw,circle,fill=black] (X) at (-1,0) {};
\node[scale=\scl,draw,circle,fill=black] (Y) at (0,0.5) {};
\node[scale=\scl,draw,circle,fill=black] (Z) at (0,-0.5) {};
\node[scale=\scl,draw,circle,fill=black] (W) at (1,0) {};
\draw[postaction={decorate},thick] (X) -- (Y);
\draw[postaction={decorate},thick] (X) -- (Z);
\draw[postaction={decorate},thick] (Y) -- (W);
\draw[postaction={decorate},thick] (Z) -- (W);
\end{tikzpicture}
  &
    \def\hordist{0.8}
\begin{tikzpicture}[>={Stealth[width=2mm,length=2mm,bend]},decoration={markings,mark=at position 0.6 with {\arrow{>}}}]
  \def\scl{0.2}
  \node[scale=1,anchor=east] (A) at ({-1*\hordist},0) {$\dots$};
  \node[scale=1,anchor=west] (B) at ({7*\hordist},0) {$\dots$};
  \node[scale=\scl,draw,circle,fill=black] (U) at ({3*\hordist},-1) {};
  \node[scale=\scl,draw,circle,fill=black] (V) at ({3*\hordist},1) {};
\foreach \i in {0,2,...,6}{
  \node[scale=\scl,draw,circle,fill=black] (X\i) at ({\i*\hordist},0) {};
  \draw[postaction={decorate},thick] (X\i) -- (U);

}
\foreach \i in {1,3,...,5}{
  \node[scale=\scl,draw,circle,fill=black] (X\i) at ({\i*\hordist},0) {};
  \draw[postaction={decorate},thick] (X\i) -- (V);
}
  \draw[postaction={decorate},thick] (A) -- (X0);
  \draw[postaction={decorate},thick] (X0) -- (X1);
  \draw[postaction={decorate},thick] (X1) -- (X2);
  \draw[postaction={decorate},thick] (X2) -- (X3);
  \draw[postaction={decorate},thick] (X3) -- (X4);
  \draw[postaction={decorate},thick] (X4) -- (X5);
  \draw[postaction={decorate},thick] (X5) -- (X6);
  \draw[postaction={decorate},thick] (X6) -- (B);
  \draw[postaction={decorate},thick] (A.south) -- (U);
  \draw[postaction={decorate},thick] (B.south) -- (U);
  \draw[postaction={decorate},thick] (A.north) -- (V);
  \draw[postaction={decorate},thick] (B.north) -- (V);
\end{tikzpicture}\\
  (IV) & (V) & (VI)\\
\end{tabular}
}
\caption{Examples of various relations: (I) is confluent from $x$ but not from $y$; (II) and (III) are confluent but not terminating; (IV) and~(VI) are locally confluent but not confluent; (V) is confluent and terminating.} \label{fig:relationexs}
\end{figure}

Now let us review some notions of confluence for binary relations. Here we generally follow standard terminology in the theory of abstract rewriting systems, as laid out for instance in~\cite{huet1980confluent}; however, instead following chip-firing terminology, we use ``stable'' in place of what would normally be called ``irreducible,'' and rather than ``normal forms'' we refer to ``stabilizations.'' We say that~$\ra$ is \emph{terminating} (also sometimes called \emph{noetherian}) if  there is no infinite sequence of relations~$x_0 \ra x_1 \ra x_2 \ra \cdots$; i.e., $\ra$ is terminating means that~$\Gamma_{\ra}$ has no infinite paths (which implies in particular that this digraph has no directed cycles). Generally speaking, the relations we are most interested in will all be terminating and it will be easy for us to establish that they are terminating. For~$x \in X$, we say that~$\ra$ is \emph{confluent from $x$} if whenever $x \raAst y_1$ and~$x \raAst y_2$, there is $y_3$ such that~$y_1 \raAst y_3$ and~$y_2 \raAst y_3$. We say~$x \in X$ is \emph{$\ra$-stable} (or just \emph{stable} if the context is clear) if there is no $y \in X$ with $x \ra y$. In graph-theoretic language, $x$ is~$\ra$-stable means that~$x$ is a sink (vertex of outdegree zero) of $\Gamma_{\ra}$. If $\ra$ is terminating, then for every~$x \in X$ there must be at least one stable $y \in X$ with $x \raAst y$. On the other hand, if $\ra$ is confluent from~$x \in X$, then there can be at most one stable $y \in X$ with $x \raAst y$. Hence if $\ra$ is terminating and is confluent from $x$, then there exists a unique stable~$y$ with~$x \raAst y$; we call this~$y$ the~\emph{$\ra$-stabilization} (or just \emph{stabilization} if the context is clear) of~$x$. We say that~$\ra$ is \emph{confluent} if it is confluent from every~$x \in X$. As we just explained, if $\ra$ is confluent and terminating then a unique stabilization of~$x$ exists for all $x \in X$. A weaker notion than confluence is that of local confluence: we say that $\ra$ is \emph{locally confluent} if for any~$x \in X$, if $x \ra y_1$ and $x \ra y_2$, then there is some $y_3$ with~$y_1 \raAst y_3$ and~$y_2 \raAst y_3$. Figure~\ref{fig:relationexs} gives some examples of relations comparing these various notions of confluence and termination. Observe that there is no example in this figure of a relation that is locally confluent and terminating but not confluent. That is no coincidence: Newman's lemma, a.k.a.~the diamond lemma, says that local confluence plus termination implies confluence.

\begin{lemma}[Diamond lemma, see~{\cite[Theorem 3]{newman1942theories} or~\cite[Lemma 2.4]{huet1980confluent}}]\label{lem:diamond}
Suppose $\ra$ is terminating. Then $\ra$ is confluent if and only if it is locally confluent.
\end{lemma}

\section{Definition of interval-firing}

In this section we formally define the interval-firing processes in their most general form. We use the notation $\mathbf{k} \in \mathbb{Z}[\Phi]^{W}$ to mean that $\mathbf{k}$ is an integer-valued function on the roots of~$\Phi$ that is invariant under the action of the Weyl group. We write $\mathbf{a} \leq \mathbf{b}$ to mean that~$\mathbf{a}(\alpha) \leq \mathbf{b}(\alpha)$ for all $\alpha \in \Phi$. We use the notation~$\mathbf{k} = k$ to mean that $\mathbf{k}$ is constantly equal to~$k$. We also use the obvious notation $a\mathbf{a}+b\mathbf{b}$ for linear combinations of these functions. We use $\mathbb{N}[\Phi]^W$ to denote the set of~$\mathbf{k} \in \mathbb{Z}[\Phi]^W$ with $\mathbf{k} \geq 0$. We write $\rho_{\mathbf{k}} \coloneqq  \sum_{i=1}^{n}\mathbf{k}(\alpha_i)\omega_i$. Since we have assumed that~$\Phi$ is irreducible, there are at most two $W$-orbits of $\Phi$: the short roots and the long roots. If $\Phi$ is simply laced then it has a single Weyl group orbit and~$\mathbf{k} = k$ for some constant~$k \in \mathbb{Z}$; otherwise, we have two constants $k_{s}, k_{l} \in \mathbb{Z}$ so that $\mathbf{k}(\alpha)=k_s$ if $\alpha$ is short and~$\mathbf{k}(\alpha) = k_l$ if $\alpha$ is long.

For $\mathbf{k}\in\mathbb{N}[\Phi]^W$, the \emph{symmetric interval-firing process} is the binary relation $\raU{\mathrm{sym},\mathbf{k}}$ on $P$ defined by
\[\lambda \raU{\mathrm{sym},\mathbf{k}} \lambda + \alpha, \; \textrm{ for $\lambda \in P$ and $\alpha\in \Phi^+$ with $\<\lambda+\frac{\alpha}{2},\alpha^\vee\>\in [-\mathbf{k}(\alpha),\mathbf{k}(\alpha)]$}\]
and the \emph{truncated interval-firing process} is the binary relation $\raU{\mathrm{tr},\mathbf{k}}$ on $P$ defined by
\[\lambda \raU{\mathrm{tr},\mathbf{k}} \lambda + \alpha, \; \textrm{ for $\lambda \in P$ and $\alpha\in \Phi^+$ with $\<\lambda+\frac{\alpha}{2},\alpha^\vee\>\in [-\mathbf{k}(\alpha)+1,\mathbf{k}(\alpha)]$}.\]

From now own we will often think about a relation $\ra$ as $\Gamma_{\ra}$. So we use the shorthand notations~$\Gamma_{\mathrm{sym},\mathbf{k}}\coloneqq \Gamma_{\raU{\mathrm{sym},\mathbf{k}}}$ and $\Gamma_{\mathrm{tr},\mathbf{k}}\coloneqq \Gamma_{\raU{\mathrm{tr},\mathbf{k}}}$.

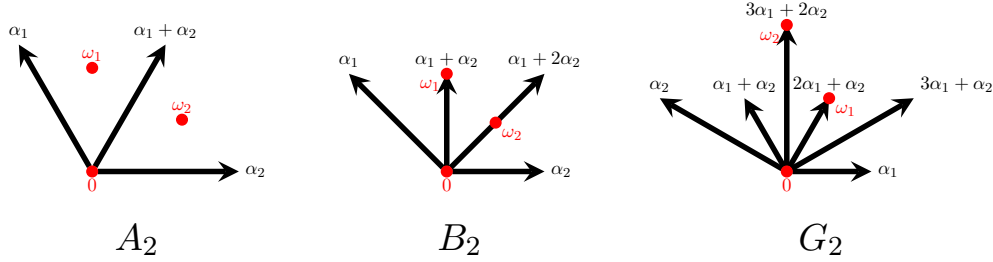
\begin{figure}
\def\scl{0.3}
\def\textscl{0.5}
\scalebox{1.3}{
\begin{tikzpicture}
  \node (A2) at (-3,0) {\begin{tikzpicture}
\node[scale=\scl,draw=white,circle] (ZERO) at (0,0) { };
\node[scale=\scl,draw=red,circle,fill=red] (Om2) at (30:{1.5/sqrt(2)}) { };
\node[scale=\scl,draw=red,circle,fill=red] (Om1) at (90:{1.5/sqrt(2)}) { };
\draw[->,>=stealth,ultra thick] (ZERO) to (120:1.5);
\draw[->,>=stealth,ultra thick] (ZERO) to (60:1.5);
\draw[->,>=stealth,ultra thick] (ZERO) to (0:1.5);
\node[anchor=south,scale=\textscl] at (120:1.5) {$\alpha_1$};
\node[anchor=west,scale=\textscl] at (0:1.5) {$\alpha_2$};
\node[anchor=south,scale=\textscl] at (60:1.5) {$\alpha_1+\alpha_2$};
\node[anchor=south,scale=\textscl] at (Om2) {\textcolor{red}{$\omega_2$}};
\node[anchor=south,scale=\textscl] at (Om1) {\textcolor{red}{$\omega_1$}};
\node at (90:2.2) {};
\node[scale=\scl,draw=red,circle,fill=red] (ZERO) at (0,0) { };
  \node[anchor=north,scale=\textscl,shift={(-0mm,-0.25mm)}] at (ZERO)  {\textcolor{red}{$0$}};
\end{tikzpicture}};
\node[anchor=north] at (A2.south) {$A_2$};

\node (B2) at (0.3,0) {\begin{tikzpicture}
\node[scale=\scl,draw=white,circle] (ZERO) at (0,0) { };
\draw[->,>=stealth,ultra thick] (ZERO) to (-1,1);
\draw[->,>=stealth,ultra thick] (ZERO) to (0,1);
\draw[->,>=stealth,ultra thick] (ZERO) to (1,1);
\draw[->,>=stealth,ultra thick] (ZERO) to (1,0);
\node[scale=\scl,draw=red,circle,fill=red] (Om2) at (0.5,0.5) { };
\node[scale=\scl,draw=red,circle,fill=red] (Om1) at (0,1) { };
\node[anchor=south,scale=\textscl] at (-1,1) {$\alpha_1$};
\node[anchor=west,scale=\textscl] at (1,0) {$\alpha_2$};
\node[anchor=south,scale=\textscl] at (0,1) {$\alpha_1+\alpha_2$};
\node[anchor=south,scale=\textscl] at (1,1) {$\alpha_1+2\alpha_2$};
\node[anchor=north west,scale=\textscl] at (Om2) {\textcolor{red}{$\omega_2$}};
\node[anchor=north east,scale=\textscl] at (Om1) {\textcolor{red}{$\omega_1$}};
\node at (90:2.2) {};
\node[scale=\scl,draw=red,circle,fill=red] (ZERO) at (0,0) { };
  \node[anchor=north,scale=\textscl,shift={(-0mm,-0.25mm)}] at (ZERO)  {\textcolor{red}{$0$}};
\end{tikzpicture}};
\node[anchor=north] at (B2.south) {$B_2$};

\node (G2) at (4,0) {\begin{tikzpicture}
\node[scale=\scl,draw=white,circle] (ZERO) at (0,0) { };
\draw[->,>=stealth,ultra thick] (ZERO) to (150:1.5);
\draw[->,>=stealth,ultra thick] (ZERO) to (120:{1.5/sqrt(3)});
\draw[->,>=stealth,ultra thick] (ZERO) to (90:1.5);
\draw[->,>=stealth,ultra thick] (ZERO) to (60:{1.5/sqrt(3)});
\draw[->,>=stealth,ultra thick] (ZERO) to (30:1.5);
\draw[->,>=stealth,ultra thick] (ZERO) to (0:{1.5/sqrt(3)});
\node[scale=\scl,draw=red,circle,fill=red] (Om2) at (90:1.5) { };
\node[scale=\scl,draw=red,circle,fill=red] (Om1) at (60:{1.5/sqrt(3)}) { };
\node[anchor=south,scale=\textscl] at (150:1.5) {$\alpha_2$};
\node[anchor=west,scale=\textscl] at (0:{1.5/sqrt(3)}) {$\alpha_1$};
\node[anchor=south,scale=\textscl] at (120:{1.5/sqrt(3)}) {$\alpha_1+\alpha_2$};
\node[anchor=south,scale=\textscl] at (60:{1.5/sqrt(3)}) {$2\alpha_1+\alpha_2$};
\node[anchor=south west,scale=\textscl] at (30:1.5) {$3\alpha_1+\alpha_2$};
\node[anchor=south,scale=\textscl] at (90:1.5) {$3\alpha_1+2\alpha_2$};
\node at (90:2.2) {};
\node[anchor=north east,scale=\textscl] at (Om2) {\textcolor{red}{$\omega_2$}};
\node[anchor=north west,scale=\textscl] at (Om1) {\textcolor{red}{$\omega_1$}};
\node[scale=\scl,draw=red,circle,fill=red] (ZERO) at (0,0) { };
  \node[anchor=north,scale=\textscl,shift={(-0mm,-0.25mm)}] at (ZERO)  {\textcolor{red}{$0$}};
\end{tikzpicture}};
\node[anchor=north] at (G2.south) {$G_2$};

\end{tikzpicture}
}
\caption{The positive roots of the rank~$2$ root systems $A_2$, $B_2$, and~$G_2$. The elements of $\Omega\cup\{0\}$ are shown in red.} \label{fig:rank2posroots}
\end{figure}

\begin{example} \label{ex:rank2graphs}
The irreducible rank~$2$ root systems are $A_2$, $B_2$ and $G_2$. The positive roots and fundamental weights for these root systems are depicted in Figure~\ref{fig:rank2posroots}. In Figures~\ref{fig:a2symtr},~\ref{fig:b2symtr}, and~\ref{fig:g2symtr} we depict the the truncated and symmetric interval-firing processes~$\Gamma_{\mathrm{tr},\mathbf{k}}$ and $\Gamma_{\mathrm{sym},\mathbf{k}}$ for $\mathbf{k}=0,1,2$ for these three root systems. Of course these graphs are infinite, so we depict the ``interesting part'' of the graphs near the origin (which is circled in black). The colors in these drawings correspond to classes of weights modulo the root lattice (hence there are three colors in the $A_2$ graphs, two in the $B_2$ graphs, and one in the $G_2$ graphs). Note that as $\mathbf{k}$ increases, the scale of the drawing is not maintained. Most, if not all, of the features of truncated and symmetric interval-firing that we care about are visible already in rank~$2$. Thus the reader is encouraged, while reading the rest of this paper, to return to these figures and understand how each of the results apply to these two dimensional examples.
\end{example}


\begin{figure}
\begin{center}
\ifthenelse{\boolean{bigimages}}{
\begin{tabular}{ >{\centering\arraybackslash} m{0.5cm} |  >{\centering\arraybackslash} m{6.25cm} |  >{\centering\arraybackslash} m{6.25cm} m{0cm}}
\large $\mathbf{k}$ & \large $\Gamma_{\mathrm{tr},\mathbf{k}}$ & \large $\Gamma_{\mathrm{sym},\mathbf{k}}$ & \\ [2ex] \hline & & \\
\end{tabular}}{{\bf IMAGE HERE}}
\end{center}

\caption[Interval-firing processes for $A_2$]{The graphs $\Gamma_{\mathrm{tr},\mathbf{k}}$ and $\Gamma_{\mathrm{sym},\mathbf{k}}$ for $\Phi=A_2$ and $\mathbf{k}=0,1,2$.}  \label{fig:a2symtr}
\end{figure}

\begin{figure}
\begin{center}
\ifthenelse{\boolean{bigimages}}{
\begin{tabular}{ >{\centering\arraybackslash} m{0.5cm} |  >{\centering\arraybackslash} m{6.25cm} |  >{\centering\arraybackslash} m{6.25cm} m{0cm}}
\large $\mathbf{k}$ & \large $\Gamma_{\mathrm{tr},\mathbf{k}}$ & \large $\Gamma_{\mathrm{sym},\mathbf{k}}$ & \\ [2ex] \hline & & \\
\end{tabular}}{{\bf IMAGE HERE}}
\end{center}

\caption[Interval-firing processes for $B_2$]{The graphs $\Gamma_{\mathrm{tr},\mathbf{k}}$ and $\Gamma_{\mathrm{sym},\mathbf{k}}$ for $\Phi=B_2$ and $\mathbf{k}=0,1,2$.}  \label{fig:b2symtr}
\end{figure}


\begin{figure}
\begin{center}
\ifthenelse{\boolean{bigimages}}{
\begin{tabular}{ >{\centering\arraybackslash} m{0.5cm} |  >{\centering\arraybackslash} m{6.25cm} |  >{\centering\arraybackslash} m{6.25cm} m{0cm}}
\large $\mathbf{k}$ & \large $\Gamma_{\mathrm{tr},\mathbf{k}}$ & \large  $\Gamma_{\mathrm{sym},\mathbf{k}}$ & \\ [2ex] \hline & & \\
\end{tabular}}{{\bf IMAGE HERE}}
\end{center}

\caption[Interval-firing processes for $G_2$]{The graphs $\Gamma_{\mathrm{tr},\mathbf{k}}$ and $\Gamma_{\mathrm{sym},\mathbf{k}}$ for $\Phi=G_2$ and $\mathbf{k}=0,1,2$.} \label{fig:g2symtr}
\end{figure} 

\begin{remark} \label{rem:chipfiringinterpretation}
Let us recall Propp's labeled chip-firing process (studied in~\cite{hopkins2017sorting}), which motivated our study of interval-firing processes. The states of labeled chip-firing are configurations of labeled chips on the infinite path graph $\mathbb{Z}$, such as:
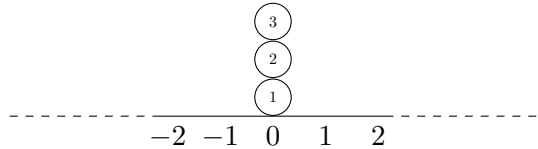
\begin{center}
\begin{tikzpicture}[scale=\chiptikzscl,block/.style={draw,circle, minimum width={width("11")+12pt},
font=\small,scale=\chipscl}]
    \foreach \x in {-2,-1,...,2} {%
      \node[anchor=north] (A\x) at (\x,0) {$\x$};
    }
    \draw (-2.2,0) -- (2.2,0);
    \draw[dashed] (-5,0) -- (-2.2,0);
     \draw[dashed] (5,0) -- (2.2,0);
    \foreach[count=\i] \a/\b in {0/1,0/2,0/3} {%
      \node[block] at (\a,{\b*\chipscl*\chipcoef-0.5*\chipscl*\chipcoef}) {$\i$};
      }
\end{tikzpicture}
\end{center}
If two chips with labels occupy the same position, we may \emph{fire} them, which sends the lesser-labeled chip one vertex to the right and the greater-labeled chip one vertex to the left. For instance, firing the chips~\chip{1} and~\chip{2} above leads to
\smallskip 
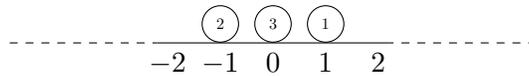
\begin{center}
\begin{tikzpicture}[scale=\chiptikzscl,block/.style={draw,circle, minimum width={width("11")+12pt},
font=\small,scale=\chipscl}]
    \foreach \x in {-2,-1,...,2} {%
      \node[anchor=north] (A\x) at (\x,0) {$\x$};
    }
    \draw (-2.2,0) -- (2.2,0);
    \draw[dashed] (-5,0) -- (-2.2,0);
     \draw[dashed] (5,0) -- (2.2,0);
    \foreach[count=\i] \a/\b in {1/1,-1/1,0/1} {%
      \node[block] at (\a,{\b*\chipscl*\chipcoef-0.5*\chipscl*\chipcoef}) {$\i$};
      }
\end{tikzpicture}
\end{center}
Firing the chips \chip{i} and \chip{j} with $i < j$ corresponds to $c \ra c + (e_i-e_j)$, where the integer vector $c \coloneqq  (c_1,\ldots,c_N) \in \mathbb{Z}^N$ is given by $c_i \coloneqq  \textrm{the position of the chip \chip{i}}$. In  this way central-firing (the subject of our sequel paper~\cite{galashin2017rootfiring2}) is the same as the labeled chip-firing process for $\Phi$ of Type A. Via this same correspondence between lattice vectors and configurations of chips, symmetric and truncated interval-firing in Type A can also be seen as ``labeled chip-firing processes'' that consist of the same chip-firing moves, which send chip~\chip{i} one vertex to the right and chip~\chip{j} one vertex to the left for any $i < j$, but where we allow these moves to be applied under different conditions: namely, when the position of chip~\chip{i} minus the position of chip~\chip{j} is either in the interval $[-k-1,k-1]$ (in the symmetric case) or in the interval $[-k,k-1]$ (in the truncated case). For example, consider the smallest non-trivial case of these interval-firing processes, which is symmetric interval-firing with $k=0$. This corresponds to the labeled chip-firing process that allows the transposition of the chips~\chip{i} and~\chip{j} with~$i < j$ when \chip{i} is one position to the left of \chip{j}. It is immediately apparent that this process is confluent; for instance, the configuration
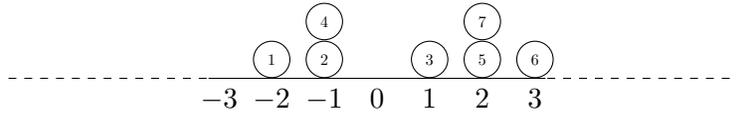
\begin{center}
\begin{tikzpicture}[scale=\chiptikzscl,block/.style={draw,circle, minimum width={width("11")+12pt},
font=\small,scale=\chipscl}]
    \foreach \x in {-3,-2,...,3} {%
      \node[anchor=north] (A\x) at (\x,0) {$\x$};
    }
    \draw (-3.2,0) -- (3.2,0);
    \draw[dashed] (-7,0) -- (-3.2,0);
     \draw[dashed] (7,0) -- (3.2,0);
    \foreach[count=\i] \a/\b in {-2/1,-1/1,1/1,-1/2,2/1,3/1,2/2} {%
      \node[block] at (\a,{\b*\chipscl*\chipcoef-0.5*\chipscl*\chipcoef}) {$\i$};
      }
\end{tikzpicture}
\end{center}
$\raU{\mathrm{sym},0}$-stabilizes to 
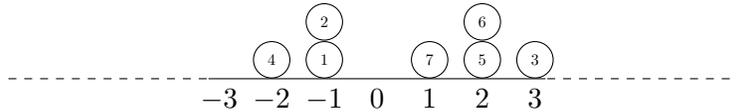
\vspace{-0.5cm}
\begin{center}
\begin{tikzpicture}[scale=\chiptikzscl,block/.style={draw,circle, minimum width={width("11")+12pt},
font=\small,scale=\chipscl}]
    \foreach \x in {-3,-2,...,3} {%
      \node[anchor=north] (A\x) at (\x,0) {$\x$};
    }
    \draw (-3.2,0) -- (3.2,0);
    \draw[dashed] (-7,0) -- (-3.2,0);
     \draw[dashed] (7,0) -- (3.2,0);
    \foreach[count=\i] \a/\b in {-1/1,-1/2,3/1,-2/1,2/1,2/2,1/1} {%
      \node[block] at (\a,{\b*\chipscl*\chipcoef-0.5*\chipscl*\chipcoef}) {$\i$};
      }
\end{tikzpicture}
\end{center}
In general the stabilization will weakly sort each collection of contiguous chips, while leaving the underlying unlabeled configuration of chips the same. The next smallest case to consider is truncated interval-firing with $k=1$. This corresponds to the labeled chip-firing process that allows both the transposition moves from the symmetric $k=0$ case, and the usual labeled chip-firing moves from the central-firing case. The reader can verify that for instance the configuration
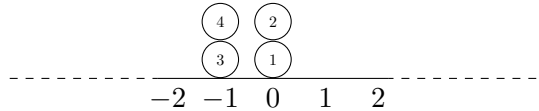
\begin{center}
\begin{tikzpicture}[scale=\chiptikzscl,block/.style={draw,circle, minimum width={width("11")+12pt},
font=\small,scale=\chipscl}]
    \foreach \x in {-2,-1,...,2} {%
      \node[anchor=north] (A\x) at (\x,0) {$\x$};
    }
    \draw (-2.2,0) -- (2.2,0);
    \draw[dashed] (-5,0) -- (-2.2,0);
     \draw[dashed] (5,0) -- (2.2,0);
    \foreach[count=\i] \a/\b in {0/1,0/2,-1/1,-1/2} {%
      \node[block] at (\a,{\b*\chipscl*\chipcoef-0.5*\chipscl*\chipcoef}) {$\i$};
      }
\end{tikzpicture}
\end{center}
$\raU{\mathrm{tr},1}$-stabilizes to 
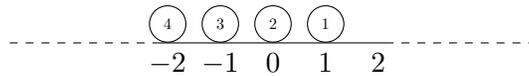
\vspace{-0.5cm}
\begin{center}
\begin{tikzpicture}[scale=\chiptikzscl,block/.style={draw,circle, minimum width={width("11")+12pt},
font=\small,scale=\chipscl}]
    \foreach \x in {-2,-1,...,2} {%
      \node[anchor=north] (A\x) at (\x,0) {$\x$};
    }
    \draw (-2.2,0) -- (2.2,0);
    \draw[dashed] (-5,0) -- (-2.2,0);
     \draw[dashed] (5,0) -- (2.2,0);
    \foreach[count=\i] \a/\b in {1/1,0/1,-1/1,-2/1} {%
      \node[block] at (\a,{\b*\chipscl*\chipcoef-0.5*\chipscl*\chipcoef}) {$\i$};
      }
\end{tikzpicture}
\end{center}
Here it is less obvious that confluence holds (although it is not too hard to prove this fact directly via a diamond lemma argument). The reader is now encouraged to experiment with this labeled chip-firing interpretation of symmetric and truncated interval-firing for higher values of $k$. Note that increasing $k$ allows for the firing of chips~\chip{i} and~\chip{j} when they are further apart.
\end{remark}

In our further treatment of the interval-firing processes we will focus on the geometric picture (on display in Example~\ref{ex:rank2graphs}) and not the chip-firing picture (discussed in Remark~\ref{rem:chipfiringinterpretation}).

To close out this section, let us demonstrate that the interval-firing processes are always terminating. This is straightforward because the collection~$\Phi^+$ of vectors we are adding is acyclic.

\begin{prop} \label{prop:interval_firing_terminating}
For~$\mathbf{k}\in\mathbb{N}[\Phi]^W$, the relations $\raU{\mathrm{sym},\mathbf{k}}$ and $\raU{\mathrm{tr},\mathbf{k}}$ are terminating.
\end{prop}
\begin{proof}
It is enough to show this for $\raU{\mathrm{sym},\mathbf{k}}$, which has more firing moves than~$\raU{\mathrm{tr},\mathbf{k}}$. For~$\lambda \in P$ define $\varphi(\lambda) \coloneqq  \<\rho_{\mathbf{k}+1}-\lambda,  \rho_{\mathbf{k}+1}-\lambda\>$; in other words, $\varphi(\lambda)$ is the length of the vector $\rho_{\mathbf{k}+1}-\lambda$. Suppose $\lambda \raU{\mathrm{sym},\mathbf{k}} \lambda + \alpha$ for $\alpha \in \Phi^+$. Then,
\begin{align*}
\varphi(\lambda) - \varphi(\lambda+\alpha) &= \< \rho_{\mathbf{k}+1}-\lambda, \rho_{\mathbf{k}+1}-\lambda\>- \< \rho_{\mathbf{k}+1}-(\lambda+\alpha), \rho_{\mathbf{k}+1}-(\lambda+\alpha)\>\\
&=2\< \rho_{\mathbf{k}+1},\alpha\>-2\<\lambda,\alpha\>-\<\alpha,\alpha\>\\
&\geq \<\alpha,\alpha\>(\mathbf{k}(\alpha)+1-\mathbf{k}(\alpha)+1-1) =\<\alpha,\alpha\>,
\end{align*}
where we use the facts that $\<\lambda,\alpha\>\leq \frac{\<\alpha,\alpha\>}{2}(\mathbf{k}(\alpha)-1)$ since $\lambda \raU{\mathrm{sym},\mathbf{k}} \lambda + \alpha$, and that $\<\rho_{\mathbf{k}+1},\alpha\>\geq\frac{\<\alpha,\alpha\>}{2}(\mathbf{k}(\alpha)+1)$ because $\alpha$ is $W$-conjugate to at least one simple root appearing with nonzero coefficient in its expansion in terms of simple roots. So each firing move causes the quantity  $\varphi(\lambda)$ to decrease by at least some fixed nonzero amount. But $\varphi(\lambda)\geq 0$ because it is the length of a vector. Thus indeed $\raU{\mathrm{sym},\mathbf{k}}$ is terminating.
\end{proof}

\section{Symmetries of interval-firing processes} \label{sec:symmetry}

In this section we study the symmetries of the two interval-firing processes. Since the set of positive roots $\Phi^+$ is an ``oriented'' set of vectors, we do not expect the directed graphs $\Gamma_{\mathrm{sym},\mathbf{k}}$ and~$\Gamma_{\mathrm{tr},\mathbf{k}}$  to have many symmetries, and certainly none coming from the Weyl group. But if we consider instead the undirected graphs $\Gamma^{\mathrm{un}}_{\mathrm{sym},\mathbf{k}}$ and~$\Gamma^{\mathrm{un}}_{\mathrm{tr},\mathbf{k}}$ (corresponding to the symmetric relations~$\baU{\mathrm{sym},\mathbf{k}}$ and~$\baU{\mathrm{tr},\mathbf{k}}$), we will see that both of these do in fact have symmetries coming from the Weyl group.

For the symmetric interval-firing process, the graph $\Gamma^{\mathrm{un}}_{\mathrm{sym},\mathbf{k}}$ is invariant under the action of the whole Weyl group~$W$. This explains the name ``symmetric'' for the process: it has the biggest possible group of symmetries. As for the truncated process, in order to understand its symmetries we need to introduce a certain subgroup of the Weyl group~$C\subseteq W$. In fact this $C$ is an abelian group and satisfies $C\simeq P/Q$. In our definition of $C$ we follow Lam and Postnikov~\cite{lam2012alcoved}\footnote{Lam and Postnikov worked in a completely dual setting to ours: that is, they described a copy of the coweight lattice modulo the coroot lattice inside of $W$; hence, they used $\theta$ instead of $\widehat{\theta}$, etc.}. The \emph{Coxeter number} of $\Phi$, another fundamental invariant of the root system, is $h\coloneqq  \<\rho,\widehat{\theta}^\vee\>+1$. (The Coxeter number is also equal to $h = 1+\sum_{i=1}^{n}a_i$ where $\theta=\sum_{i=1}^{n}a_i\alpha_i$). Lam and Postnikov~\cite[\S5]{lam2012alcoved} defined the subgroup $C\coloneqq \{w\in W\colon \rho-w(\rho) \in hP\}$ of the Weyl group  and explained (using the affine Weyl group, which we will not discuss here) that $C$ is naturally isomorphic to $P/Q$: the isomorphism is explicitly given by~$w \mapsto \omega \in \Omega^0_m$ if and only if~$\rho-w(\rho) = h\omega$. (Since $\rho-w(\rho) \in Q$ for any $w\in W$, a consequence of this description of the isomorphism is that $h\cdot(P/Q) =\{0\}$.) As they mention, this subgroup was also studied before by Verma~\cite{verma1975affine}, but in spite of its significance it does not seem to have any name other than $C$ in the root system literature. Lam and Postnikov gave another characterization~\cite[Proposition 6.4]{lam2012alcoved} of~$C$ that will be useful for us: $C=\{w\in W\colon w(\{\alpha^\vee_0,\alpha^\vee_1,\ldots,\alpha^\vee_n\}) = \{\alpha^\vee_0,\alpha^\vee_1,\ldots,\alpha^\vee_n\}\}$, where we use the suggestive notation $\alpha^\vee_0 \coloneqq  -\widehat{\theta}^\vee$.

\begin{thm} \label{thm:symmetry}
Let $\mathbf{k}\in\mathbb{N}[\Phi]^W$. Set $\Gamma \coloneqq  \Gamma^{\mathrm{un}}_{\mathrm{sym},\mathbf{k}}$ or $\Gamma \coloneqq  \Gamma^{\mathrm{un}}_{\mathrm{tr},\mathbf{k}}$. Then,
\begin{itemize}
\item if $\Gamma = \Gamma^{\mathrm{un}}_{\mathrm{sym},\mathbf{k}}$, the linear map $v\mapsto w(v)$ is an automorphism of $\Gamma$ for all $w \in W$;
\item if $\Gamma = \Gamma^{\mathrm{un}}_{\mathrm{tr},\mathbf{k}}$, the affine map $v \mapsto w(v-\frac{1}{h}\rho)+\frac{1}{h}\rho$ is an automorphism of $\Gamma$ for all~$w \in C\subseteq W$.
\end{itemize}
\end{thm}
\begin{proof}
If $\Gamma = \Gamma^{\mathrm{un}}_{\mathrm{sym},\mathbf{k}}$ set $c\coloneqq 0$, and if $\Gamma = \Gamma^{\mathrm{un}}_{\mathrm{tr},\mathbf{k}}$ set $c\coloneqq 1$. Consider the hyperplane arrangement $\mathcal{H} \coloneqq  \left\{H_{\alpha^\vee,\frac{c}{2}}\colon \alpha \in \Phi^{+}\right\}$ with hyperplanes $H_{\alpha^\vee,\frac{c}{2}} \coloneqq 
  \left\{v\in V\colon \<v,\alpha^\vee\> = \frac{c}{2}\right\}$.
  
First we claim that if for $w\in W$ and $u\in V$ the affine map $\varphi\colon v\mapsto w(v-u)+u$ is an automorphism of $\mathcal{H}$ which maps $P$ to $P$, then it is an automorphism of $\Gamma$ (by an automorphism of the hyperplane arrangement, we mean an invertible affine map $\varphi$ such that $\varphi$ permutes the hyperplanes in $\mathcal{H}$). Indeed, observe that there is an edge in $\Gamma$ between $\lambda$ and $\mu$ if and only if there is some $\alpha \in \Phi^{+}$ such that $\mu = \lambda +\alpha$ and $\mathrm{max}(\{|\<\mu,\alpha^\vee\>-\frac{c}{2}|,|\<\lambda,\alpha^\vee\>-\frac{c}{2}|\}) \leq \mathbf{k}(\alpha)+1-\frac{c}{2}$. So suppose there is an edge between~$\lambda$ and~$\mu$ in the $\alpha$ direction. Any $\varphi$ of this form will satisfy~$\varphi(\mu)-\varphi(\lambda)=w(\alpha)$ and $\varphi(H_{\alpha^\vee,\frac{c}{2}}) = H_{\pm w(\alpha)^\vee,\frac{c}{2}}$ (where the sign~$\pm$ is chosen so that $\pm w(\alpha) \in \Phi^{+}$). Moreover, since all Weyl group elements are orthogonal, and, in particular, preserve distances, the distance from $\mu$ to~$H_{\alpha^\vee,\frac{c}{2}}$ will be the same as the distance from $\varphi(\mu)$ to $H_{\pm w(\alpha)^\vee,\frac{c}{2}}$, and ditto for $\lambda$. But $|\<\mu,\alpha^\vee\>-\frac{c}{2}|$ is precisely the distance from $\mu$ to~$H_{\alpha^\vee,\frac{c}{2}}$, and ditto for $\lambda$. Hence indeed we will get that $\varphi(\mu) = \varphi(\lambda)+w(\alpha)$ and that
\begin{align*}
\mathrm{max}\left(\left\{\left|\<\varphi(\mu),(\pm w(\alpha))^\vee\>-\frac{c}{2}\right|,\left|\<\varphi(\lambda),(\pm w(\alpha))^\vee\>-\frac{c}{2}\right|\right\}\right) \\\leq \mathbf{k}(\alpha)+1-\frac{c}{2}
=  \mathbf{k}(\pm w(\alpha))+1-\frac{c}{2},
\end{align*}
which means there is an edge in $\Gamma$ between $\varphi(\lambda)$ and $\varphi(\mu)$ in the $\pm w(\alpha)$ direction. To see that conversely if there is an edge between $\varphi(\lambda)$ and $\varphi(\mu)$ in $\Gamma$, there is one between~$\lambda$ and $\mu$, use that $\varphi$ is invertible and $\varphi^{-1}$ is of the same form.

In the case $c=0$, the hyperplane arrangement $\mathcal{H}$ is just the Coxeter arrangement of~$\Phi$ and it is easy to see that every $w\in W$ is an automorphism of~$\mathcal{H}$.

Now consider the case $c=1$, in which case $\mathcal{H}$ is (a scaled version of) the \emph{$\Phi^\vee$-Linial arrangement}; see for instance~\cite{postnikov2000deformations} and~\cite{athanasiadis2000deformations}. We claim that $\varphi\colon v \mapsto w(v-\frac{c}{h}\rho)+\frac{c}{h}\rho$ is an automorphism of $\mathcal{H}$ for all $w \in C$. So suppose $x \in H_{\alpha^\vee,\frac{c}{2}}$; we want to show that $\varphi(x) \in H_{\pm w(\alpha)^\vee,\frac{c}{2}}$ where the sign $\pm$ is chosen so that $\pm w(\alpha)$ is positive. (The reverse implication will then follow from consideration of $\varphi^{-1}(v) = w^{-1}(v-\frac{c}{h}\rho)+\frac{c}{h}\rho$.) We have
\begin{equation} \label{eqn:onhyperplane}
\<\varphi(x),w(\alpha)^\vee\> = \frac{c}{2}-\left\<\frac{c}{h}\rho,\alpha^\vee\right\>+\left\<\frac{c}{h}\rho,w(\alpha)^\vee\right\>.
\end{equation}
Write $\alpha^\vee = \sum_{i=1}^{n}a_i\alpha^\vee_i$, with the convention $a_0\coloneqq 0$. By a result of Lam-Postnikov mentioned above, there is a permutation $\pi\colon \{0,1,\ldots,n\}\to \{0,1,\ldots,n\}$ such that $w(\alpha^\vee_i) = \alpha^\vee_{\pi(i)}$ (with the aforementioned convention $\alpha^\vee_0 \coloneqq  -\widehat{\theta}^\vee$ where $\widehat{\theta}^\vee$ is the highest root of $\Phi^\vee$). Thus, $w(\alpha)^\vee=\sum_{i=1}^{n} a_i\alpha^\vee_{\pi(i)}$. 

We will consider two cases. First suppose that $a_{\pi^{-1}(0)} = 0$. Then $w(\alpha)^\vee$ is clearly a positive root, so $\pm = +$; moreover, we have $\<\frac{c}{h}\rho,\alpha^\vee\>=\<\frac{c}{h}\rho,w(\alpha)^\vee\>=\frac{c}{h}\cdot\sum_{i=1}^{n}a_i$. So from~\eqref{eqn:onhyperplane} we get that $\<\varphi(x),w(\alpha)^\vee\>=\frac{c}{2}$, that is, $\varphi(x) \in H_{\pm w(\alpha)^\vee,\frac{c}{2}}$, as desired.

Now suppose that $a_{\pi^{-1}(0)} \neq 0$. We claim that this forces $a_{\pi^{-1}(0)} = 1$: indeed, otherwise the height of $w(\alpha)^\vee$ would be strictly less than $-(h-1)$, which is impossible because $-\widehat{\theta}^{\vee}$ has height $-(h-1)$ and is the root in $\Phi^\vee$ of smallest height. So indeed we have $a_{\pi^{-1}(0)} = 1$. Note also that in this case the height of $w(\alpha)^\vee$ is a negative root and hence $w(\alpha)^\vee$ is negative, so $\pm= -$. Then we compute
\[-\left\<\frac{c}{h}\rho,\alpha^\vee\right\>+\left\<\frac{c}{h}\rho,w(\alpha)^\vee\right\>=
-\left\<\frac{c}{h}\rho,\alpha_{\pi^{-1}(0)}^\vee\right\>+\left\<\frac{c}{h}\rho,\alpha^\vee_0\right\> = -\frac{c}{h}-\left(\frac{c}{h}(h-1)\right) = -c.\]
Thus from~\eqref{eqn:onhyperplane} we get that $\<\varphi(x),-w(\alpha)^\vee\> = -\frac{c}{2}+c = \frac{c}{2}$, that is, $\varphi(x) \in H_{\pm w(\alpha)^\vee,\frac{c}{2}}$, as desired.

Finally, the description of $C$ given above says that $\varphi(0) = w(0-\frac{c}{h}\rho)+\frac{c}{h}\rho = c\omega$ for some $\omega \in \Omega^0_m$. Hence indeed $\varphi$ maps $P$ to $P$.
\end{proof}

\section{Sinks of symmetric interval-firing and the map~\texorpdfstring{$\eta$}{eta}}

Recall that our overall strategy for proving confluence of the interval-firing processes is to show that they get ``trapped'' inside certain permutohedra, and then to analyze where these processes must terminate. In order to carry out this strategy, we need to understand what are the possible final points we terminate at, i.e., what are the stable points of these processes.

In this section we describe the $\raU{\mathrm{sym},\mathbf{k}}$-stable points, i.e., the sinks of $\Gamma_{\mathrm{sym},\mathbf{k}}$. We will show in particular that there is a way to consistently label the sinks of $\Gamma_{\mathrm{sym},\mathbf{k}}$ across all values of~$\mathbf{k}$. 

In order to define this labeling we need to review some basic facts about parabolic subgroups and parabolic cosets. Recall that the Weyl group $W$ is generated by the \emph{simple reflections}~$s_i \coloneqq  s_{\alpha_i}$ for $i=1,\ldots,n$. For any $w \in W$ we use $\ell(w)$ to denote the \emph{length} of~$w$, which is the length of the shortest representation of $w$ as a product of simple reflections. An \emph{inversion} of $w$ is a positive root $\alpha \in \Phi^{+}$ for which $w(\alpha)$ is negative. The length $\ell(w)$ is equal to the number of inversions of $w$. The identity is the only Weyl group element of length zero. The simple reflections are the only Weyl group elements of length one: $s_i$ sends $\alpha_i$ to~$-\alpha_i$ and permutes $\Phi^{+}\setminus\{\alpha_i\}$. A \emph{(right) descent} of $w \in W$ is a simple reflection~$s_i$ such that $\ell(ws_i) <\ell(w)$. The reflection $s_i$ is a descent of $w$ if and only if $\alpha_i$ is an inversion of $w$.

Recall that for $I\subseteq[n]$ we use $W_I$ to denote the corresponding \emph{parabolic subgroup} of $W$, that is, the subgroup of~$GL(V)$ generated by simple reflections $s_i$ for $i \in I$. Note that $W_I$ is (isomorphic to) the Weyl group of $\Phi_I$. For $\lambda \in P$ we define the parabolic permutohedron $\Pi_I(\lambda) \coloneqq \mathrm{ConvexHull}\, W_I(\lambda)$ and $\Pi^Q_I(\lambda)\coloneqq \Pi_I(\lambda)\cap(Q+\lambda)$. An important property of parabolic subgroups is the existence of distinguished coset representatives: each (left) coset $wW_I$ in $W$ contains a unique element of minimal length. We use $W^I$ for the set of minimal length coset representatives of~$W_I$. There is even an explicit description: $W^I \coloneqq  \{w\in W\colon \textrm{$s_i$ is not a descent of $w$ for all $i \in I$}\}$ (see for instance~\cite[\S2.4]{bjorner2005coxeter}).

Recall that for any $\lambda \in P$ we use $\lambda_{\mathrm{dom}}$ to denote the dominant element of $W(\lambda)$. For a dominant weight $\lambda = \sum_{i=1}^n c_i \omega_i \in P_{\geq 0}$, we define $I^{0}_{\lambda} \coloneqq  \{i \in [n] \colon c_i = 0\}$. And then for any weight $\lambda \in P$ we define $I^{0}_{\lambda} \coloneqq I^{0}_{\lambda_{\mathrm{dom}}}$.

\begin{prop} \label{prop:stabilizer}
For $\lambda \in P_{\geq 0}$, the stabilizer of $\lambda$ in $W$ is $W_{I^{0}_{\lambda}}$.
\end{prop}
\begin{proof}
This (straightforward proposition) is~\cite[Lemma 10.2B]{humphreys1972lie}.
\end{proof}

\begin{cor} \label{cor:cosets}
For any $\lambda \in P$, $\{w \in W\colon w^{-1}(\lambda) \in P_{\geq 0}\}$ is a coset of~$W_{I^{0}_{\lambda}}$.
\end{cor}
\begin{proof}
First let us show that if $w^{-1}(\lambda)$ is dominant then $(ww')^{-1}(\lambda)$ is dominant for any $w' \in W_{I^{0}_{\lambda}}$. This is clear: $(ww')^{-1}(\lambda) = (w')^{-1}(w^{-1}(\lambda)) =  (w')^{-1}(\lambda_{\mathrm{dom}}) = \lambda_{\mathrm{dom}}$ since $w'$ is in the stabilizer of $\lambda_{\mathrm{dom}}$ by Proposition~\ref{prop:stabilizer}. Next let us show that if~$w^{-1}(\lambda)$ is dominant and $(w')^{-1}(\lambda)$ is dominant then $w' = ww''$ for some $w'' \in W_{I^{0}_{\lambda}}$. This is also clear: $w^{-1}(w'(\lambda_{\mathrm{dom}})) = w^{-1}(\lambda) = \lambda_{\mathrm{dom}}$, so $w^{-1}w'$ is in the stabilizer of~$\lambda_{\mathrm{dom}}$, that is, $w^{-1}w'= w''$ for some $w'' \in W_{I^{0}_{\lambda}}$ thanks to Proposition~\ref{prop:stabilizer}, as claimed.
\end{proof}

In light of Corollary~\ref{cor:cosets}, for $\lambda \in P$ we define $w_{\lambda}$ to be the minimal length element of~$\{w \in W\colon w^{-1}(\lambda) \in P_{\geq 0}\}$. Hence, for $\lambda \in P_{\geq 0}$ we have (by the Orbit-Stabilizer Theorem) that $W^{I^{0}_{\lambda}} = \{w_\mu\colon \mu \in W(\lambda)\}$ and $w_{\mu}\neq w_{\mu'}$ for $\mu\neq \mu' \in W(\lambda)$. Another way to think about $w_{\lambda}$: $\lambda$ may belong to the closure of many chambers, but there will be a unique chamber $wC_0$ with $w$ of minimal length such that $\lambda$ belongs to the closure of~$wC_0$ and this is when $w=w_{\lambda}$. Then for $\mathbf{k} \in \mathbb{N}[\Phi]^W$, we define the map $\eta_{\mathbf{k}}\colon P \to P$ by setting~$\eta_{\mathbf{k}}(\lambda) \coloneqq  \lambda + w_{\lambda}(\rho_{\mathbf{k}})$ for all $\lambda \in P$ (where, as above, we have~$\rho_{\mathbf{k}} \coloneqq  \sum_{i=1}^{n} \mathbf{k}(\alpha_i) \omega_i$).

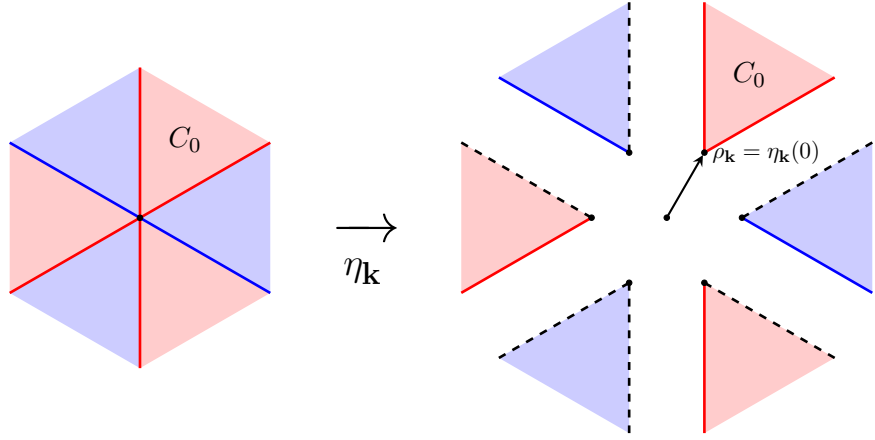
\begin{figure}
\begin{center}
  \def\originscale{0.2}
  \def\lw{1}
  \def\tikzscale{1}
\begin{tikzpicture}
\node at (0,0) {\begin{tikzpicture}[scale=-\tikzscale]
\fill[color=red!20]  (0,0) -- (30:-2) -- (90:-2) -- cycle;
\fill[color=blue!20]  (0,0) -- (150:-2) -- (90:-2) -- cycle;
\fill[color=blue!20]  (0,0) -- (-30:-2) -- (30:-2) -- cycle;
\fill[color=red!20]  (0,0) -- (-90:-2) -- (-30:-2) -- cycle;
\fill[color=red!20]  (0,0) -- (210:-2) -- (150:-2) -- cycle;
\fill[color=blue!20]  (0,0) -- (210:-2) -- (-90:-2) -- cycle;

\draw[color=blue,line width=\lw] (0,0) -- (150:2);
\draw[color=blue,line width=\lw] (0,0) -- (-30:2);
\draw[color=red,line width=\lw] (0,0) -- (30:2);
\draw[color=red,line width=\lw] (0,0) -- (90:2);
\draw[color=red,line width=\lw] (0,0) -- (-90:2);
\draw[color=red,line width=\lw] (0,0) -- (210:2);

\node[draw,circle,fill=black,scale=\originscale] at (0,0) { };
\node at (60:-1.2) {$C_0$};
\end{tikzpicture}};

\node at ({\tikzscale*3},-0.5) {\huge $\xrightarrow[\eta_{\mathbf{k}}]{}$};

\node at ({7*\tikzscale},0) {\begin{tikzpicture}[scale=\tikzscale]
\begin{scope}[shift={(60:1)}]
\fill[color=red!20]  (0,0) -- (30:2) -- (90:2) -- cycle;
\draw[color=red,line width=\lw] (0,0) -- (30:2);
\draw[color=red,line width=\lw] (0,0) -- (90:2);
\node[draw,circle,fill=black,scale=\originscale] (rho) at (0,0) { };
\node at (60:1.2) {$C_0$};
\end{scope}

\begin{scope}[shift={(120:1)}]
\fill[color=blue!20]  (0,0) -- (150:2) -- (90:2) -- cycle;
\draw[color=blue,line width=\lw] (0,0) -- (150:2);
\draw[color=black,line width=\lw,dashed] (0,0) -- (90:2);
\node[draw,circle,fill=black,scale=\originscale] at (0,0) { };
\end{scope}

\begin{scope}[shift={(0:1)}]
\fill[color=blue!20]  (0,0) -- (-30:2) -- (30:2) -- cycle;
\draw[color=blue,line width=\lw] (0,0) -- (-30:2);
\draw[color=black,line width=\lw,dashed] (0,0) -- (30:2);
\node[draw,circle,fill=black,scale=\originscale] at (0,0) { };
\end{scope}

\begin{scope}[shift={(-60:1)}]
\fill[color=red!20]  (0,0) -- (-90:2) -- (-30:2) -- cycle;
\draw[color=red,line width=\lw] (0,0) -- (-90:2);
\draw[color=black,line width=\lw,dashed] (0,0) -- (-30:2);
\node[draw,circle,fill=black,scale=\originscale] at (0,0) { };
\end{scope}

\begin{scope}[shift={(180:1)}]
\fill[color=red!20]  (0,0) -- (210:2) -- (150:2) -- cycle;
\draw[color=red,line width=\lw] (0,0) -- (210:2);
\draw[color=black,line width=\lw,dashed] (0,0) -- (150:2);
\node[draw,circle,fill=black,scale=\originscale] at (0,0) { };
\end{scope}

\begin{scope}[shift={(240:1)}]
\fill[color=blue!20]  (0,0) -- (210:2) -- (-90:2) -- cycle;
\draw[color=black,line width=\lw,dashed] (0,0) -- (210:2);
\draw[color=black,line width=\lw,dashed] (0,0) -- (-90:2);
\node[draw,circle,fill=black,scale=\originscale] at (0,0) { };
\end{scope}

\draw[->,>=stealth,thick] (0,0) to (rho);
\node[draw,circle,fill=black,scale=\originscale] at (0,0) { };
\node[anchor=west,scale=0.8] at (rho) {$\rho_{\mathbf{k}}=\eta_{\mathbf{k}}(0)$};
\end{tikzpicture}};
\end{tikzpicture}
\end{center}
\caption[The map~$\eta$]{A graphical depiction of the piecewise-linear map $\eta_{\mathbf{k}}$.} \label{fig:eta}
\end{figure}

This map $\eta_{\mathbf{k}}$ will be of crucial importance for us in our investigation of both the symmetric and truncated interval-firing processes and the relationship between these two processes. Figure~\ref{fig:eta} gives a graphical depiction of $\eta_{\mathbf{k}}$: as we can see, this map ``dilates'' space by translating the chambers radially outwards; a point not inside any chamber travels in the same direction as the chamber closest to the fundamental chamber among those chambers whose closure the point lies in. The following proposition lists some very basic properties of $\eta_{\mathbf{k}}$.

\begin{prop} \label{prop:etafacts} $ $
\begin{itemize}
\item For any $\mathbf{k},\mathbf{m} \in \mathbb{N}[\Phi]^W$, we have $\eta_{\mathbf{k}+\mathbf{m}} = \eta_{\mathbf{m}}(\eta_{\mathbf{k}})$.
\item For any $\mathbf{k} \in \mathbb{N}[\Phi]^W$, the map $\eta_{\mathbf{k}}\colon P \to P$ is injective.
\end{itemize}
\end{prop}
\begin{proof}
For the first bullet point: let $\lambda \in P$. Set $\lambda' \coloneqq  \eta_{\mathbf{k}}(\lambda) = \lambda + w_{\lambda}(\rho_{\mathbf{k}})$. Observe that $\lambda'_{\mathrm{dom}} = w_{\lambda}^{-1}(\eta_{\mathbf{k}}(\lambda)) = \lambda_{\mathrm{dom}} + \rho_{\mathbf{k}}$. Hence, $I^{0}_{\lambda'} \subseteq I^{0}_{\lambda}$. This means the cosets of~$W_{I^{0}_{\lambda}}$ are unions of cosets of~$W_{I^{0}_{\lambda'}}$. But we just saw that $w_{\lambda} \in w_{\lambda'}W_{I^{0}_{\lambda'}}$, because~$w_{\lambda}^{-1}(\lambda')$ is dominant. So $w_{\lambda}$ must be the minimal length element of $w_{\lambda'}W_{I^{0}_{\lambda'}}$ (since it is the minimal length element of a superset of $w_{\lambda'}W_{I^{0}_{\lambda'}}$). Hence $w_{\lambda'} = w_{\lambda}$. This means that $\eta_{\mathbf{m}}(\eta_{\mathbf{k}}(\lambda)) = \lambda + w_{\lambda}(\rho_{\mathbf{k}}) + w_{\lambda}(\rho_{\mathbf{m}}) = \lambda + w_{\lambda}(\rho_{\mathbf{k+m}}) = \eta_{\mathbf{k}+\mathbf{m}}(\lambda)$ and thus the claim is proved.

For the second bullet point: suppose $\lambda,\mu \in P$ with $\eta_{\mathbf{k}}(\lambda) = \eta_{\mathbf{k}}(\mu)$. First of all, since~$\eta_{\mathbf{k}}(\lambda)_{\mathrm{dom}} = \lambda_{\mathrm{dom}} + \rho_{\mathbf{k}}$ and similarly for $\mu$, we have $\lambda_{\mathrm{dom}}=\mu_{\mathrm{dom}}$. Let~$m \gg 0 \in \mathbb{Z}$ be some very large constant. From the first bullet point we know~$\eta_{\mathbf{k}+m}(\lambda) = \eta_{\mathbf{k}+m}(\mu)$ and hence $\lambda + w_{\lambda}(\rho_{\mathbf{k}+m}) = \mu + w_{\mu}(\rho_{\mathbf{k}+m})$. But $\rho_{\mathbf{k}+m}$ is inside the fundamental chamber~$C_0$, and hence $w(\rho_{\mathbf{k}+m}) = w'(\rho_{\mathbf{k}+m})$ if and only if $w=w'$. Moreover, by taking $m$ large enough we can guarantee that $w(\rho_{\mathbf{k}+m})$ and $w'(\rho_{\mathbf{k}+m})$ are very far away from one another for $w \neq w'$. Hence $\lambda + w_{\lambda}(\rho_{\mathbf{k}+m}) = \mu + w_{\mu}(\rho_{\mathbf{k}+m})$ in fact forces $w_{\lambda} = w_{\mu}$. But~$w_{\lambda} = w_{\mu}$ together with~$\lambda_{\mathrm{dom}} = \mu_{\mathrm{dom}}$ means $\lambda = \mu$ and thus the claim is proved.
\end{proof}

In light of Proposition~\ref{prop:etafacts} it makes sense to set~$\eta \coloneqq  \eta_1$ so that $\eta_k = \eta^k$. Now we proceed to explain how $\eta_{\mathbf{k}}$ labels the sinks of $\Gamma_{\mathrm{sym},\mathbf{k}}$. 

For a dominant weight~$\lambda = \sum_{i=1}^n c_i \omega_i \in P_{\geq 0}$, define $I^{0,1}_{\lambda} \coloneqq  \{i \in [n] \colon c_i \in \{0,1\}\}$. And for any weight $\lambda \in P$ define $I^{0,1}_{\lambda} \coloneqq I^{0,1}_{\lambda_{\mathrm{dom}}}$.

\begin{prop} \label{prop:posparaboliccosets}
Let $\lambda \in P$ with $\<\lambda,\alpha^\vee\> \neq -1$ for all~$\alpha \in \Phi^{+}$. Then $w_{\lambda} (\Phi^{+}_{I^{0,1}_{\lambda}})$ is a subset of positive roots.
\end{prop}
\begin{proof}
It suffices to show that $w_{\lambda}(\alpha_i)$ is positive for all $i \in I^{0,1}_{\lambda}$. Suppose that $w_{\lambda}(\alpha_i)$ is negative for some $i \in I^{0,1}_{\lambda}$, i.e., $s_i$ is a descent of $w_{\lambda}$. Note~$\<\lambda_{\mathrm{dom}},\alpha_i^\vee\> \in \{0,1\}$. If $\<\lambda_{\mathrm{dom}},\alpha_i^\vee\> = 1$, then $\<\lambda_{\mathrm{dom}},-\alpha_i^\vee\> = -1$ so~$\<\lambda,-w_{\lambda}(\alpha_i)^\vee\> = -1$, which contradicts that  $\<\lambda,\alpha^\vee\> \neq -1$ for all $\alpha \in \Phi^{+}$. But since $w_{\lambda}$ is the minimal length representative of $w_{\lambda}W_{I^{0}_{\lambda}}$, it cannot have any descents $s_j$ with $j\in I^{0}_{\lambda}$. Hence we cannot have that~$\<\lambda_{\mathrm{dom}},\alpha_i^\vee\> = 0$ either. Thus it must be that $w_{\lambda}(\alpha_i)$ is positive for all $i \in I^{0,1}_{\lambda}$.
\end{proof}

\begin{prop} \label{prop:sinksminlen}
For a dominant weight $\mu \in P_{\geq 0}$, we have that
\[ W^{I^{0,1}_{\mu}} = \{w_{\lambda}\colon \lambda \in P, \lambda_{\mathrm{dom}}=\mu, \<\lambda,\alpha^\vee\> \neq -1 \textrm{ for all $\alpha \in \Phi^{+}$}\}.\]
\end{prop}
\begin{proof}
Let $\lambda \in P$ with $\lambda_{\mathrm{dom}}=\mu$ and first suppose that $\<\lambda,\alpha^\vee\>=-1$ for some $\alpha \in \Phi^{+}$. Then we have~$\<w^{-1}_{\lambda}(\lambda),w^{-1}_{\lambda}(\alpha)^\vee\> = -1$. But since $w^{-1}_{\lambda}(\lambda) = \lambda_{\mathrm{dom}}$ is dominant, this means~$w^{-1}_{\lambda}(\alpha)$ is a negative root; moreover, the only way $\<\lambda_{\mathrm{dom}},w^{-1}_{\lambda}(\alpha)^\vee\> = -1$ is possible is if all the simple coroots $\alpha_i^\vee$ appearing in the expansion of $-w^{-1}_{\lambda}(\alpha)^\vee$ have $i \in I^{0,1}_{\lambda}$. This implies that~$w_{\lambda}(\alpha_i)$ is negative for some $i \in I^{0,1}_{\lambda}$. But then $s_i$ would be a descent of $w_{\lambda}$, and hence $w_{\lambda}$ cannot be the minimal length element of $w_{\lambda}W_{I^{0,1}_{\lambda}}$.

If $\lambda \in P$ with $\lambda_{\mathrm{dom}}=\mu$ satisfies $\<\lambda,\alpha^\vee\>\neq-1$ for all $\alpha \in \Phi^{+}$, then we have seen in Proposition~\ref{prop:posparaboliccosets} that $w_{\lambda}$ has no descents $s_i$ with $i \in I^{0,1}_{\mu}$ and hence indeed~$w_{\lambda} \in W^{I^{0,1}_{\mu}}$. On the other hand, since $W_{I^{0}_{\mu}} \subseteq W_{I^{0,1}_{\mu}}$, the cosets of $W_{I^{0,1}_{\mu}}$ are unions of cosets of $W_{I^{0}_{\mu}}$ and hence the minimal length element of any coset of~$W_{I^{0,1}_{\mu}}$ must be of the form $w_{\lambda}$ for some $\lambda \in P$ with $\lambda_{\mathrm{dom}} = \mu$.
\end{proof}

\begin{lemma} \label{lem:symsinks}
For any $\mathbf{k} \in \mathbb{N}[\Phi]^W$, the sinks of $\Gamma_{\mathrm{sym},\mathbf{k}}$ are
\[\{\eta_{\mathbf{k}}(\lambda)\colon \lambda \in P, \<\lambda,\alpha^\vee\> \neq -1 \textrm{ for all $\alpha \in \Phi^{+}$}\}\]
\end{lemma}
\begin{proof}
First suppose that $\lambda \in P$ satisfies $\<\lambda,\alpha^\vee\> \neq -1$ for all $\alpha \in \Phi^{+}$. Let $\alpha \in \Phi^{+}$. If~$\alpha \in w_{\lambda}(\Phi_{I^{0,1}_{\lambda}})$, then $\<\eta_{\mathbf{k}}(\lambda),\alpha^\vee\> = \<\lambda_{\mathrm{dom}} + \rho_{\mathbf{k}},w_{\lambda}^{-1}(\alpha)^\vee\> \geq \mathbf{k}(\alpha)$ since $w_{\lambda}^{-1}(\alpha) \in \Phi^{+}$ by Proposition~\ref{prop:posparaboliccosets}. So now consider $\alpha \notin w_{\lambda}(\Phi_{I^{0,1}_{\lambda}})$. Then $w_{\lambda}^{-1}(\alpha)$ may be positive or negative, but $|\<\lambda_{\mathrm{dom}},w_{\lambda}(\alpha)^\vee\>| \geq 2$ (because $\lambda_{\mathrm{dom}}$ has an $\omega_i$ coefficient of at least~$2$ for some $i \notin I^{0,1}_{\lambda}$ such that $\alpha_i^\vee$ appears in the expansion of $\pm w_{\lambda}(\alpha)^\vee$). Hence 
\[|\<\eta_{\mathbf{k}}(\lambda),\alpha^\vee\>| = |\<\lambda_{\mathrm{dom}} + \rho_{\mathbf{k}},w_{\lambda}^{-1}(\alpha)^\vee\>| \geq \mathbf{k}(\alpha)+2,\] 
which means that $\<\eta_{\mathbf{k}}(\lambda),\alpha^\vee\> \notin [-\mathbf{k}(\alpha)-1,\mathbf{k}(\alpha)-1]$. Thus $\eta_{\mathbf{k}}(\lambda)$ is a sink of~$\Gamma_{\mathrm{sym},\mathbf{k}}$.

Now suppose $\mu$ is a sink of $\Gamma_{\mathrm{sym},\mathbf{k}}$. Since $\<\mu,\alpha^\vee\> \notin [-\mathbf{k}(\alpha)-1,\mathbf{k}(\alpha)-1]$ for~$\alpha \in \Phi^{+}$, in particular $|\<\mu,\alpha^\vee\>| \geq \mathbf{k}(\alpha)$ for all $\alpha\in \Phi^{+}$. This means~$\<\mu_{\mathrm{dom}},\alpha^\vee\> \geq \mathbf{k}(\alpha)$ for all~$\alpha \in \Phi^{+}$. Hence $\mu_{\mathrm{dom}} = \mu' + \rho_{\mathbf{k}}$ for some dominant $\mu' \in P_{\geq 0}$. Suppose to the contrary that $w_{\mu}$ is not the minimal length element of $w_{\mu}W_{I^{0,1}_{\mu'}}$. Then there exists a descent $s_i$ of~$w_{\mu}$ with~$i \in I^{0,1}_{\mu'}$. But then
\[\<\mu,-w_{\mu}(\alpha_i)^\vee\> = \<\mu_{\mathrm{dom}},-\alpha_i^\vee\> = -\<\mu',\alpha_i^\vee\> -\<\rho_{\mathbf{k}},\alpha_i^\vee\> \geq -\mathbf{k}(\alpha_i)-1,\] 
and also $\<\mu,-w_{\mu}(\alpha_i)^\vee\> =-\< \mu_{\mathrm{dom}},\alpha_i^\vee\> \leq 0$.  This would imply that $\mu$ is not a sink of~$\Gamma_{\mathrm{sym},\mathbf{k}}$, since $-w_{\mu}(\alpha_i) \in \Phi^{+}$. So $w_{\mu}$ must be the minimal length element of~$w_{\mu}W_{I^{0,1}_{\mu'}}$. Thanks to Proposition~\ref{prop:sinksminlen}, this means $w_{\mu} = w_{\lambda}$ for some $\lambda \in P$ with $\lambda_{\mathrm{dom}} = \mu'$ and~$\<\lambda,\alpha^\vee\> \neq -1$ for all $\alpha \in \Phi^{+}$. Moreover, $\mu = w_{\mu}(\mu_{\mathrm{dom}}) = \lambda + w_{\lambda}(\rho_{\mathbf{k}}) = \eta_{\mathbf{k}}(\lambda)$, as claimed.
\end{proof}

\section{Traverse lengths of permutohedra}

Our goal will now be to describe the connected components of $\Gamma_{\mathrm{sym},\mathbf{k}}$, with the eventual aim of establishing confluence of $\raU{\mathrm{sym},\mathbf{k}}$. (By \emph{connected component} of a directed graph, we mean a connected component of its underlying undirected graph.) We will show over the course of the next several sections that the connected components are contained in certain permutohedra; from this confluence will follow easily. First we need to discuss traverse lengths.

\begin{definition}
For a root $\alpha\in\Phi$, an \emph{$\alpha$-string} of length $\ell$ is a subset of $P$ of the form $\{\mu,\mu-\alpha,\mu-2\alpha,\dots,\mu-\ell\alpha\}$ for some weight $\mu \in P$. For a dominant weight $\lambda\in P_{\geq 0}$, an \emph{$\alpha$-traverse} in the discrete permutohedron $\Pi^Q(\lambda)$ is a maximal (as a set) $\alpha$-string that belongs to $\Pi^Q(\lambda)$. Concretely, it is an $\alpha$-string~$\{\mu,\mu-\alpha,\mu-2\alpha,\dots,\mu-\ell\alpha\}\subseteq \Pi^Q(\lambda)$ such that $\mu+\alpha,\, \mu-(\ell+1)\alpha\not\in \Pi^Q(\lambda)$. Finally, for a dominant weight $\lambda \in P_{\geq 0}$, the \emph{traverse length} $\mathbf{l}_\lambda \in \mathbb{Z}[\Phi]^W$ is given by
\[ \mathbf{l}_{\lambda}(\alpha) \coloneqq  \textrm{the minimal length $\ell$ of an $\alpha$-traverse in } \Pi^Q(\lambda).\]
Clearly, by the $W$-symmetry of permutohedra, the traverse length is $W$-invariant and hence really does belong to $\mathbb{Z}[\Phi]^W$.
\end{definition}

\begin{lemma} \label{lem:traversesym}
For $\lambda\in P$ and $\alpha\in \Phi$, any $\alpha$-traverse $\{\mu,\mu-\alpha,\dots,\mu-\ell\alpha\} \subseteq \Pi^Q(\lambda)$ is symmetric with respect to the reflection $s_\alpha$, i.e., $s_\alpha(\mu - i\alpha) = \mu-(\ell-i)\alpha$ for all~$i=0,\dots,l$. Its length is $\ell=\<\mu,\alpha^\vee\>$. In particular, $\<\mu,\alpha^\vee\>\geq 0$.
\end{lemma}
\begin{proof} 
By the $W$-symmetry of discrete permutohedra, we have $s_\alpha(\Pi^Q(\lambda)) = \Pi^Q(\lambda)$, which implies the first sentence.
The second sentence then follows from 
\[\mu-\ell\alpha = s_\alpha(\mu) = \mu - \<\mu,\alpha^\vee\>\,\alpha.\]
The last sentence is clear because the length $\ell$ must be nonnegative.
\end{proof}

Lemma~\ref{lem:traversesym} implies the following reformulation of the definition of~$\mathbf{l}_\lambda $.

\begin{cor} \label{cor:traverselenreform}
For $\lambda \in P$, the traverse length $\mathbf{l}_\lambda$ is given by
\[\mathbf{l}_{\lambda}(\alpha) = \mathrm{min}( \{\<\mu,\alpha^\vee\>\colon \mu\in \Pi^Q(\lambda),\, \mu + \alpha\not\in\Pi^Q(\lambda)\}).\]
\end{cor}

Corollary~\ref{cor:traverselenreform} explains the connection of traverse length to interval-firing: we are going to prove that interval-firing processes get ``trapped'' inside of permutohedra because the traverse lengths of these permutohedra are large (and hence if $\mu$ is inside such a permutohedron but $\mu +\alpha$ is not, $\<\mu,\alpha^\vee\>$ must be so large that it is outside the fireability interval of our process). To do this we need a formula for traverse length. In most cases, the traverse length of a permutohedron in a given direction~$\alpha$ is realized on some edge of the permutohedron in direction~$\alpha$. However, there are some strange exceptions to this general rule, for which we need the concept of ``funny'' weights.

\begin{definition} \label{def:funny}
If $\Phi$ is simply laced, then there are no funny weights. So suppose~$\Phi$ is not simply laced. Then there is a unique long simple root $\alpha_l$ and short simple root~$\alpha_s$ with $\<\alpha_l,\alpha^\vee_s\> \neq 0$. We say the dominant weight $\lambda = \sum_{i=1}^{n}c_i\omega_i \in P_{\geq 0}$ is \emph{funny} if~$c_s = 0$ and $c_l \geq 1$ and $c_i \geq c_l$ for all $i$ such that $\alpha_i$ is long.
\end{definition}

\begin{example}
With the numbering of simple roots as in Figure~\ref{fig:dynkinclassification}, if~$\Phi = B_n$ then $\lambda= \sum_{i=1}^{n}c_i\omega_i \in P_{\geq 0}$ is funny if $c_1,\ldots,c_{n-2} \geq c_{n-1} \geq 1$ and $c_n = 0$. If~$\Phi=C_n$, then $\lambda$ is funny if $c_{n-1}=0$ and $c_n \geq 1$.
\end{example}

For a dominant weight $\lambda = \sum_{i=1}^{n}c_i\omega_i \in P_{\geq 0}$, define $\mathbf{m}_{\lambda}\in\mathbb{Z}[\Phi]^W$ by setting \[\mathbf{m}_{\lambda}(\alpha) \coloneqq  \mathrm{min}(\{c_i\colon \alpha \in W(\alpha_i)\}).\]

\begin{thm} \label{thm:traverseformula}
For a dominant weight $\lambda \in P_{\geq 0}$, we have
\[\mathbf{l}_{\lambda}(\alpha) = \begin{cases}\mathbf{m}_{\lambda}(\alpha)-1 &\textrm{if $\alpha$ is long and $\lambda$ is funny},\\ \mathbf{m}_{\lambda}(\alpha) &\textrm{otherwise}.  \end{cases}\]
\end{thm}
\begin{proof}
Let $\lambda = \sum_{i=1}^{n}c_i\omega_i \in P_{\geq 0}$. The $\alpha_i$-traverse $\{\lambda,\lambda-\alpha,\ldots,\lambda-\ell\alpha=s_i(\lambda)\}$, which is contained in the edge $[\lambda,s_i(\lambda)]$ of the permutohedron $\Pi(\lambda)$, has length equal to $\ell = \<\lambda,\alpha_i^\vee\> = c_i$. By the $W$-symmetry of the traverse length (and because any root is $W$-conjugate to some simple root), it follows that $\mathbf{l}_{\lambda} \leq \mathbf{m}_{\lambda}$.

We will show that in most of the cases (except the case with long roots and funny weights) we actually have $\mathbf{l}_{\lambda} = \mathbf{m}_{\lambda}$. We need to show that the length of any $\alpha$-traverse in $\Pi^Q(\lambda)$ is greater than or equal to $\mathbf{m}_{\lambda}(\alpha)$, i.e., for $\mu \in \Pi^{Q}(\lambda)$ such that $\mu +\alpha \notin \Pi^Q(\lambda)$, we have $\<\mu,\alpha^\vee\> \geq \mathbf{m}_{\lambda}(\alpha)$.

If $\mathbf{m}_{\lambda}(\alpha) = 0$, then we automatically get $\mathbf{l}_{\lambda}(\alpha) = \mathbf{m}_{\lambda}(\alpha)=0$, because $\mathbf{l}_{\lambda}(\alpha) \geq 0$. So let us assume that $\mathbf{m}_{\lambda}(\alpha) \geq 1$.

Let $\mu \in \Pi^Q(\lambda)$ be such that $\mu + \alpha \notin \Pi^Q(\lambda)$. Since $\mu +\alpha \in Q+\lambda$, we deduce that $\mu + \alpha \notin \Pi(\lambda)$. This means that the line segment $[\mu,\mu+\alpha]$ must ``exit'' the permutohedron $\Pi(\lambda)$ at some point $v \in V$, i.e., there exists a unique point $v = \mu+t\alpha$, where $t \in \mathbb{R}$, with $v \in \Pi(\lambda)$ but $\mu + q\alpha \notin \Pi(\lambda)$ for any $q > t$. We have $0 \leq t < 1$.

Let $F$ be the minimal (by inclusion) face of $\Pi(\lambda)$ that contains the point $v$. The minimal value of the linear form $\<\cdot,\alpha^\vee\>$ on the face $F$ should be reached at a vertex $\nu$ of $F$. By the $W$-symmetry of $\Pi(\lambda)$, we assume without loss of generality that this minimum is achieved at $\nu=\lambda$. So we have $\<\lambda,\alpha^\vee\> \leq \<v,\alpha^\vee\>$.

If $\lambda$ is strictly in the fundamental chamber, then any edge of $\Pi(\lambda)$ coming out of~$\lambda$ must be in the direction of a negative simple root. This is not true for general~$\lambda \in P_{\geq 0}$, but the edges of $\Pi(\lambda)$ coming out of $\lambda$ that are not in the direction of a negative simple root must immediately leave the dominant chamber. Hence if we let~$x \in V$ be some generic point in the interior of the face $F$ very close to $\lambda$, by acting by $W_{I^{0}_{\lambda}}$ we can transport $x$ to the dominant chamber while fixing $\lambda$. Thus, we may assume that the affine span of~$F$ is spanned by simple roots. So let $I \subseteq [n]$ be the minimal set of indices such that the face $F$ belongs to the affine subspace $\lambda + \Span_{\mathbb{R}}(\{\alpha_i\colon i \in I\})$.

Let $\alpha = \sum_{i=1}^{n}a_i\alpha_i$, where the $a_i$ are either all nonnegative or all nonpositive. Then we have~$\alpha^\vee = \sum_{i=1}^{n} \widetilde{a}_i\alpha^\vee_i$ where $\widetilde{a}_i = \frac{\<\alpha_i,\alpha_i\>}{\<\alpha,\alpha\>}a_i$. Note that these $\widetilde{a}_i$ are also integers.

Any root $\alpha$ is $W$-conjugate to at least one simple root that appears with nonzero coefficient in its expansion in terms of the simple roots. So there exists $j \in [n]$ such that $\alpha_j \in W(\alpha)$ and~$a_j = \widetilde{a}_j \neq 0$. We have $c_j \geq \mathbf{m}_{\lambda}(\alpha) \geq 1$.

We have $\lambda \in \Pi(\lambda)$ and $\lambda + \alpha \notin \Pi(\lambda)$. So $\<\lambda,\alpha^\vee\> \geq 0$, because $\<\lambda,\alpha^\vee\>$ is the length of the $\alpha$-traverse that starts at $\lambda$, which is always nonnegative. Therefore we have~$\<\lambda,\alpha^\vee\> = \sum_{i=1}^{n} \widetilde{a}_i c_i \geq 0$; moreover, all nonzero terms in this expression have the same sign and at least one term $\widetilde{a}_jc_j$ is nonzero. It follows that $a_1,\ldots,a_n \geq 0$, i.e., $\alpha$ is a positive root.

We have $\mu = v - t\alpha = (\lambda - \sum_{i \in I}b_i\alpha_i) - t\alpha$ for real numbers $0 \leq t < 1$ and $b_i \geq 0$, $i \in I$. Thus $\<\mu,\alpha^\vee\> = \<v,\alpha^\vee\> - t\<\alpha,\alpha^\vee\> = \<v,\alpha^\vee\> - 2t \geq \<\lambda,\alpha^\vee\> - 2t > \<\lambda,\alpha^\vee\> - 2$. Moreover, since both $\<\mu,\alpha^\vee\>$ and $\<\lambda,\alpha^\vee\>-2$ are integers, and the first is strictly greater than the second, we get
\[ \<\mu,\alpha^\vee\> \geq \<\lambda,\alpha^\vee\>- 1 = \left(\sum_{i=1}^{n}\widetilde{a}_ic_i\right) - 1.\]
We already noted that the last expression involves at least one nonzero term $\widetilde{a}_jc_j$ such that $\alpha_j \in W(\alpha)$. So $\widetilde{a}_jc_j\geq c_j \geq \mathbf{m}_{\lambda}(\alpha)$ and thus $\<\mu,\alpha^\vee\> \geq \mathbf{m}_{\lambda}(\alpha) - 1$.

We need to prove just a slightly stronger inequality  $\<\mu,\alpha^\vee\> \geq \mathbf{m}_{\lambda}(\alpha)$.

If $\sum_{\alpha_i \in W(\alpha)} \widetilde{a}_i \geq 2$, we get
\[ \<\mu,\alpha^\vee\> \geq \sum_{i=1}^{n}\widetilde{a}_ic_i - 1 \geq \sum_{\alpha_i \in W(\alpha)} \widetilde{a}_ic_i - 1 \geq 2\mathbf{m}_{\lambda}(\alpha)-1 \geq \mathbf{m}_{\lambda}(\alpha),\]
as needed. So we now assume that $\sum_{\alpha_i \in W(\alpha)} \widetilde{a}_i =1$. Note that this means $a_j = \widetilde{a}_j=1$.

If we had $c_j > \mathbf{m}_{\lambda}(\alpha)$, then we would get 
\[\<\mu,\alpha^\vee\> \geq \widetilde{a}_jc_j -1 \geq c_j-1 \geq \mathbf{m}_{\lambda}(\alpha)\] 
and we would also be done. So we now assume that $c_j = \mathbf{m}_{\lambda}(\alpha)$.

Since $\alpha$ does not belong to the subspace spanned by the $\alpha_i$ for $i \in I$, there is $r\in[n]$ with $r\notin I$ such that $a_r \geq 1$.

If $a_r=1$, then, from the fact that $\lambda-\mu = (\sum_{i\in I}b_i\alpha_i) + t\alpha$ belongs to the root lattice $Q$ and thus is an \emph{integer} linear combination of the simple roots, we deduce that in fact~$t \in \mathbb{Z}$ and thus $t = 0$. In this case get $\<\mu,\alpha^\vee\> \geq \<\lambda,\alpha^\vee\> \geq a_jc_j \geq \mathbf{m}_{\lambda}(\alpha)$, as needed. So we now assume that $a_r \geq 2$.

Then note that $\alpha_r \notin W(\alpha)$, because we assumed~$\sum_{\alpha_i \in W(\alpha)} \widetilde{a}_i  = \sum_{\alpha_i \in W(\alpha)} a_i = 1$.

If there is $q \in [n]$ such that $a_q \notin W(\alpha)$, $\widetilde{a}_q \geq 1$ and $c_q \geq 1$, we have
\[ \<\mu,\alpha^\vee\> = \left( \sum_{i=1}^{n} \widetilde{a}_ic_i \right) - 1 \geq \widetilde{a}_jc_j + \widetilde{a}_qc_q - 1 \geq  \widetilde{a}_jc_j \geq \mathbf{m}_{\lambda}(\alpha),\]
as needed.

Thus, the only possibility which is not covered by the above discussion is when:
\begin{enumerate}[(1)]
\item \label{cond:funny1} There is exactly one nonzero term $a_j\alpha_j$ in the expansion $\alpha = \sum_{i=1}^{n}a_i\alpha_i$ such that $\alpha_j \in W(\alpha)$. For this term, $a_j=1$ and $c_j = \mathbf{m}_{\lambda}(\alpha) \geq 1$.
\item \label{cond:funny2} There is at least one more more nonzero term $a_i\alpha_i$ in that expansion. For all such terms, $\alpha_i \notin W(\alpha)$, $a_i \geq 2$, and $c_i = 0$.
\end{enumerate}
We claim that these conditions imply that $\alpha$ is a long root. This is easy to check by hand for $\Phi=B_n$, $C_n$, or $G_2$. One does not need to check Type $F_4$ separately, because in this case there are two long simple roots and two short simple roots, but the expansion of $\alpha$ involves either only one short simple root or only one long simple root. We leave it as an exercise for the reader to find a uniform root theoretic argument of the fact that conditions~\eqref{cond:funny1} and~\eqref{cond:funny2} above imply that $\alpha$ is long.

Also, we claim that conditions~\eqref{cond:funny1} and~\eqref{cond:funny2} above imply that $\lambda$ is a funny weight. Indeed, it is a well-known fact that for any root $\alpha = \sum_{i=1}^{n}a_i\alpha_i$, the set of~$i \in [n]$ for which $a_i \neq 0$ must be a connected subset of the Dynkin diagram (see for instance~\cite[Chapter~VI, \S1.6, Corollary~3]{bourbaki2002lie}). Hence the $\alpha_j$ in condition~\eqref{cond:funny1} must be the long simple root $\alpha_l$, and one of the $\alpha_i$ in condition~\eqref{cond:funny2} must be the short simple root $\alpha_s$ (with notation as in Definition~\ref{def:funny}). Note also that $\mathbf{m}_{\lambda}(\alpha)=c_l$ forces $c_i \geq c_l$ for all $i$ such that $\alpha_i$ is long.

In this ``long and funny'' case we can only get the (slightly) weaker inequality: 
\[\<\mu,\alpha^\vee\> \geq \mathbf{m}_{\lambda}(\alpha)-1.\]

It remains to show that this last inequality is tight in this ``long and funny'' case. Let us concentrate on the $2$-dimensional face of the permutohedron $\Pi(\lambda)$ contained in the affine subspace $\lambda + \Span_{\mathbb{R}}(\{\alpha_l,\alpha_s\})$ (with notation as in Definition~\ref{def:funny}).

This face is equivalent to the $2$-dimensional $W'$-permutohedron $\Pi_{W'}(\lambda')$ corresponding to the sub-root system $\Phi'$ of rank $2$ with simple roots $\alpha_l$ and $\alpha_s$, and fundamental weights $\omega'_1$ (corresponding to $\alpha_l$) and $\omega'_2$ (corresponding to $\alpha_s$),  where $W'$ is the Weyl group of $\Phi'$, and $\lambda'=\mathbf{m}_{\lambda}(\alpha)\omega'_1 + 0\cdot\omega'_2$.

The $2$-dimensional root system $\Phi'$ must be equal to either $B_2$ or $G_2$. In this situation there in fact is a $\mu \in \Pi_{W'}^Q(\lambda')$ with $\mu+\alpha \notin \Pi_{W'}^Q(\lambda')$ for some long $\alpha \in \Phi'$ such that $\<\mu,\alpha^\vee\> = \mathbf{m}_{\lambda}(\alpha)-1$: indeed, we can take $\alpha\coloneqq \alpha_l$ and $\mu \coloneqq  (\mathbf{m}_{\lambda}(\alpha)-1)\omega'_1$ for $B_2$ or $\alpha\coloneqq \alpha_l$ and $\mu \coloneqq  (\mathbf{m}_{\lambda}(\alpha)-1)\omega'_1 + \omega'_2$ for $G_2$.

 This finishes the proof of the theorem.
\end{proof}

\section{The permutohedron non-escaping lemma}

We need to place some restrictions on our parameter $\mathbf{k}$ so that funny weights do not occur in our analysis of the relevant permutohedra traverse lengths. For this we have the notion of ``goodness.''

\begin{definition}
If $\Phi$ is simply laced, then every $\mathbf{k} \in \mathbb{N}[\Phi]^W$ is good. So suppose~$\Phi$ is not simply laced and let $\mathbf{k} \in \mathbb{N}[\Phi]^W$. Then there exist $k_s, k_l \in \mathbb{Z}$ with~$\mathbf{k}(\alpha) = k_s$ if~$\alpha$ is short and~$\mathbf{k}(\alpha) = k_l$ if~$\alpha$ is long. We say $\mathbf{k}$ is \emph{good} if $k_s = 0 \Rightarrow k_l=0$. Note in particular that if $\mathbf{k}=k\geq 0$ is constant, then it is good.
\end{definition}

Now we can prove the following permutohedron non-escaping lemma, which says that certain discrete permutohedra ``trap'' the symmetric interval-firing process inside of them. 

\begin{lemma} \label{lem:permtrap}
Let $\mathbf{k} \in \mathbb{N}[\Phi]^W$ be good and let $\Gamma \coloneqq  \Gamma^{\mathrm{un}}_{\mathrm{sym},\mathbf{k}}$. Let $\lambda \in P_{\geq 0}$. Then there is no directed edge $(\mu,\mu')$ in $\Gamma$ with $\mu \in \Pi^Q(\eta_{\mathbf{k}}(\lambda))$ and $\mu' \notin \Pi^Q(\eta_{\mathbf{k}}(\lambda))$.
\end{lemma}
\begin{proof}
First suppose $\Phi$ is not simply laced and $k_s=0$. Then also $k_l =0$, i.e., $\mathbf{k} = 0$, since $\mathbf{k}$ is good. Hence $\rho_{\mathbf{k}} = 0$, so $\eta_{\mathbf{k}}(\lambda) = \lambda$. If $\mu \in  \Pi^Q(\lambda)$ but $\mu + \alpha \notin  \Pi^Q(\lambda)$, then by Corollary~\ref{cor:traverselenreform} we have~$\<\mu,\alpha^\vee\> \geq \mathbf{l}_{\lambda}(\alpha)$. Note that by definition $ \mathbf{l}_{\lambda}(\alpha) \geq 0$. But this means $\<\mu,\alpha^\vee\> \geq \mathbf{k}(\alpha)$, so indeed $(\mu,\mu+\alpha)$ cannot be a directed edge of $\Gamma$.

Now suppose either $\Phi$ is simply laced or $\Phi$ is not simply laced but $k_s \geq 1$. Then note that $\rho_{\mathbf{k}}$ is not funny. Hence by Theorem~\ref{thm:traverseformula} we conclude that~$\mathbf{l}_{\eta_{\mathbf{k}}(\lambda)}(\alpha) \geq \mathbf{k}(\alpha)$. If~$\mu \in  \Pi^Q(\eta_{\mathbf{k}}(\lambda))$ but $\mu + \alpha \notin  \Pi^Q(\eta_{\mathbf{k}}(\lambda))$, then $\<\mu,\alpha^\vee\> \geq \mathbf{l}_{\eta_{\mathbf{k}}(\lambda)}(\alpha)$ by Corollary~\ref{cor:traverselenreform}. This means $\<\mu,\alpha^\vee\> \geq \mathbf{k}(\alpha)$, so indeed $(\mu,\mu+\alpha)$ cannot be a directed edge of $\Gamma$.
\end{proof}

\begin{figure}
\begin{center}
\begin{tikzpicture}
  
\node[scale=\scl,draw,circle,fill=red] (5) at (0,-1) {};
\node[scale=\scl,draw,circle,fill=red] (6) at (-1,0) {};
\node[scale=\scl,draw,circle,fill=red] (7) at (1,0) {};
\node[scale=\scl,draw,circle,fill=red] (8) at (0,1) {};
\fill[red!20] (5.center)--(7.center)--(8.center)--(6.center)--cycle;
\node[scale=\scl,draw,circle,fill=red] (5) at (0,-1) {};
\node[scale=\scl,draw,circle,fill=red] (6) at (-1,0) {};
\node[scale=\scl,draw,circle,fill=red] (7) at (1,0) {};
\node[scale=\scl,draw,circle,fill=red] (8) at (0,1) {};
\node[anchor=north west,red] at (7) {$\Pi(\rho_{\mathbf{k}})$};

\node[scale=\scl,draw,circle,fill=black] (1) at (0,0) {};
\node[scale=\scl,draw,circle,fill=black] (2) at (-1,1) {};
\node[scale=\scl,draw,circle,fill=black] (3) at (1,1) {};
\node[scale=\scl,draw,circle,fill=black] (4) at (0,2) {};

\node[anchor=north] at (1.south) {$0$};
\node[anchor=east] at (2.west) {$\alpha_1$};
\node[anchor=west] at (3.east) {$\alpha_1+2\alpha_2$};

\draw[->,>=stealth,thick] (1)--(2);
\draw[->,>=stealth,thick] (2)--(4);
\draw[->,>=stealth,thick] (1)--(3);
\draw[->,>=stealth,thick] (3)--(4);

\node[anchor=south] at (8.north) {$\omega_1$};

\draw[->,>=stealth,thick] (5)--(6);
\draw[->,>=stealth,thick]  (6)--(8);
\draw[->,>=stealth,thick] (5)--(7);
\draw[->,>=stealth,thick]  (7)--(8);

\node (9) at (-2,0) {};
\node (10) at (-1,-1) {};
\node (11) at (1,-1) {};
\node (12) at (2,0) {};

\draw[->,>=stealth,thick] (9) -- (2);
\draw[->,>=stealth,thick] (10) -- (1);
\draw[->,>=stealth,thick] (11) -- (1);
\draw[->,>=stealth,thick] (12) -- (3);

\end{tikzpicture}
\end{center}
\caption[Failure of permutohedron non-escaping for not good~$\mathbf{k}$]{The graph $\Gamma_{\mathrm{sym},\mathbf{k}}$ in Example~\ref{ex:notgood}. The permutohedron $\Pi(\rho_{\mathbf{k}})$ is shown in red.} \label{fig:notgood}
\end{figure}
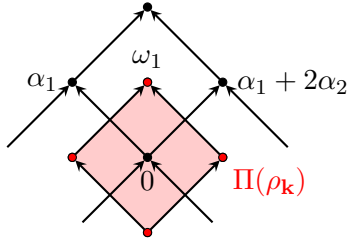

\begin{example} \label{ex:notgood}
Lemma~\ref{lem:permtrap} is false in general without the goodness assumption. Indeed, suppose $\Phi=B_2$ and $\mathbf{k} \in \mathbb{N}[\Phi]^W$ is given by $k_s\coloneqq 0$ and $k_l\coloneqq 1$. Then Figure~\ref{fig:notgood} depicts (a portion of) the graph $\Gamma_{\mathrm{sym},\mathbf{k}}$. In this picture we only show elements of the root lattice~$Q$. The permutohedron~$\Pi(\rho_{\mathbf{k}}) = \Pi(\omega_1)$ is shown in red. Observe that although $0 \in \Pi^Q(\rho_{\mathbf{k}})$and~$\alpha_1 \notin \Pi^Q(\rho_{\mathbf{k}})$, we have an edge $(0,\alpha_1)$ in $\Gamma_{\mathrm{sym},\mathbf{k}}$.
\end{example}

We also need a ``lower-dimensional'' version of the permutohedron non-escaping lemma that says that these interval-firing processes get trapped inside of permutohedra of parabolic subgroups of $W$. This is established in the following lemma and corollary.

\begin{lemma} \label{lem:permtrapsmall}
Let $\mathbf{k} \in \mathbb{N}[\Phi]^W$ and $\Gamma \coloneqq  \Gamma^{\mathrm{un}}_{\mathrm{sym},\mathbf{k}}$. Let~$\lambda \in P_{\geq 0}$. Then if $(\mu,\mu+\alpha)$ is a directed edge in $\Gamma$ with~$\mu \in \Pi_{I^{0,1}_\lambda}^Q(\eta_{\mathbf{k}}(\lambda))$, we have~$\alpha \in \Phi_{I^{0,1}_\lambda}$.
\end{lemma}
\begin{proof}
Write $\eta_{\mathbf{k}}(\lambda) = \sum_{i=1}^{n}c_i\omega_i$. Assume to the contrary that there exists an edge $(\mu,\mu+\alpha)$ in $\Gamma$ such that $\mu \in  \Pi_{I^{0,1}_\lambda}^Q(\eta_{\mathbf{k}}(\lambda))$ but $\alpha$ does not belong to $\Span_{\mathbb{R}}(\{\alpha_i\colon i \in I\})$.

Note that $\alpha$ is a root (positive or negative) with $-\mathbf{k}(\alpha)-1\leq \<\mu,\alpha^\vee\> \leq \mathbf{k}(\alpha)+1$. Let $\beta = \pm\alpha \in \Phi^{+}$ be the positive root. Then~$\<\mu,\beta^\vee\> \leq \mathbf{k}(\alpha)+1 \leq \mathbf{k}(\beta) + 1$.

Since the point $\mu$ belongs to $\Pi_{I^{0,1}_{\lambda}}(\eta_{\mathbf{k}}(\lambda))$, we deduce that the same inequality $\<\nu,\beta^\vee\> \leq \mathbf{k}(\beta)+1$ holds for some vertex $\nu$ of $\Pi_{I^{0,1}_{\lambda}}(\eta_{\mathbf{k}}(\lambda))$. We have $\nu= w(\eta_{\mathbf{k}}(\lambda))$ where $w \in W_{I^{0,1}_{\lambda}}$. Hence we have that~$\<w(\lambda),\beta^\vee\> = \<\lambda,w^{-1}(\beta)^\vee\> \leq \mathbf{k}(\beta)+1$ for some~$w \in W_{I^{0,1}_{\lambda}}$.

The action of the parabolic subgroup $W_{I^{0,1}_{\lambda}}$ on $\beta^\vee$ does not change the coefficients~$b_j$ of the expansion $\beta^\vee = \sum_{i=1}^{n}b_i\alpha_i^\vee$ for all $j \notin I$, and at least one of these coefficients~$b_j$ should be strictly positive (because $\beta^\vee$ is a positive coroot that does not belong to~$\Span_{\mathbb{R}}(\{\alpha_i\colon i \in I\})$). So the expansion $w^{-1}(\beta)^\vee = \sum_{i=1}^{n}b'_i\alpha_i^\vee$ contains some strictly positive coefficient, which means that $w^{-1}(\beta)^\vee$ is a positive coroot and thus we have~$b'_i \geq 0$ for all $i$.

Moreover, any coroot is $W$-conjugate to some simple coroot that appears in its expansion with nonzero coefficient. These observations mean that we can find $j\notin I$ such that $b'_j=b_j \geq 1$, and also (possibly the same) $i$ such that $b'_i\geq 1$ and $\alpha_i \in W(\alpha)$. Note that for this $i$ we have $\mathbf{k}(\alpha_i) = \mathbf{k}(w^{-1}(\beta)) = \mathbf{k}(\beta) = \mathbf{k}(\alpha)$.

If $i=j$, we get $\<\lambda,w^{-1}(\beta)^\vee\> \geq \<\lambda,b'_j\alpha_j^\vee\> \geq \<\lambda,\alpha_j^\vee\> = c_j \geq \mathbf{k}(\alpha_j)+2 = \mathbf{k}(\alpha)+2$ (because for $j \notin I$, $c_j\geq \mathbf{k}(\alpha_j)+2$). But this contradicts~$\<\lambda,w^{-1}(\beta)^\vee\> \leq \mathbf{k}(\alpha)+1$.

On the other hand, if $i \neq j$, we get 
\begin{align*}
\<\lambda,w^{-1}(\beta)^\vee\> \geq \<\lambda,b'_i\alpha_i^\vee+b'_j\alpha_j^\vee\> &\geq \<\lambda,\alpha_i^\vee\> + \<\lambda,\alpha_j^\vee\> = c_i + c_j \\
&\geq \mathbf{k}(\alpha_i) + (\mathbf{k}(\alpha_j) + 2) \geq \mathbf{k}(\alpha_i) + 2 = \mathbf{k}(\alpha)+2.
\end{align*}
Again, we get a contradiction.
\end{proof}

\begin{cor} \label{cor:permtrap}
Let $\mathbf{k} \in \mathbb{N}[\Phi]^W$ be good and $\Gamma \coloneqq  \Gamma^{\mathrm{un}}_{\mathrm{sym},\mathbf{k}}$. Let~$\lambda \in P_{\geq 0}$. Then there is no directed edge $(\mu,\mu')$ in $\Gamma$ with $\mu \in \Pi_{I^{0,1}_\lambda}^Q(\eta_{\mathbf{k}}(\lambda))$ and~$\mu' \notin \Pi_{I^{0,1}_\lambda}^Q(\eta_{\mathbf{k}}(\lambda))$.
\end{cor}
\begin{proof}
This follows by combining Lemmas~\ref{lem:permtrap} and~\ref{lem:permtrapsmall}: if we have a directed edge $(\mu,\mu+\alpha)$ with $\mu \in \Pi_{I^{0,1}_\lambda}^Q(\eta_{\mathbf{k}}(\lambda))$, then $\alpha \in \Phi_{I^{0,1}_\lambda}$ by Lemma~\ref{lem:permtrapsmall}; hence this firing move is equivalent (via projection) to the same move for the sub-root system $ \Phi_{I^{0,1}_\lambda}$; so by Lemma~\ref{lem:permtrap} applied to that sub-root system, we have $\mu+\alpha \in \Pi_{I^{0,1}_\lambda}^Q(\eta_{\mathbf{k}}(\lambda))$.
\end{proof}

\section{Confluence of symmetric interval-firing} \label{sec:sym_conf}

Now, as promised, we are ready to show that connected components of $\Gamma_{\mathrm{sym},\mathbf{k}}$ are contained inside permutohedra.

\begin{thm} \label{thm:permcc}
Let $\mathbf{k} \in \mathbb{N}[\Phi]^W$ be good. Let $\lambda \in P$ with $\<\lambda,\alpha^\vee\> \neq -1$ for all $\alpha \in \Phi^{+}$. Let $Y_{\lambda}\coloneqq \{\mu\in P\colon \mu \baAstU{\mathrm{sym},\mathbf{k}} \eta_{\mathbf{k}}(\lambda) \}$ be the connected component of $\Gamma_{\mathrm{sym},\mathbf{k}}$ containing the sink~$\eta_{\mathbf{k}}(\lambda)$. Then $Y_{\lambda}$ is contained in $w_{\lambda}\Pi_{I^{0,1}_{\lambda}}^Q(\eta_{\mathbf{k}}(\lambda_{\mathrm{dom}}))$.
\end{thm}
\begin{proof}
First suppose that $\lambda$ is dominant. By Corollary~\ref{cor:permtrap} there is no edge $(\mu,\mu')$  in~$\Gamma_{\mathrm{sym},\mathbf{k}}$ where one of $\mu,\mu'$ is in $\Pi_{I^{0,1}_\lambda}^Q(\eta_{\mathbf{k}}(\lambda))$ and the other is not, which implies that~$Y_{\lambda}$ is contained in $\Pi_{I^{0,1}_{\lambda}}^Q(\eta_{\mathbf{k}}(\lambda))$. Now suppose $\lambda$ is not dominant. By the preceding argument, the result is true for~$\lambda_{\mathrm{dom}}$. But then we have $Y_{\lambda} = w_{\lambda} Y_{\lambda_{\mathrm{dom}}}$ by the $W$-symmetry of~$\Gamma^{\mathrm{un}}_{\mathrm{sym},\mathbf{k}}$, i.e., by Theorem~\ref{thm:symmetry}.
\end{proof}

And now we can prove half of Theorem~\ref{thm:confluence_intro}.

\begin{cor} \label{cor:symconfluence}
Let $\mathbf{k} \in \mathbb{N}[\Phi]^W$ be good. Then $\raU{\mathrm{sym},\mathbf{k}}$ is confluent (and terminating).
\end{cor}
\begin{proof}
We already saw in Proposition~\ref{prop:interval_firing_terminating} that $\raU{\mathrm{sym},\mathbf{k}}$ is terminating. Thus, every connected component of $\Gamma_{\mathrm{sym},\mathbf{k}}$ contains at least one sink, and $\raU{\mathrm{sym},\mathbf{k}}$ is confluent as long as every connected component contains a unique sink. 

So suppose that two sinks belong to the same connected component of $\Gamma_{\mathrm{sym},\mathbf{k}}$. By Lemma~\ref{lem:symsinks}, we know that these sinks must be of the form $\eta_{\mathbf{k}}(\lambda)$ and $\eta_{\mathbf{k}}(\lambda')$ for~$\lambda,\lambda'\in P$ with $\<\lambda,\alpha^\vee\>\neq -1$ and~$\<\lambda',\alpha^\vee\> \neq -1$ for all $\alpha \in \Phi^{+}$. 

By Theorem~\ref{thm:permcc}, $\eta_{\mathbf{k}}(\lambda) \in w_{\lambda'}\Pi_{I^{0,1}_{\lambda'}}^Q(\eta_{\mathbf{k}}(\lambda'_{\mathrm{dom}}))$ and vice-versa. In particular we have that $\eta_{\mathbf{k}}(\lambda_{\mathrm{dom}}) \in \Pi^Q(\eta_{\mathbf{k}}(\lambda'_{\mathrm{dom}}))$ and $\eta_{\mathbf{k}}(\lambda'_{\mathrm{dom}}) \in \Pi^Q(\eta_{\mathbf{k}}(\lambda_{\mathrm{dom}}))$. Proposition~\ref{prop:perm_containment} then says that~$\eta_{\mathbf{k}}(\lambda_{\mathrm{dom}}) - \eta_{\mathbf{k}}(\lambda'_{\mathrm{dom}})$ and $\eta_{\mathbf{k}}(\lambda'_{\mathrm{dom}}) - \eta_{\mathbf{k}}(\lambda_{\mathrm{dom}})$ are both in $Q_{\geq 0}$, which is possible only if $\eta_{\mathbf{k}}(\lambda_{\mathrm{dom}}) = \eta_{\mathbf{k}}(\lambda'_{\mathrm{dom}})$. That is, thanks to the injectivity of~$\eta_{\mathbf{k}}$ established in Proposition~\ref{prop:etafacts}, we must have $\lambda_{\mathrm{dom}} = \lambda'_{\mathrm{dom}}$. 

But then the fact that $\eta_{\mathbf{k}}(\lambda) \in w_{\lambda'}\Pi_{I^{0,1}_{\lambda}}^Q(\eta_{\mathbf{k}}(\lambda_{\mathrm{dom}}))$ means that $\eta_{\mathbf{k}}(\lambda)$ is a vertex of $ w_{\lambda'}\Pi_{I^{0,1}_{\lambda}}(\eta_{\mathbf{k}}(\lambda_{\mathrm{dom}}))$, i.e., $\eta_{\mathbf{k}}(\lambda) = w_{\lambda'}w(\eta_{\mathbf{k}}(\lambda_{\mathrm{dom}}))$ for some $w \in W_{I^{0,1}_{\lambda}}$. Note that~$(w_{\lambda'}w)^{-1}(\eta_{\mathbf{k}}(\lambda))$ is dominant. We have seen in the the proof of Proposition~\ref{prop:etafacts} that this means $(w_{\lambda'}w)^{-1}(\lambda)$ is dominant as well, or in other words, that $w_{\lambda'}w = w_{\lambda}w'$ for some $w' \in W_{I^{0}_{\lambda}}$. This shows that~$w_{\lambda} \in w_{\lambda'}W_{I^{0,1}_{\lambda}}$. By Proposition~\ref{prop:sinksminlen}, $w_{\lambda}$ and~$w_{\lambda'}$ must both be the minimal length elements of the cosets of $W_{I^{0,1}_{\lambda}}$ they belong to. So~$w_{\lambda} = w_{\lambda'}$. That~$\lambda_{\mathrm{dom}} = \lambda'_{\mathrm{dom}}$ and $w_{\lambda} = w_{\lambda'}$ implies that~$\lambda = \lambda'$, and consequently that~$\eta_{\mathbf{k}}(\lambda)=\eta_{\mathbf{k}}(\lambda')$, as required.
\end{proof}

\begin{remark}
As far as we know, Theorem~\ref{thm:permcc} and Corollary~\ref{cor:symconfluence} may be true even in the case where $\mathbf{k}$ is not good. Indeed, it appears that $\raU{\mathrm{sym},\mathbf{k}}$ is confluent for all~$\mathbf{k} \in \mathbb{N}[\Phi]^W$ and to prove this it would be sufficient, thanks to the diamond lemma (Lemma~\ref{lem:diamond}), to prove it for root systems of rank~$2$, of which there are only four: $A_1 \oplus A_1$, $A_2$, $B_2$, $G_2$. All~$\mathbf{k}$ are good for simply laced root systems, so in fact one would need only check $B_2$ and $G_2$.
\end{remark}

\section{Full-dimensional components, saturated components, and Cartan matrix chip-firing as a limit} \label{sec:cartanmatrix}

Let $\mathbf{k}\in\mathbb{N}[\Phi]^W$ be good. For $\lambda \in P$, recall the notation $Y_{\lambda}\coloneqq \{\mu\in P\colon \mu \baAstU{\mathrm{sym},\mathbf{k}} \eta_{\mathbf{k}}(\lambda) \}$ for the connected component of $\Gamma_{\mathrm{sym},\mathbf{k}}$ containing the sink~$\eta_{\mathbf{k}}(\lambda)$ from the last section. By the results of the last section, all these components are distinct. In this section, we take a moment to highlight certain special components~$Y_{\lambda}$, namely:
\begin{itemize}
\item those which are \emph{full-dimensional} in the sense that their affine hulls are the whole vector space: $\mathrm{AffineHull}\, Y_{\lambda}=V$;
\item those which are full-dimensional and \emph{saturated} in the sense that they contain all lattice points in their convex hulls: $Y_{\lambda} = (\mathrm{ConvexHull}\, Y_{\lambda})\cap(Q+\eta_{\mathbf{k}}(\lambda))$.
\end{itemize}

For the full-dimensional components: by a result we will prove later (Corollary~\ref{cor:symccsweylorbit}), we have that $Y_{\lambda}$ always contains $W(\eta_{\mathbf{k}}(\lambda))$ for $\lambda \in P_{\geq 0}$ with $I^{0,1}_{\lambda} = [n]$. Hence by Theorem~\ref{thm:permcc} we see that the full-dimensional connected components of~$\Gamma_{\mathrm{sym},\mathbf{k}}$ are exactly $Y_{\lambda}$ for $\lambda \in P_{\geq 0}$ with $I^{0,1}_{\lambda}=[n]$, i.e., for those~$\lambda = \sum_{i=1}^{n}c_i\omega_i \in P$ with~$c_i\in \{0,1\}$ for all $i\in [n]$. Clearly there are $2^n$ such full-dimensional components. (Strictly speaking we do not have $\mathrm{AffineHull}\, Y_{\lambda}=V$ when $\lambda =0$ and $\mathbf{k}=0$, but to make our description of full-dimensional components consistent across all values of $\mathbf{k}$ it is best to nevertheless consider this component full-dimensional.)

For the full-dimensional and saturated components: by that same Corollary~\ref{cor:symccsweylorbit}, we see that $Y_{\lambda}$ being full-dimensional and saturated is equivalent to having this component satisfy~$Y_{\lambda}=\Pi^Q(\eta_{\mathbf{k}}(\lambda_{\mathrm{dom}}))$. Recall that $\Omega_m^0$ denotes the set of minuscule weights together with zero; then we have the following:
\begin{prop} \label{prop:saturatedccs}
Let $\mathbf{k} \in \mathbb{N}[\Phi]^W$ be good. Let $\lambda \in P$ be a weight with $\<\lambda,\alpha^\vee\> \neq -1$ for all $\alpha \in \Phi^{+}$. Let $Y_{\lambda}\coloneqq \{\mu\in P\colon \mu \baAstU{\mathrm{sym},\mathbf{k}} \eta_{\mathbf{k}}(\lambda) \}$ be the connected component of $\Gamma_{\mathrm{sym},\mathbf{k}}$ containing the sink~$\eta_{\mathbf{k}}(\lambda)$. Then $Y_{\lambda}$ is equal to $\Pi^{Q}(\eta_{\mathbf{k}}(\lambda))$ if and only if $\lambda \in \Omega_m^0$.
\end{prop}
\begin{proof}
First note that if $\lambda$ is a sink of $\Gamma_{\mathrm{sym},\mathbf{k}}$ then so is $\lambda_{\mathrm{dom}}$ and by the confluence of~$\raU{\mathrm{sym},\mathbf{k}}$ there cannot be two sinks in a single connected component of $\Gamma_{\mathrm{sym},\mathbf{k}}$, so it suffices to prove this proposition for dominant $\lambda \in P_{\geq 0}$ with $I^{0,1}_{\lambda} = [n]$. (Observe that if~$\lambda \in \Omega_m^0$ then certainly it is of this form.) 

By the polytopal characterization of minuscule weights, there exists a dominant weight~$\mu \in P_{\geq 0}$ with $\mu \in \Pi^Q(\lambda)$ but $\mu \neq \lambda$ if and only if $\lambda\notin \Omega_m^0$. Hence by Proposition~\ref{prop:perm_containment} there exists~$\mu \in P_{\geq 0}$ with $\eta_{\mathbf{k}}(\mu) \in \Pi^Q(\eta_{\mathbf{k}}(\lambda))$ but $\eta_{\mathbf{k}}(\mu) \neq \eta_{\mathbf{k}}(\lambda)$ if and only if $\lambda \notin \Omega_m^0$. By applying $W$, we see that there is a sink $\eta_{\mathbf{k}}(\mu)$ of $\Gamma_{\mathrm{sym},\mathbf{k}}$ with~$\eta_{\mathbf{k}}(\mu) \in \Pi^Q(\eta_{\mathbf{k}}(\lambda))$ but $\eta_{\mathbf{k}}(\mu) \notin W(\eta_{\mathbf{k}}(\lambda))$ if and only if $\lambda \notin \Omega_m^0$. Finally, by the permutohedron non-escaping lemma, Lemma~\ref{lem:permtrap}, this means precisely that $ \Pi^Q(\eta_{\mathbf{k}}(\lambda))$ is its own connected component if and only if~$\lambda \in \Omega_m^0$.
\end{proof}

\begin{remark}
Proposition~\ref{prop:saturatedccs} fails when $\mathbf{k}$ is not good, as can be seen in Example~\ref{ex:notgood} above: in this example, $0 \in \Pi^Q(\rho_{\mathbf{k}})$ but $0$ does not belong to the connected component of $\Gamma_{\mathrm{sym},\mathbf{k}}$  containing $\rho_{\mathbf{k}}=\eta_{\mathbf{k}}(0)=\omega_1$.
\end{remark}

So we see that the full-dimensional and saturated components of $\Gamma_{\mathrm{sym},\mathbf{k}}$ are exactly the~$Y_{\omega}$ for $\omega \in \Omega_m^0$. There are $f$ of these, where we recall that $f \coloneqq \#P/Q$ is the index of connection of $\Phi$. In some sense $P/Q$ is the ``sandpile group'' in our setting, and in fact we have that $P/Q\simeq \mathrm{coker}(\mathbf{C}^t)$, where~$\mathbf{C}$ is the Cartan matrix of~$\Phi$. Hence, these full-dimensional and saturated components suggest that interval-firing may possibly be connected to Cartan matrix chip-firing. The next remark explains that indeed there is some connection.

\begin{remark} \label{rem:bkr}
Let us explain how Cartan matrix chip-firing (which, as mentioned, has been investigated by Benkart-Klivans-Reiner~\cite{benkart2016chip}) can be realized as a certain ``limit'' of symmetric interval-firing. Note that a Cartan matrix is always an M-matrix (see~\cite[Proposition 4.1]{benkart2016chip}). By associating to each vector $c=(c_1,\ldots,c_n)\in\mathbb{Z}^n$ the weight $\lambda = \sum_{i=1}^{n}c_i\omega_i\in P$, we can view Cartan matrix chip-firing as the relation~$\raU{\mathbf{C}}$ on $P$ defined by
\[ \lambda \raU{\mathbf{C}} \lambda - \alpha_i\; \textrm{ for $\lambda \in P$ and simple root $\alpha_i$, $i\in [n]$ with $\<\lambda,\alpha_i^\vee\>\geq 2$}.\] 
For $\lambda=\sum_{i=1}^{n}c_i\omega_i \in P$ and $k\in\mathbb{Z}_{\geq 0}$ set $B_k(\lambda) \coloneqq  \{\sum_{i=1}^{n}c'_i\omega_i \in P\colon \sum_{i=1}^{n} |c_i-c'_i| \leq k\}$.  In other words, $B_k(\lambda)$ consists of those $\mu$ which are within weight lattice distance $k$ of~$\lambda$. Note that for all $\lambda \in B_k(\rho_k)$, we have that $\<\lambda,\alpha^\vee\>\geq k$ if~$\alpha \in \Phi^+$ is not a simple root. In other words, for $\lambda \in B_k(\rho_k)$, if~$\lambda \raU{\mathrm{sym},k} \lambda + \alpha$, then~$\alpha = \alpha_i$ is some simple root. Moreover, for $\lambda \in B_k(\rho_k)$ we have $\<\lambda,\alpha_i^\vee\>\geq 0$ for any simple root~$\alpha_i$. Hence, for~$\lambda \in B_k(\rho_k)$ the symmetric interval-firing relation reduces to 
\[ \lambda \raU{\mathrm{sym},k} \lambda + \alpha_i\; \textrm{ for a simple root $\alpha_i$, $i\in [n]$ with $\<\lambda,\alpha_i^\vee\>\leq k-1$}.\] 
Define $\Psi_k\colon P\to P$ by $\Psi_k(\lambda) \coloneqq  -\lambda+\rho_{k+1}$ (so $\Psi_k$ is just a ``reflection plus translation''). Then for $\lambda \in \Psi_k^{-1}(B_k(\rho_k)) = B_k(\rho)$ we have
\[ \Psi_k(\lambda) \raU{\mathrm{sym},k} \Psi_k(\lambda - \alpha_i)\; \textrm{ for a simple root $\alpha_i$, $i\in [n]$ with $\<\Psi_k(\lambda),\alpha_i^\vee\>\geq 2$}.\] 
Thus the restriction of $\Psi_{k}^{-1}(\Gamma_{\mathrm{sym},k})$ to $B_k(\rho)$ is exactly the same as the restriction of $\Gamma_{\raU{\mathbf{C}}}$ to $B_k(\rho)$. But every~$\lambda \in P$ belongs to~$B_k(\rho)$ as $k\to \infty$. In this way, we can recover Cartan matrix chip-firing as a certain~$k\to \infty$ limit of symmetric interval-firing.

Benkart-Klivans-Reiner~\cite[Theorem~1.1]{benkart2016chip} showed that the \emph{recurrent configurations} for Cartan matrix chip-firing are $\rho-\omega$ for $\omega \in \Omega_m^0$. Observe~$\Psi_k(\rho-\omega) = \eta_k(\omega)$, so these recurrent configurations correspond exactly to the sinks of our full-dimensional and saturated components. In the same way, the $2^n$ stable configurations in~$\mathbb{Z}_{\geq 0}^n$ for Cartan matrix chip-firing correspond to the sinks of our full-dimensional components. 

We should stress, however, that confluence is much easier to establish for Cartan matrix chip-firing than for our interval-firing processes: for Cartan matrix chip-firing, confluence holds locally, which ultimately has to do with the fact that simple roots are pairwise non-acute. On the other hand, when firing arbitrary positive roots confluence need not hold locally because two positive roots may form an acute angle. Hence while Cartan matrix chip-firing describes the limiting behavior of our interval-firing process, it does not explain why the system is confluent from every initial point. Indeed, we could have also obtained Cartan matrix chip-firing by taking the same $k\to\infty$ limit of the root-firing process which has $\lambda \to \lambda +\alpha$ for $\lambda\in P$, $\alpha\in\Phi^+$ when $\<\lambda,\alpha^\vee\>+ 1 \in [-k+2,k]$, but that process is not confluent.
\end{remark}

\section{Confluence of truncated interval-firing} \label{sec:tr_conf}


So far in this paper we have mostly focused on symmetric interval-firing. We now finally turn to truncated interval-firing. In this section we prove the confluence of~$\raU{\mathrm{tr},\mathbf{k}}$. Let us start by describing the sinks of $\Gamma_{\mathrm{tr},\mathbf{k}}$.

\begin{lemma} \label{lem:trsinks}
For any $\mathbf{k} \in \mathbb{N}[\Phi]^W$, the sinks of $\Gamma_{\mathrm{tr},\mathbf{k}}$ are~$\{\eta_{\mathbf{k}}(\lambda)\colon \lambda \in P\}$.
\end{lemma}
\begin{proof}
Let $\lambda \in P$.  Let $\alpha \in \Phi^{+}$. Note that since $w_{\lambda} \in W^{I^{0}_{\lambda}}$, $w_{\lambda}$ does not have a descent $s_i$ with $I^{0}_{\lambda}$ and thus $w_{\lambda}$ has no inversions in~$\Phi_{I^{0}_{\lambda}}$. Thus if~$\alpha \in w_{\lambda}(\Phi_{I^{0}_{\lambda}})$, then $\<\eta_{\mathbf{k}}(\lambda),\alpha^\vee\> = \<\lambda_{\mathrm{dom}} + \rho_{\mathbf{k}},w_{\lambda}^{-1}(\alpha)^\vee\> \geq \mathbf{k}(\alpha)$, since $w_{\lambda}^{-1}(\alpha) \in \Phi^{+}$. So now consider $\alpha \notin w_{\lambda}(\Phi_{I^{0}_{\lambda}})$. Then $w_{\lambda}^{-1}(\alpha)$ may be positive or negative, but $|\<\lambda_{\mathrm{dom}},w_{\lambda}(\alpha)^\vee\>| \geq 1$ (because $\lambda_{\mathrm{dom}}$ has an $\omega_i$ coefficient of at least~$1$ for some $i \notin I^{0}_{\lambda}$ such that $\alpha_i^\vee$ appears in the expansion of $\pm w_{\lambda}(\alpha)^\vee$). Hence 
\[|\<\eta_{\mathbf{k}}(\lambda),\alpha^\vee\>| = |\<\lambda_{\mathrm{dom}} + \rho_{\mathbf{k}},w_{\lambda}^{-1}(\alpha)^\vee\>| \geq \mathbf{k}(\alpha)+1,\] 
which means that $\<\eta_{\mathbf{k}}(\lambda),\alpha^\vee\> \notin [-\mathbf{k}(\alpha),\mathbf{k}(\alpha)-1]$. So indeed $\eta_{\mathbf{k}}(\lambda)$ is a sink of $\Gamma_{\mathrm{tr},\mathbf{k}}$.

Now suppose $\mu$ is a sink of $\Gamma_{\mathrm{tr},\mathbf{k}}$. Since $\<\mu,\alpha^\vee\> \notin [-\mathbf{k}(\alpha),\mathbf{k}(\alpha)-1]$ for all $\alpha \in \Phi^{+}$, in particular $|\<\mu,\alpha^\vee\>| \geq \mathbf{k}(\alpha)$ for all $\alpha\in \Phi^{+}$. This means that $\<\mu_{\mathrm{dom}},\alpha^\vee\> \geq \mathbf{k}(\alpha)$ for all~$\alpha \in \Phi^{+}$. Hence $\mu_{\mathrm{dom}} = \mu' + \rho_{\mathbf{k}}$ for some dominant $\mu' \in P_{\geq 0}$. Suppose to the contrary that $w_{\mu}$ is not the minimal length element of $w_{\mu}W_{I^{0}_{\mu'}}$. Then there exists a descent $s_i$ of~$w_{\mu}$ with $i \in I^{0}_{\mu'}$. But then
\[\<\mu,-w_{\mu}(\alpha_i)^\vee\> = \<\mu_{\mathrm{dom}},-\alpha_i^\vee\> = -\<\mu',\alpha_i^\vee\> -\<\rho_{\mathbf{k}},\alpha_i^\vee\> \geq -\mathbf{k}(\alpha_i),\] 
and also $\<\mu,-w_{\mu}(\alpha_i)^\vee\> =-\< \mu_{\mathrm{dom}},\alpha_i^\vee\> \leq 0$. This would mean $\mu$ is not a sink of $\Gamma_{\mathrm{tr},\mathbf{k}}$, since $-w_{\mu}(\alpha_i) \in \Phi^{+}$. So $w_{\mu}$ must be the minimal length element of~$w_{\mu}W_{I^{0}_{\mu'}}$. This means $w_{\mu} = w_{\lambda}$ for some~$\lambda \in P$ with $\lambda_{\mathrm{dom}} = \mu'$. And $\mu = w_{\mu}(\mu_{\mathrm{dom}}) = \lambda + w_{\lambda}(\rho_{\mathbf{k}}) = \eta_{\mathbf{k}}(\lambda)$, as claimed.
\end{proof}

We now proceed to prove the confluence of truncated interval-firing. In some sense our proof of confluence here is less satisfactory than the one for symmetric interval-firing because we heavily rely on the diamond lemma, and reduction to rank~$2$, which is a kind of ``trick'' that obscures the underlying polytopal geometry (and requires us at one point to use the classification of rank~$2$ root systems). But we also do crucially use the permutohedron non-escaping lemma in the following lemma, which says that ``small'' permutohedra close to the origin are connected components of truncated interval-firing.

\begin{lemma} \label{lem:trccs}
Let $\mathbf{k} \in \mathbb{N}[\Phi]^W$ be good. Then for all $\omega \in \Omega^m_0$, the (translated) discrete permutohedron $\Pi^Q(\rho_{\mathbf{k}})+\omega$ is a connected component of $\Gamma_{\mathrm{tr},\mathbf{k}}$ and the unique sink of this connected component is $\rho_{\mathbf{k}}+\omega$.
\end{lemma}
\begin{proof}
First let us prove a preliminary result: for any $\lambda \in P$ and $\omega \in \Omega^m_0$, we have that~$(\lambda-\omega)_{\mathrm{dom}} = \lambda_{\mathrm{dom}} - w(w^{-1}_{\lambda}(\omega))$ for some $w \in W_{I^{0}_{\lambda}}$. Indeed, since~$\omega$ is minuscule or zero,  we have that $\<-w'(\omega),\alpha^\vee\> \in \{-1,0,1\}$ for any $\alpha\in \Phi$ and any~$w'\in W$. Therefore~$w^{-1}_{\lambda}(\lambda-\omega) = \lambda_{\mathrm{dom}} - w^{-1}_{\lambda}(\omega)$ may not be dominant, but the only $\alpha_i$ for which we have~$\< \lambda_{\mathrm{dom}} - w^{-1}_{\lambda}(\omega),\alpha_i^\vee\> < 0$ must have $i \in I^{0}_{\lambda}$. Hence, if we let~$w \in W_{I^{0}_{\lambda}}$ be such that $\<w(w^{-1}_{\lambda}(\omega)),\alpha_i^\vee\> \geq 0$ for all $i \in I^{0}_{\lambda}$, then~$(\lambda-\omega)_{\mathrm{dom}} = \lambda_{\mathrm{dom}} - w(w^{-1}_{\lambda}(\omega))$ as claimed.

Now let us show that for any $\omega \in \Omega_m^0$, the only sink of $\Gamma_{\mathrm{tr},\mathbf{k}}$ in $\Pi^Q(\rho_{\mathbf{k}})+\omega$ is $\rho_{\mathbf{k}}+\omega$. Suppose $\eta_{\mathbf{k}}(\lambda) \in \Pi^Q(\rho_{\mathbf{k}})+\omega$ for some $\lambda \in P$. This means $\eta_{\mathbf{k}}(\lambda) -\omega \in \Pi^Q(\rho_{\mathbf{k}})$, which means that $(\eta_{\mathbf{k}}(\lambda)-\omega)_{\mathrm{dom}} = \lambda_{\mathrm{dom}}+\rho_{\mathbf{k}} - w(w_{\lambda}(\omega)) \in \Pi^Q(\rho_{\mathbf{k}})$ for some $w \in W_{I^{0}_{\lambda}}$ (we are using that $w_{\eta_{\mathbf{k}}(\lambda)}=w_{\lambda}$, which we have seen before, and that $W_{I^{0}_{\lambda}} \subseteq W_{I^{0}_{\lambda}}+\rho_{\mathbf{k}}$). Hence Proposition~\ref{prop:perm_containment} tells us that
\[\rho_{\mathbf{k}}-( \lambda_{\mathrm{dom}}+\rho_{\mathbf{k}} - w(w_{\lambda}(\omega))) = -(\lambda_{\mathrm{dom}}-\omega) + ( w(w_{\lambda}(\omega)) - \omega) \in Q_{\geq 0}.\]
Now, since $\lambda_{\mathrm{dom}} \in (Q+\omega)\cap P_{\geq 0}$, we know that $\lambda_{\mathrm{dom}}-\omega \in Q_{\geq 0}$ (by one characterization of minuscule weights mentioned in~\S\ref{sec:rootsystemdefs}). Also, $\omega-w(w_{\lambda}(\omega)) \in Q_{\geq 0}$ by Proposition~\ref{prop:perm_containment}. Hence we conclude that $\lambda_{\mathrm{dom}}=\omega$ and $w(w_{\lambda}(\omega)) = \omega$. But since we have~$w \in W_{I^{0}_{\lambda}}$, we conclude that~$w(w_{\lambda}(\omega))=w_{\lambda}(\omega)$, and thus $w_{\lambda}(\omega)=\omega$, which forces $w_{\lambda}$ to be the identity, i.e., we have $\lambda =\omega$. So indeed the only sink of $\Gamma_{\mathrm{tr},\mathbf{k}}$ in~$\Pi^Q(\rho_{\mathbf{k}})+\omega$ is~$\rho_{\mathbf{k}}+\omega$.

Let us prove the lemma first for $\omega = 0$. Since $\raU{\mathrm{tr},\mathbf{k}}$ is terminating by Proposition~\ref{prop:interval_firing_terminating}, any $\raU{\mathrm{tr},\mathbf{k}}$-firing sequence starting at some $\mu \in \Pi^Q(\rho_{\mathbf{k}})$ has to terminate somewhere. By the permutohedron non-escaping lemma, Lemma~\ref{lem:permtrap}, such a sequence must terminate somewhere inside $\Pi^Q(\rho_{\mathbf{k}})$; and since $\rho_{\mathbf{k}}$ is the only sink in~$\Pi^Q(\rho_{\mathbf{k}})$, it must terminate at $\rho_{\mathbf{k}}$. So indeed $\Pi^Q(\rho_{\mathbf{k}})$ is a connected component of~$\Gamma_{\mathrm{tr},\mathbf{k}}$.

Finally, let $\omega \in \Omega_m$ be arbitrary, and let $w \in C$ be the element corresponding to~$\omega$ under the isomorphism $C \simeq P/Q$. Then by the description of this isomorphism in~\S\ref{sec:symmetry} we get~$w(0-\rho/h) +\rho/h= \omega$, and hence $w(\Pi^Q(\rho_{\mathbf{k}})-\rho/h)+\rho/h = \Pi^Q(\rho_{\mathbf{k}})+\omega$. So from the symmetry of $\Gamma^{\mathrm{un}}_{\mathrm{tr},\mathbf{k}}$ described in Theorem~\ref{thm:symmetry}, we get that $\Pi^Q(\rho_{\mathbf{k}})+\omega$ is also a connected component of $\Gamma_{\mathrm{tr},\mathbf{k}}$.
\end{proof}

Now we consider truncated interval-firing for rank~$2$ root systems.

\begin{prop} \label{prop:rank2tr}
Suppose $\Phi$ is of rank~$2$. Let $\mathbf{k} \in \mathbb{N}[\Phi]^W$. Let $\lambda \in P$ be such that $\<\lambda,\alpha^\vee\> \in [-\mathbf{k}(\alpha),\mathbf{k}(\alpha)]$ and  $\<\lambda,\beta^\vee\> \in [-\mathbf{k}(\beta),\mathbf{k}(\beta)]$ for two linearly independent roots $\alpha,\beta \in \Phi$.  Suppose that either $\Phi$ is simply laced or one of $\alpha$ and $\beta$ is short and the other is long. Let $\omega \in \Omega_m^0$ be such that $\rho_{\mathbf{k}} - \lambda \in Q+\omega$. Then $\lambda \in \Pi^Q(\rho_\mathbf{k}) + \omega$.
\end{prop}
\begin{proof}
First let us show $\lambda_{\mathrm{dom}} = c_1\omega_1 + c_2\omega_2$ with $c_1 \in [0,\mathbf{k}(\alpha_1)]$ and $c_2 \in [0,\mathbf{k}(\alpha_2)]$. Observe that $\<\lambda_{\mathrm{dom}},w_{\lambda}(\alpha)^\vee\> \in [-\mathbf{k}(\alpha),\mathbf{k}(\alpha)]$ and similarly for $\beta$. By replacing $\alpha$ with~$-\alpha$ and $\beta$ with $-\beta$ if necessary, we can assume $\<\lambda_{\mathrm{dom}},w_{\lambda}(\alpha)^\vee\> \in [0,\mathbf{k}(\alpha)]$ and similarly for $\beta$, and since $\lambda_{\mathrm{dom}}$ is dominant, we are free to assume that $w_{\lambda}(\alpha)^\vee$ is positive and similarly for $\beta$. Note that $w_{\lambda}(\alpha)^\vee$ and $w_{\lambda}(\beta)^\vee$ are both nonnegative integer combinations of the simple coroots $\alpha_1^\vee$ and $\alpha_2^\vee$. Then, since $\alpha$ and $\beta$ are linearly independent, and since either $\Phi$ is simply laced, in which case $\mathbf{k}(\alpha)=\mathbf{k}(\beta)=k$, or one of $\alpha,\beta$ is short (say e.g. $\mathbf{k}(\alpha)=k_s$) and the other is long (say e.g. $\mathbf{k}(\beta)=k_l$), we can conclude in fact that $\<\lambda_{\mathrm{dom}},\alpha_1^\vee\> \in [0,\mathbf{k}(\alpha_1)]$ and $\<\lambda_{\mathrm{dom}},\alpha_2^\vee\> \in [0,\mathbf{k}(\alpha_2)]$.

So indeed, $\lambda_{\mathrm{dom}} = c_1\omega_1 + c_2\omega_2$ with $c_1 \in [0,\mathbf{k}(\alpha_1)]$ and $c_2 \in [0,\mathbf{k}(\alpha_2)]$. If $c_1 =\mathbf{k}(\alpha_1)$ and $c_2 = \mathbf{k}(\alpha_2)$, then $\lambda_{\mathrm{dom}} = \rho_{\mathbf{k}}$ and the proposition is obvious in this case (note that we will have~$\omega=0$). So assume without loss of generality that $c_2 \leq \mathbf{k}(\alpha_2)-1$. 

Let $\lambda' \coloneqq  \lambda - \omega$. We want to show $\lambda' \in \Pi^Q(\rho_\mathbf{k})$. As we have seen in the proof of Lemma~\ref{lem:trccs}, we have $\lambda'_{\mathrm{dom}} = \lambda_{\mathrm{dom}} - w(\omega)$ for some $w \in W$. So let $w\in W$ be such that $\lambda'_{\mathrm{dom}} = \lambda_{\mathrm{dom}} - w(\omega)$ and write $\lambda'_{\mathrm{dom}} = c'_1\omega_1 + c'_2\omega_2$. Since~$\<-w(\omega),\alpha^\vee\> \in \{-1,0,1\}$ for any $\alpha\in \Phi$, we have $c'_1 \leq \mathbf{k}(\alpha_1)+1$ and $c'_2 \leq \mathbf{k}(\alpha_2)$. First suppose $c'_1\leq \mathbf{k}(\alpha_1)$. Together with $c'_2\leq \mathbf{k}(\alpha_2)$, this implies that $\rho_{\mathbf{k}}-\lambda'_{\mathrm{dom}} \in P_{\geq 0}$, and hence $\rho_{\mathbf{k}}-\lambda'_{\mathrm{dom}} \in Q_{\geq 0}$. Thus we conclude $\lambda' \in \Pi^Q(\rho_\mathbf{k})$ by Proposition~\ref{prop:perm_containment}. 

So suppose that $c'_1 = \mathbf{k}(\alpha_1)+1$. This means $\<-w(\omega),\alpha_1^\vee\> =1$. Note that this implies $\omega \neq 0$, and hence $\omega$ must be a minuscule weight. But $G_2$ has no minuscule weights, so we may from now on assume that $\Phi \neq G_2$. Since $\omega$ is the only dominant element of $W(\omega)$, we also have that~$\<-w(\omega),\alpha_2^\vee\> \leq 0$ and hence~$c'_2 \leq \mathbf{k}(\alpha_2)-1$. Write $\rho_{\mathbf{k}}-\lambda'_{\mathrm{dom}} = a_1\alpha_1 + a_2\alpha_2$ for some integers $a_1,a_2\in \mathbb{Z}$. Then~$c'_1 = \mathbf{k}(\alpha_1)+1$ and $c'_2 \leq \mathbf{k}(\alpha_2)-1$ translate to
\begin{align*}
2a_1+\<\alpha_2,\alpha^\vee_1\>a_2 &= -1; \\
\<\alpha_1,\alpha^\vee_2\>a_1 + 2a_2 &\geq 1.
\end{align*}
By the classification of rank~$2$ root systems we have $\<\alpha_2,\alpha^\vee_1\>,\<\alpha_1,\alpha^\vee_2\> \in \{-1,-2\}$ with at least one of them equal to $-1$. It is then not hard to check that all integer solutions $a_1,a_2\in \mathbb{Z}$ to the above system of inequalities must have $a_1,a_2\geq 0$. Hence we conclude~$\rho_{\mathbf{k}}-\lambda'_{\mathrm{dom}} \in Q_{\geq 0}$, and thus $\lambda' \in \Pi^Q(\rho_\mathbf{k})$ by Proposition~\ref{prop:perm_containment}.
\end{proof}

\begin{cor} \label{cor:rank2trconfluence}
Suppose $\Phi$ is of rank~$2$. Let $\mathbf{k}\in\mathbb{N}[\Phi]^W$ be good. Then $\raU{\mathrm{tr},\mathbf{k}}$ is confluent (and terminating).
\end{cor}
\begin{proof}
We know $\raU{\mathrm{tr},\mathbf{k}}$ is terminating thanks to Proposition~\ref{prop:interval_firing_terminating}. Hence by the diamond lemma, Lemma~\ref{lem:diamond}, it is enough to prove that $\raU{\mathrm{tr},\mathbf{k}}$ is locally confluent.

First let us prove this when $\Phi$ is simply laced. Suppose $\lambda \raU{\mathrm{tr},\mathbf{k}} \lambda+\alpha$ and $\lambda \raU{\mathrm{tr},\mathbf{k}} \lambda+\beta$ for $\alpha,\beta \in \Phi^{+}$. Then by Proposition~\ref{prop:rank2tr} we have that $\lambda \in \Pi^Q(\rho_\mathbf{k}) + \omega$ where $\omega\in \Omega_m^0$ is such that $\rho_{\mathbf{k}} - \lambda \in Q+\omega$. But by Lemma~\ref{lem:trsinks}, $\Pi^Q(\rho_\mathbf{k}) + \omega$ is a connected component of~$\Gamma_{\mathrm{tr},\mathbf{k}}$ with unique sink $\rho_{\mathbf{k}}+\omega$; since $\raU{\mathrm{tr},\mathbf{k}}$ is terminating this means that any $\raU{\mathrm{tr},\mathbf{k}}$-firing sequence starting at $\lambda$ eventually terminates at $\rho_{\mathbf{k}}+\omega$. Hence we can bring~$\lambda+\alpha$ and~$\lambda+\beta$ back together again via $\raU{\mathrm{tr},\mathbf{k}}$-firings.

Note that confluence for $\Phi=A_1\oplus A_1$ (for any $\mathbf{k} \in\mathbb{N}[\Phi]^W$) reduces to confluence for~$\Phi=A_1$, which is trivial. Thus in fact we have proved confluence for all simply laced root systems of rank~$2$, including those which are not irreducible.

So assume~$\Phi$ is not simply laced. Suppose $\lambda \raU{\mathrm{tr},\mathbf{k}} \lambda+\alpha$ and $\lambda \raU{\mathrm{tr},\mathbf{k}} \lambda+\beta$ for~$\alpha,\beta \in \Phi^{+}$. If one of $\alpha$ and $\beta$ is short and the other is long, then we can apply Proposition~\ref{prop:rank2tr} and Lemma~\ref{lem:trsinks} as above to conclude that we can bring $\lambda+\alpha$ and $\lambda+\beta$ back together again via $\raU{\mathrm{tr},\mathbf{k}}$-firings. So suppose $\alpha$ and $\beta$ have the same length. Then let~$\widetilde{\Phi}$ be the set of all roots in $\Phi$ with the same length as $\alpha$ and $\beta$. This $\widetilde{\Phi}$ will again be a rank~$2$ root system, and by construction a simply laced one. Hence by the result for simply laced root systems, we know that truncated interval-firing is confluent for $\widetilde{\Phi}$; so in particular we can bring $\lambda+\alpha$ and $\lambda+\beta$ back together again via $\raU{\mathrm{tr},\mathbf{k}}$-firings.
\end{proof}

The confluence of truncated interval-firing for all root systems follows easily from confluence for rank~$2$ root systems. The following finishes the proof of Theorem~\ref{thm:confluence_intro}.

\begin{cor} \label{cor:trconfluence}
Let  $\mathbf{k}\in\mathbb{N}[\Phi]^W$ be good. Then $\raU{\mathrm{tr},\mathbf{k}}$ is confluent (and terminating).
\end{cor}
\begin{proof}
We know $\raU{\mathrm{tr},\mathbf{k}}$ is terminating thanks to Proposition~\ref{prop:interval_firing_terminating}. Hence by the diamond lemma, Lemma~\ref{lem:diamond}, it is enough to prove that $\raU{\mathrm{tr},\mathbf{k}}$ is locally confluent. Suppose that~$\lambda \raU{\mathrm{tr},\mathbf{k}} \lambda+\alpha$ and $\lambda \raU{\mathrm{tr},\mathbf{k}} \lambda+\beta$ for $\alpha,\beta \in \Phi^{+}$. Restricting $\Phi$ to the span of $\alpha$ and $\beta$ gives a rank~$2$ sub-root system, for which we have proved confluence in Corollary~\ref{cor:rank2trconfluence} (as remarked in the proof of that corollary, we in fact proved confluence for \emph{all} rank~$2$ root systems, including those which are not irreducible). Hence we can bring $\lambda+\alpha$ and $\lambda+\beta$ back together just with truncated interval-firing moves inside that rank~$2$ sub-root system.
\end{proof}

\begin{remark}
Our method of proof of confluence for $\raU{\mathrm{tr},\mathbf{k}}$ fails when $\mathbf{k}$ is not good; for instance, Lemma~\ref{lem:trccs} is not true for general $\mathbf{k}$, as can be seen in Example~\ref{ex:notgood}: here $0 \in \Pi^Q(\rho_{\mathbf{k}})$ but $0$ does not belong to the connected component of $\Gamma_{\mathrm{tr},\mathbf{k}}$ containing~$\rho_{\mathbf{k}}$. However, we can actually deduce that $\raU{\mathrm{tr},\mathbf{k}}$ is confluent for all $\mathbf{k}\in\mathbb{N}[\Phi]^W$ from Corollary~\ref{cor:trconfluence}. Indeed, if $\mathbf{k}\in\mathbb{N}[\Phi]^W$ is not good, then $k_s=0$. But if~$k_s=0$ then we will never be able to fire any short root. In other words, if $k_s=0$ then truncated interval-firing reduces to truncated interval-firing with respect to the long roots only; and the long roots form a simply laced root system, for which $\raU{\mathrm{tr},\mathbf{k}}$ is known to be confluent from Corollary~\ref{cor:trconfluence}.
\end{remark}

\begin{remark}
It appears that when $\Phi = A_2$ there are no intervals $[a,b]$ for which the relation $\lambda \to \lambda+\alpha$ for $\lambda \in P$, $\alpha\in\Phi^+$ with $\<\lambda,\alpha^\vee\>+ 1 \in [a,b]$ is confluent besides the symmetric and truncated intervals (and this probably would not be too hard to prove). If so, then the same would be true for all irreducible simply laced root systems (except for~$A_1$) because any irreducible root system of rank $3$ or greater contains an~$A_2$ sub-root system. This observation also severely restricts possible intervals defining confluent processes for all root systems, including the non-simply laced ones (although note that central-firing is confluent for $\Phi = B_2$).
\end{remark}

\begin{remark} \label{rem:hyperplanes}
To any root-firing process $\ra$ on $P$ let us associate the hyperplane arrangement which contains the hyperplane $H = \{v\in V\colon \<v,\alpha^\vee\>= c\}$ whenever we have a firing move $\lambda \to \lambda + \alpha$ with $\<\lambda+\frac{\alpha}{2}, \alpha^\vee\>=c$; i.e., we include a hyperplane orthogonal to $\alpha$ at the \emph{midpoint} between $\lambda$ and $\lambda +\alpha$. As mentioned in the introduction, under this correspondence the symmetric and truncated interval-firing processes correspond to the (extended) Catalan and Shi hyperplane arrangements~\cite{postnikov2000deformations, athanasiadis2000deformations}. The confluence of symmetric and truncated interval-firing seems like it might have something to do with the freeness of the Catalan and Shi arrangements.  \emph{Freeness} is a certain deep algebraic property of hyperplane arrangements introduced by Terao~\cite{terao1980arrangements}. Freeness of the (extended) Catalan and Shi hyperplane arrangements of a root system was conjectured by Edelman and Reiner~\cite{edelman1996free} and proven by Yoshinaga~\cite{yoshinaga2004characterization} building on work of Terao~\cite{terao2002multiderivations}. Vic Reiner suggested that we look at other free deformations of Coxeter arrangements as a possible source of other confluent root-firing processes. We found one such process which, experimentally, appears confluent: for $\mathbf{k}\in\mathbb{N}[\Phi]$ consider the relation $\lambda \ra \lambda+\alpha$ for $\lambda \in P$, $\alpha\in\Phi^+$ with $\<\lambda,\alpha^\vee\>+ 1 \in [-\mathbf{k}(\alpha)+1,\mathbf{k}(\alpha)]$ if $\alpha$ is long and $\<\lambda,\alpha^\vee\>+ 1 \in [-\mathbf{k}(\alpha),\mathbf{k}(\alpha)]$ if $\alpha$ is short. In other words, we use either the truncated or symmetric intervals depending on which Weyl group orbit our root lies in. This process corresponds to a \emph{Shi-Catalan} hyperplane arrangement, as studied by Abe and Terao~\cite{abe2011freeness}. Other free variants of Coxeter arrangements include the \emph{ideal subarrangements} of Coxeter arrangements~\cite{abe2016freeness, abe2016free}, but we have not been able to obtain confluent root-firing processes from these. Note that the freeness of the corresponding hyperplane arrangement certainly does not imply confluence of the root-firing process: for instance, reversing the direction of all the arrows for the truncated interval-firing process yields a process which is not confluent but which corresponds to the same Shi hyperplane arrangement. Nevertheless, it would be very interesting to understand the connection between freeness and confluence further.
\end{remark}

\begin{remark}
Under the correspondence between root-firing processes and hyperplane arrangements discussed in Remark~\ref{rem:hyperplanes}, the central-firing process corresponds not to the central Coxeter arrangement, but rather to the affine \emph{Linial arrangement}. The Linial arrangement has many interesting combinatorial properties (see e.g.~\cite{postnikov2000deformations} and~\cite{athanasiadis2000deformations}), but is not free.
\end{remark}

\part{Ehrhart-like polynomials} \label{part:ehrhart}

\section{Ehrhart-like polynomials: introduction}

Continue to fix a root system $\Phi$ in vector space $V$ as in the previous part (and retain all the notation from that part). In this part, we investigate the set of weights with given symmetric or truncated interval-firing stabilization. Thus, for good $\mathbf{k}\in \mathbb{N}[\Phi]^W$, we define the \emph{stabilization maps} $s^{\mathrm{sym}}_{\mathbf{k}}\colon P \to P$ and $s^{\mathrm{tr}}_{\mathbf{k}}\colon P \to P$ by
\begin{align*}
s^{\mathrm{sym}}_{\mathbf{k}}(\mu) = \lambda &\Leftrightarrow \textrm{ the $\raU{\mathrm{sym},\mathbf{k}}$-stabilization of $\mu$ is $\eta_{\mathbf{k}}(\lambda)$}; \\
s^{\mathrm{tr}}_{\mathbf{k}}(\mu) = \lambda &\Leftrightarrow \textrm{ the $\raU{\mathrm{sym},\mathbf{k}}$-stabilization of $\mu$ is $\eta_{\mathbf{k}}(\lambda)$}.
\end{align*}
These functions are well-defined since the symmetric and truncated interval-firing processes are confluent and terminating (Corollaries~\ref{cor:symconfluence} and~\ref{cor:trconfluence}), the stable points of these processes must have the form $\eta_{\mathbf{k}}(\lambda)$ for some $\lambda \in P$ (Lemmas~\ref{lem:symsinks} and~\ref{lem:trsinks}), and the map $\eta_{\mathbf{k}}$ is injective (Proposition~\ref{prop:etafacts}).

Looking at Example~\ref{ex:rank2graphs}, one can see that the set $(s^{\mathrm{sym}}_{\mathbf{k}})^{-1}(\lambda)$ (or $(s^{\mathrm{tr}}_{\mathbf{k}})^{-1}(\lambda)$) of weights with interval-firing stabilization $\eta_{\mathbf{k}}(\lambda)$ looks ``the same'' across all values of~$\mathbf{k}$ except that it gets ``dilated'' as $\mathbf{k}$ is scaled. In analogy with the \emph{Ehrhart polynomial}~\cite{ehrhart1977polynomes} of a convex lattice polytope, which counts the number of lattice points in dilations of the polytope, let us define for all $\lambda \in P$ and all good $\mathbf{k}\in\mathbb{N}[\Phi]^W$ the quantities:
\begin{align*}
L^{\mathrm{sym}}_{\lambda}(\mathbf{k}) &\coloneqq \#(s^{\mathrm{sym}}_{\mathbf{k}})^{-1}(\lambda); \\
L^{\mathrm{tr}}_{\lambda}(\mathbf{k}) &\coloneqq \#(s^{\mathrm{tr}}_{\mathbf{k}})^{-1}(\lambda).
\end{align*}
Our aim is to show that $L^{\mathrm{sym}}_{\lambda}(\mathbf{k})$ and $L^{\mathrm{tr}}_{\lambda}(\mathbf{k})$ are polynomials in~$\mathbf{k}$. By ``polynomial in~$\mathbf{k}$'' we mean that, if $\Phi$ is simply laced, then these $L^{\mathrm{sym}}_{\lambda}(\mathbf{k})$ and $L^{\mathrm{tr}}_{\lambda}(\mathbf{k})$ are single-variable polynomials in $k$, where $\mathbf{k}(\alpha)=k$ for all $\alpha \in \Phi$; and if $\Phi$ is non-simply laced, then they are two-variable polynomials in $k_s$ and $k_l$, where $\mathbf{k}(\alpha)=k_s$ if $\alpha$ is short and $\mathbf{k}(\alpha)=k_l$ if $\alpha$ is long.

We are able to show that the $L^{\mathrm{sym}}_{\lambda}(\mathbf{k})$ are polynomials for all root systems~$\Phi$ (Theorem~\ref{thm:symehrhart}), and we are able to show that the $L^{\mathrm{tr}}_{\lambda}(\mathbf{k})$ are polynomials assuming that~$\Phi$ is simply laced (Theorem~\ref{thm:trehrhart}). In fact, we show that all these polynomials have integer coefficients. Moreover, we conjecture that for all $\Phi$ that these $L^{\mathrm{sym}}_{\lambda}(\mathbf{k})$ and $L^{\mathrm{tr}}_{\lambda}(\mathbf{k})$ are polynomials with \emph{nonnegative} integer coefficients.

We refer to these $L^{\mathrm{sym}}_{\lambda}(\mathbf{k})$ and $L^{\mathrm{tr}}_{\lambda}(\mathbf{k})$ as the \emph{symmetric} and \emph{truncated Ehrhart-like polynomials} because they count the size of some discrete subset of lattice points as that set is somehow ``dilated.'' But it is important to note that the sets $(s^{\mathrm{sym}}_{\mathbf{k}})^{-1}(\lambda)$ and $(s^{\mathrm{tr}}_{\mathbf{k}})^{-1}(\lambda)$ are in general not the set of lattice points of any convex polytope, or indeed, any convex set. This can already be seen in rank~$2$ (see Example~\ref{ex:rank2graphs}). Nevertheless, for some special~$\lambda$ (namely, $\lambda \in \Omega_m^0$) the polynomials $L^{\mathrm{sym}}_{\lambda}(\mathbf{k})$ and $L^{\mathrm{tr}}_{\lambda}(\mathbf{k})$ are (essentially) genuine Ehrhart polynomials; and so we do use Ehrhart theory to prove the polynomiality of~$L^{\mathrm{sym}}_{\lambda}(\mathbf{k})$ and~$L^{\mathrm{tr}}_{\lambda}(\mathbf{k})$. Note that, because they apparently have nonnegative integer coefficients, these polynomials are (as we explain below) most similar to the Ehrhart polynomials \emph{of zonotopes}. 

\section{Symmetric Ehrhart-like polynomials}\label{sec:sym_Ehrhart}

The \emph{Ehrhart polynomial} $L_{\mathcal{P}}(k)$ of a convex lattice polytope $\mathcal{P}$ is a single-variable polynomial in $k$ which satisfies 
\[L_{\mathcal{P}}(k) = \textrm{the number of lattice points in $k\mathcal{P}$ (the $k$th dilate of $\mathcal{P}$)}\]
for all $k \geq 1$. Such polynomials were first investigated by Ehrhart~\cite{ehrhart1977polynomes}, who proved that they exist for all lattice polytopes. A famous result of Stanley~\cite[Example~3.1]{stanley1980decompositions} says that the Ehrhart polynomial of a lattice \emph{zonotope} (i.e., a Minkowski sum of line segments) has nonnegative integer coefficients. A standard way to prove this result is to inductively \emph{pave} the zonotope (see~\cite[\S9.2]{beck2015computing}); this decomposition of a zonotope goes back to Shephard~\cite{shephard1974combinatorial}. In the following theorem we apply this same paving technique to a slightly more general setting: namely, we show that if~$\mathcal{P}$ is any fixed convex lattice polytope, and $\mathcal{Z}$ is a lattice zonotope, then for $k \geq 1$ the number of lattice points in~$\mathcal{P}+k\mathcal{Z}$ is a polynomial with nonnegative integer coefficients in~$k$. Stanley's result corresponds to taking $\mathcal{P}$ to be a point. Although the proof is, as mentioned, standard, we have not found this theorem in the Ehrhart theory literature; and it turns out that this result is just what we need to prove that the symmetric Ehrhart-like polynomials~$L^{\mathrm{sym}}_{\lambda}(\mathbf{k})$ exist.

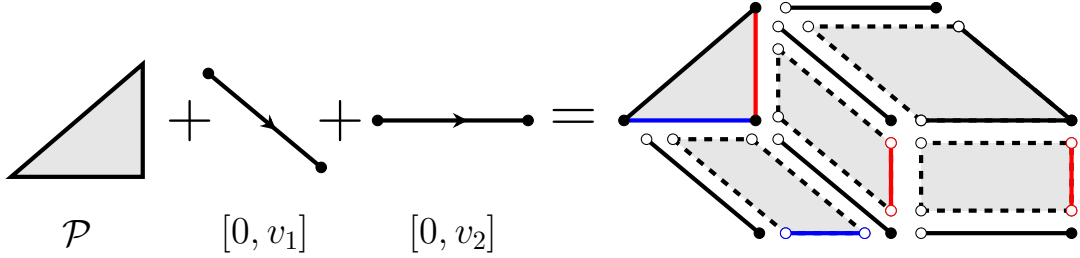
\begin{figure}
\begin{tikzpicture}[>={Stealth[width=2mm,length=2mm,bend]},decoration={markings,mark=at position 0.6 with {\arrow{>}}}]
\def \ax {0}
\def \ay {0}
\def \bx {1.75}
\def \by {0}
\def \cx {1.75}
\def \cy {1.5}
\def \vx {1.5}
\def \vy {-1.25}
\def \wx {2}
\def \wy {0}
\def \spacing {0.2}
\def \circsize {0.07}

\node at (0,0){\begin{tikzpicture}
\draw [ultra thick,fill=gray!20] (\ax,\ay) -- (\bx,\by) -- (\cx,\cy) -- cycle;
\end{tikzpicture}};

\node at (0,-1.5){\Large$\mathcal{P}$};

\node at (1.5,0){\Huge $+$};

\node at (2.5,0){\begin{tikzpicture}
\draw [draw=black,fill=black] (0,0) circle (\circsize);
\draw [draw=black,fill=black] (\vx,\vy) circle (\circsize);
\draw[ultra thick,draw=black,postaction={decorate}]  (0,0) -- (\vx,\vy);
\end{tikzpicture}};

\node at (2.5,-1.5){\Large$[0,v_1]$};

\node at (3.5,0){\Huge $+$};

\node at (5,0){\begin{tikzpicture}
\draw [draw=black,fill=black] (0,0) circle (\circsize);
\draw [draw=black,fill=black] (\wx,\wy) circle (\circsize);
\draw[ultra thick,draw=black,postaction={decorate}]  (0,0) -- (\wx,\wy);
\end{tikzpicture}};

\node at (5,-1.5){\Large$[0,v_2]$};

\node at (6.6,0){\Huge $=$};

\node at (10.25,0){\begin{tikzpicture}
\fill [fill=gray!20] (\ax,\ay) -- (\bx,\by) -- (\cx,\cy) -- cycle;
\draw[ultra thick,draw=blue]  (\ax,\ay) -- (\bx,\by);
\draw[ultra thick,draw=red]  (\bx,\by) -- (\cx,\cy);
\draw[ultra thick,draw=black]  (\cx,\cy) -- (\ax,\ay);

\draw [draw=black,fill=black] (\ax,\ay) circle (\circsize);
\draw [draw=black,fill=black] (\bx,\by) circle (\circsize);
\draw [draw=black,fill=black] (\cx,\cy) circle (\circsize);

\draw [ultra thick,draw=black] (\ax+\spacing*\vx,\ay+\spacing*\vy) -- ({\ax+(1+\spacing)*\vx},{\ay+(1+\spacing)*\vy});
\draw [draw=black,fill=white] (\ax+\spacing*\vx,\ay+\spacing*\vy) circle (\circsize);
\draw [draw=black,fill=black] ({\ax+(1+\spacing)*\vx},{\ay+(1+\spacing)*\vy}) circle (\circsize);

\draw [ultra thick,draw=black] (\bx+\spacing*\vx,\by+\spacing*\vy) -- ({\bx+(1+\spacing)*\vx},{\by+(1+\spacing)*\vy});
\draw [draw=black,fill=white] (\bx+\spacing*\vx,\by+\spacing*\vy) circle (\circsize);
\draw [draw=black,fill=black] ({\bx+(1+\spacing)*\vx},{\by+(1+\spacing)*\vy}) circle (\circsize);

\draw [ultra thick,draw=black] (\cx+\spacing*\vx,\cy+\spacing*\vy) -- ({\cx+(1+\spacing)*\vx},{\cy+(1+\spacing)*\vy});
\draw [draw=black,fill=white] (\cx+\spacing*\vx,\cy+\spacing*\vy) circle (\circsize);
\draw [draw=black,fill=black] ({\cx+(1+\spacing)*\vx},{\cy+(1+\spacing)*\vy}) circle (\circsize);

\draw [ultra thick,dashed,fill=gray!20] ({\ax+\spacing*\vx+\spacing*(\bx-\ax)},{\ay+\spacing*\vy+\spacing*(\by-\ay)}) -- ({\ax+(1+\spacing)*\vx+\spacing*(\bx-\ax)},{\ay+(1+\spacing)*\vy+\spacing*(\by-\ay)}) -- ({\bx+(1+\spacing)*\vx+\spacing*(\ax-\bx)},{\by+(1+\spacing)*\vy+\spacing*(\ay-\by)}) -- ({\bx+\spacing*\vx+\spacing*(\ax-\bx)},{\by+\spacing*\vy+\spacing*(\ay-\by)}) -- cycle;
\draw [draw=black,fill=white] ({\ax+\spacing*\vx+\spacing*(\bx-\ax)},{\ay+\spacing*\vy+\spacing*(\by-\ay)}) circle (\circsize);
\draw [draw=black,fill=white] ({\ax+(1+\spacing)*\vx+\spacing*(\bx-\ax)},{\ay+(1+\spacing)*\vy+\spacing*(\by-\ay)}) circle (\circsize);
\draw [draw=black,fill=white] ({\bx+(1+\spacing)*\vx+\spacing*(\ax-\bx)},{\by+(1+\spacing)*\vy+\spacing*(\ay-\by)}) circle (\circsize);
\draw [draw=black,fill=white] ({\bx+\spacing*\vx+\spacing*(\ax-\bx)},{\by+\spacing*\vy+\spacing*(\ay-\by)}) circle (\circsize);

\draw[ultra thick,draw=blue]  ({\ax+(1+\spacing)*\vx+\spacing*(\bx-\ax)},{\ay+(1+\spacing)*\vy+\spacing*(\by-\ay)}) -- ({\bx+(1+\spacing)*\vx+\spacing*(\ax-\bx)},{\by+(1+\spacing)*\vy+\spacing*(\ay-\by)});
\draw [draw=blue,fill=white] ({\ax+(1+\spacing)*\vx+\spacing*(\bx-\ax)},{\ay+(1+\spacing)*\vy+\spacing*(\by-\ay)}) circle (\circsize);
\draw [draw=blue,fill=white] ({\bx+(1+\spacing)*\vx+\spacing*(\ax-\bx)},{\by+(1+\spacing)*\vy+\spacing*(\ay-\by)}) circle (\circsize);

\draw [ultra thick,dashed,fill=gray!20] ({\cx+\spacing*\vx+\spacing*(\bx-\cx)},{\cy+\spacing*\vy+\spacing*(\by-\cy)}) -- ({\cx+(1+\spacing)*\vx+\spacing*(\bx-\cx)},{\cy+(1+\spacing)*\vy+\spacing*(\by-\cy)}) -- ({\bx+(1+\spacing)*\vx+\spacing*(\cx-\bx)},{\by+(1+\spacing)*\vy+\spacing*(\cy-\by)}) -- ({\bx+\spacing*\vx+\spacing*(\cx-\bx)},{\by+\spacing*\vy+\spacing*(\cy-\by)}) -- cycle;

\draw [draw=black,fill=white] ({\cx+\spacing*\vx+\spacing*(\bx-\cx)},{\cy+\spacing*\vy+\spacing*(\by-\cy)}) circle (\circsize);
\draw [draw=black,fill=white] ({\cx+(1+\spacing)*\vx+\spacing*(\bx-\cx)},{\cy+(1+\spacing)*\vy+\spacing*(\by-\cy)}) circle (\circsize);
\draw [draw=black,fill=white] ({\bx+(1+\spacing)*\vx+\spacing*(\cx-\bx)},{\by+(1+\spacing)*\vy+\spacing*(\cy-\by)}) circle (\circsize);
\draw [draw=black,fill=white] ({\bx+\spacing*\vx+\spacing*(\cx-\bx)},{\by+\spacing*\vy+\spacing*(\cy-\by)}) circle (\circsize);
\draw[ultra thick,draw=red]  ({\cx+(1+\spacing)*\vx+\spacing*(\bx-\cx)},{\cy+(1+\spacing)*\vy+\spacing*(\by-\cy)}) -- ({\bx+(1+\spacing)*\vx+\spacing*(\cx-\bx)},{\by+(1+\spacing)*\vy+\spacing*(\cy-\by)});
\draw [draw=red,fill=white] ({\cx+(1+\spacing)*\vx+\spacing*(\bx-\cx)},{\cy+(1+\spacing)*\vy+\spacing*(\by-\cy)})  circle (\circsize);
\draw [draw=red,fill=white] ({\bx+(1+\spacing)*\vx+\spacing*(\cx-\bx)},{\by+(1+\spacing)*\vy+\spacing*(\cy-\by)}) circle (\circsize);

\draw[ultra thick,draw=black] ({\cx+\spacing*\wx},{\cy+\spacing*\wy})-- ({\cx+(1+\spacing)*\wx},{\cy+(1+\spacing)*\wy});
\draw [draw=black,fill=white] ({\cx+\spacing*\wx},{\cy+\spacing*\wy})  circle (\circsize);
\draw [draw=black,fill=black] ({\cx+(1+\spacing)*\wx},{\cy+(1+\spacing)*\wy})  circle (\circsize);

\draw [ultra thick,dashed,fill=gray!20] ({\cx+(1+\spacing)*\vx+\spacing*\wx},{\cy+(1+\spacing)*\vy+\spacing*\wy}) -- ({\cx+\wx+(1+\spacing)*\vx+\spacing*\wx},{\cy+\wy+(1+\spacing)*\vy+\spacing*\wy}) -- ({\cx+\wx+(0+\spacing)*\vx+\spacing*\wx},{\cy+\wy+(0+\spacing)*\vy+\spacing*\wy})--({\cx+(0+\spacing)*\vx+\spacing*\wx},{\cy++(0+\spacing)*\vy+\spacing*\wy})--cycle;

\draw [draw=black,fill=white] ({\cx+(1+\spacing)*\vx+\spacing*\wx},{\cy+(1+\spacing)*\vy+\spacing*\wy}) circle (\circsize);
\draw [draw=black,fill=white] ({\cx+\wx+(1+\spacing)*\vx+\spacing*\wx},{\cy+\wy+(1+\spacing)*\vy+\spacing*\wy}) circle (\circsize);
\draw [draw=black,fill=white] ({\cx+\wx+(0+\spacing)*\vx+\spacing*\wx},{\cy+\wy+(0+\spacing)*\vy+\spacing*\wy}) circle (\circsize);
\draw [draw=black,fill=white] ({\cx+(0+\spacing)*\vx+\spacing*\wx},{\cy++(0+\spacing)*\vy+\spacing*\wy}) circle (\circsize);

\draw [ultra thick,draw=black] ({\cx+(1+\spacing)*\vx+\spacing*\wx},{\cy+(1+\spacing)*\vy+\spacing*\wy}) -- ({\cx+\wx+(1+\spacing)*\vx+\spacing*\wx},{\cy+\wy+(1+\spacing)*\vy+\spacing*\wy});
\draw [ultra thick,draw=black] ({\cx+\wx+(1+\spacing)*\vx+\spacing*\wx},{\cy+\wy+(1+\spacing)*\vy+\spacing*\wy})-- ({\cx+\wx+(0+\spacing)*\vx+\spacing*\wx},{\cy+\wy+(0+\spacing)*\vy+\spacing*\wy});
\draw [draw=black,fill=white] ({\cx+(1+\spacing)*\vx+\spacing*\wx},{\cy+(1+\spacing)*\vy+\spacing*\wy})  circle (\circsize);
\draw [draw=black,fill=black] ({\cx+\wx+(1+\spacing)*\vx+\spacing*\wx},{\cy+\wy+(1+\spacing)*\vy+\spacing*\wy}) circle (\circsize);
\draw [draw=black,fill=white] ({\cx+\wx+(0+\spacing)*\vx+\spacing*\wx},{\cy+\wy+(0+\spacing)*\vy+\spacing*\wy}) circle (\circsize);

\draw[ultra thick,draw=black] ({\bx+(1+\spacing)*\vx+\spacing*\wx},{\by+(1+\spacing)*\vy+\spacing*\wy}) -- ({\bx+\wx+(1+\spacing)*\vx+\spacing*\wx},{\by+\wy+(1+\spacing)*\vy+\spacing*\wy});
\draw [draw=black,fill=white] ({\bx+(1+\spacing)*\vx+\spacing*\wx},{\by+(1+\spacing)*\vy+\spacing*\wy})  circle (\circsize);
\draw [draw=black,fill=black] ({\bx+\wx+(1+\spacing)*\vx+\spacing*\wx},{\by+\wy+(1+\spacing)*\vy+\spacing*\wy})  circle (\circsize);

\draw [ultra thick,dashed,fill=gray!20] ({\bx+(1+\spacing)*\vx+\spacing*\wx+\spacing*(\cx-\bx)},{\by+(1+\spacing)*\vy+\spacing*\wy+\spacing*(\cy-\by)}) -- ({\bx+\wx+(1+\spacing)*\vx+\spacing*\wx+\spacing*(\cx-\bx)},{\by+\wy+(1+\spacing)*\vy+\spacing*\wy+\spacing*(\cy-\by)}) -- ({\cx+\wx+(1+\spacing)*\vx+\spacing*\wx+\spacing*(\bx-\cx)},{\cy+\wy+(1+\spacing)*\vy+\spacing*\wy+\spacing*(\by-\cy)}) --  ({\cx+(1+\spacing)*\vx+\spacing*\wx+\spacing*(\bx-\cx)},{\cy+(1+\spacing)*\vy+\spacing*\wy+\spacing*(\by-\cy)}) -- cycle;

\draw [draw=black,fill=white] ({\bx+(1+\spacing)*\vx+\spacing*\wx+\spacing*(\cx-\bx)},{\by+(1+\spacing)*\vy+\spacing*\wy+\spacing*(\cy-\by)}) circle (\circsize);
\draw [draw=black,fill=white] ({\bx+\wx+(1+\spacing)*\vx+\spacing*\wx+\spacing*(\cx-\bx)},{\by+\wy+(1+\spacing)*\vy+\spacing*\wy+\spacing*(\cy-\by)}) circle (\circsize);
\draw [draw=black,fill=white] ({\cx+\wx+(1+\spacing)*\vx+\spacing*\wx+\spacing*(\bx-\cx)},{\cy+\wy+(1+\spacing)*\vy+\spacing*\wy+\spacing*(\by-\cy)}) circle (\circsize);
\draw [draw=black,fill=white] ({\cx+(1+\spacing)*\vx+\spacing*\wx+\spacing*(\bx-\cx)},{\cy+(1+\spacing)*\vy+\spacing*\wy+\spacing*(\by-\cy)}) circle (\circsize);

\draw[ultra thick,draw=red]  ({\bx+\wx+(1+\spacing)*\vx+\spacing*\wx+\spacing*(\cx-\bx)},{\by+\wy+(1+\spacing)*\vy+\spacing*\wy+\spacing*(\cy-\by)}) -- ({\cx+\wx+(1+\spacing)*\vx+\spacing*\wx+\spacing*(\bx-\cx)},{\cy+\wy+(1+\spacing)*\vy+\spacing*\wy+\spacing*(\by-\cy)});

\draw [draw=red,fill=white] ({\bx+\wx+(1+\spacing)*\vx+\spacing*\wx+\spacing*(\cx-\bx)},{\by+\wy+(1+\spacing)*\vy+\spacing*\wy+\spacing*(\cy-\by)}) circle (\circsize);
\draw [draw=red,fill=white] ({\cx+\wx+(1+\spacing)*\vx+\spacing*\wx+\spacing*(\bx-\cx)},{\cy+\wy+(1+\spacing)*\vy+\spacing*\wy+\spacing*(\by-\cy)}) circle (\circsize);
\end{tikzpicture}};

\end{tikzpicture}
\caption[Paving a polytope plus zonotope]{Paving a polytope plus zonotope.} \label{fig:paving}
\end{figure}

\begin{thm} \label{thm:polypluszoneehrhart}
Let $\Lambda$ be a lattice in $V$. Let $\mathcal{P}$ be any convex lattice polytope in~$V$. Let~$v_1,\ldots,v_m \in \Lambda$ be lattice elements. Then for any~$\mathbf{k}=(k_1,\ldots,k_m) \in \mathbb{Z}_{\geq 0}^m$ the quantity
\[\#(\mathcal{P}+k_1[0,v_1] + \cdots + k_m[0,v_m]) \cap \Lambda\]
is given by a polynomial in the $k_1,\ldots,k_m$ with nonnegative integer coefficients.
\end{thm}
\begin{proof}
For $X = \{u_1,\ldots,u_\ell\}\subseteq V$ linearly independent, a \emph{half-open parallelepiped} with edge set~$X$ is a convex set $\mathcal{Z}^{h.o.}_{X}$ of the form
\[\mathcal{Z}^{h.o.}_{X} = \sum_{i=1}^{\ell} \begin{cases}[0,u_i) &\textrm{if $\varepsilon = 1$}; \\ (0,u_i] &\textrm{if $\varepsilon = -1$}, \end{cases}\]
for some choice of sign vector $(\varepsilon_1,\ldots,\varepsilon_\ell)\in \{-1,1\}^{\ell}$. For $X\subseteq \{v_1,\ldots,v_m\}$ let us use~$\mathbf{k}X \coloneqq  \{k_iv_i\colon v_i \in X\}$. 

The key idea for this theorem: $\mathcal{P}+k_1[0,v_1] + \cdots + k_m[0,v_m]$ can be inductively decomposed (or ``paved'') into disjoint pieces that are (up to translation) of the form
\[ F +\mathcal{Z}^{h.o.}_{\mathbf{k}X},\]
where $X\subseteq  \{v_1,\ldots,v_m\}$ is linearly independent and $F$ is an open face of the polytope~$\mathcal{P}$ which is affinely independent from $\Span_{\mathbb{R}}(X)$. Figure~\ref{fig:paving} shows how this is done. Here by ``open face'' of $\mathcal{P}$ we mean a face minus its relative boundary. Note that vertices have empty relative boundary and hence vertices are open faces. (But observe that~Figure~\ref{fig:paving} is slightly misleading in that we should technically show the whole polytope $\mathcal{P}$ decomposed into its open faces as well; instead the figure shows these pieces grouped into a single bigger piece.) The proof, by induction on $m$, that this is possible works in exactly the same way as for paving a zonotope (see~\cite[Lemma 9.1]{beck2015computing}), so we do not go into the details. Then note that
\[ \#\left( \left(F + \mathcal{Z}^{h.o.}_{\mathbf{k}X} \right) \cap \Lambda \right) = \#\left( \left(F + \mathcal{Z}^{h.o.}_{X} \right) \cap \Lambda \right) \cdot \prod_{v_i\in X} k_i\]
precisely because $F$ is affinely independent from $\mathcal{Z}^{h.o.}_{\mathbf{k}X}$. Hence the desired polynomial in $k_1,\ldots,k_m$ indeed exists: it is a sum over the pieces of this decomposition of~$\#\left( \left(F + \mathcal{Z}^{h.o.}_{X} \right) \cap \Lambda \right) \cdot \prod_{v_i\in X} k_i$. (We are implicitly using the fact that this decomposition can be realized in a uniform way across all values of~$\mathbf{k}$).
\end{proof}

\begin{cor} \label{cor:permehrhart}
For any~$\lambda \in P_{\geq 0}$, for all $\mathbf{k} \in \mathbb{N}[\Phi]^W$ the quantity $\#\Pi^Q(\lambda+\rho_{\mathbf{k}})$ is given by a polynomial with nonnegative integer coefficients in~$\mathbf{k}$.
\end{cor}
\begin{proof}
We are free to translate $\Pi^Q(\lambda+\rho_{\mathbf{k}})$ so that it contains the origin; i.e., clearly $\#\Pi^Q(\lambda+\rho_{\mathbf{k}})$ is the number of $Q$-points in $\Pi(\lambda+\rho_{\mathbf{k}})-\lambda-\rho_{\mathbf{k}}$. One easy consequence of Proposition~\ref{prop:perm_containment} is that $\Pi(\lambda+\mu)=\Pi(\lambda)+\Pi(\mu)$ for dominant weights $\lambda,\mu\in P_{\geq 0}$. Hence, because~$\lambda$ is dominant, we have
\[\Pi(\lambda+\rho_{\mathbf{k}}) -\lambda - \rho_{\mathbf{k}} = (\Pi(\lambda)  -\lambda) + (\Pi(\rho_{\mathbf{k}}) -\rho_{\mathbf{k}}).\]
It is well known that the \emph{regular permutohedron} $\Pi(\rho)$ is a zonotope. In Type~A, a standard way to prove this fact is to compute the Newton polytope of the Vandermonde determinant in two ways (see~\cite[Theorem 9.4]{beck2015computing}). The same argument, but with \emph{Weyl's denominator formula} (see~\cite[\S24.3]{humphreys1972lie}) in place of the Vandermonde determinant, establishes that $\Pi(\rho)=\sum_{\alpha\in\Phi^+}[-\alpha/2,\alpha/2]$. It is then a simple exercise to show that $\Pi(\rho_{\mathbf{k}})=\sum_{\alpha\in\Phi^+}\mathbf{k}(\alpha)[-\alpha/2,\alpha/2]$. Hence,
\[\Pi(\lambda+\rho_{\mathbf{k}}) -\lambda - \rho_{\mathbf{k}} = (\Pi(\lambda)  -\lambda) + \sum_{\alpha \in \Phi^+} \mathbf{k}(\alpha)[0,-\alpha],\]
and so the desired polynomial indeed exists thanks to Theorem~\ref{thm:polypluszoneehrhart}.
\end{proof}

We are now ready to prove the first part of Theorem~\ref{thm:Ehrhart_intro}.
\begin{thm} \label{thm:symehrhart}
For any $\lambda \in P$, for good~$\mathbf{k} \in\mathbb{N}[\Phi]^W$  the quantity $L^{\mathrm{sym}}_{\lambda}(\mathbf{k})$ is given by a polynomial with integer coefficients in~$\mathbf{k}$.
\end{thm}
\begin{proof}
First of all, if $\lambda$ has $\<\lambda,\alpha^\vee\> = -1$ for some $\alpha \in \Phi^{+}$ then clearly we can take $L^{\mathrm{sym}}_{\lambda}(\mathbf{k})\coloneqq 0$ because, thanks to Lemma~\ref{lem:symsinks}, $\eta_{\mathbf{k}}(\lambda)$ cannot be a sink of~$\Gamma_{\mathrm{sym},\mathbf{k}}$ in this case. So now assume that $\lambda$ satisfies $\<\lambda,\alpha^\vee\>\neq -1$ for all $\alpha \in \Phi^{+}$. If~$I^{0,1}_{\lambda} \neq [n]$, then, by Theorem~\ref{thm:permcc}, the connected component of $\Gamma_{\mathrm{sym},\mathbf{k}}$ containing the sink $\eta_{\mathbf{k}}(\lambda)$ is contained in $w_\lambda \Pi^{Q}_{I^{0,1}_{\lambda}}(\lambda_{\mathrm{dom}})$, which is contained in an affine translate of the strict subspace $\Span_{\mathbb{R}}(w_{\lambda} \Phi_{I^{0,1}_{\lambda}})$. By induction on rank we know the theorem is true for the sub-root system $w_{\lambda} \Phi_{I^{0,1}_{\lambda}}$. Hence, the desired polynomial $L^{\mathrm{sym}}_{\lambda}(\mathbf{k})$ is just the corresponding polynomial for the orthogonal projection of $\lambda$ onto $\Span_{\mathbb{R}}(w_{\lambda} \Phi_{I^{0,1}_{\lambda}})$. (Here we use the fact that the map $\eta_{\mathbf k}$ respects this projection: but this is clear because the projection of $\lambda$ and the projection of $w_{\lambda}(\rho_{\mathbf{k}})$ onto the weight lattice of~$w_{\lambda} \Phi_{I^{0,1}_{\lambda}}$ are both dominant with respect to the choice of $w_{\lambda} \Phi^{+}_{I^{0,1}_{\lambda}}$ as positive roots, which is a subset of~$\Phi^{+}$ by Proposition~\ref{prop:posparaboliccosets}.)

So now assume~$I^{0,1}_{\lambda} = [n]$. This means that $\lambda$ is dominant. Let $\mathbf{k} \in \mathbb{N}[\Phi]^W$ be good. Set~$S \coloneqq  \{\mu\in P\colon \textrm{$\<\mu,\alpha^\vee\>\neq -1$ for all $\alpha \in \Phi^{+}$}, \, \eta_{\mathbf{k}}(\mu) \in \Pi^Q(\eta_{\mathbf{k}}(\lambda))\}$; i.e.,~$S$ is the set of all labels of sinks of $\Gamma_{\mathrm{sym},\mathbf{k}}$ that are inside of $\Pi^Q(\eta_{\mathbf{k}}(\lambda))$. 

We claim that in fact~$S = \{\mu\in P\colon \textrm{$\<\mu,\alpha^\vee\>\neq -1$ for all $\alpha \in \Phi^{+}$}, \, \mu \in \Pi^Q(\lambda)\}$. Indeed, for~$\mu\in P$ with $\<\mu,\alpha^\vee\>\neq -1$ for all $\alpha \in \Phi^{+}$, we have~$ \eta_{\mathbf{k}}(\mu) \in \Pi^Q(\eta_{\mathbf{k}}(\lambda))$ if and only if~$ \eta_{\mathbf{k}}(\mu)_{\mathrm{dom}} = \eta_{\mathbf{k}}(\mu_{\mathrm{dom}}) \in \Pi^Q(\eta_{\mathbf{k}}(\lambda))$. By Proposition~\ref{prop:perm_containment}, we have that~$\eta_{\mathbf{k}}(\mu_{\mathrm{dom}}) \in \Pi^Q(\eta_{\mathbf{k}}(\lambda))$ if and only if $(\lambda+\rho_{\mathbf{k}}) - (\mu_{\mathrm{dom}}+\rho_{\mathbf{k}}) = \lambda-\mu_{\mathrm{dom}} \in Q_{\geq 0}$, which, again by Proposition~\ref{prop:perm_containment}, is if and only if $\mu_{\mathrm{dom}} \in \Pi^Q(\lambda)$, that is, if and only if~$\mu \in \Pi^Q(\lambda)$. Note that this second description of $S$ is independent of $\mathbf{k}$. Also note that for all $\mu\neq \lambda \in S$, either $I^{0,1}_{\mu}\neq [n]$ or $\mu=\mu_{\mathrm{dom}}$, and in the latter case we have that~$\mu$ is strictly less than $\lambda$ in root order. Now, the permutohedron non-escaping lemma, Lemma~\ref{lem:permtrap}, says that
\[\Pi^Q(\eta_{\mathbf{k}}(\lambda)) = \bigcup_{\mu \in S} (s^{\mathrm{sym}}_{\mathbf{k}})^{-1}(\mu).\]
Hence, rewriting, and taking cardinalities, we get
\[ \#(s^{\mathrm{sym}}_{\mathbf{k}})^{-1}(\lambda) = \#\Pi^Q(\eta_{\mathbf{k}}(\lambda)) - \sum_{\mu \neq \lambda \in S} \#(s^{\mathrm{sym}}_{\mathbf{k}})^{-1}(\mu).\]
The quantity $\#\Pi^Q(\eta_{\mathbf{k}}(\lambda))$ is a polynomial in~$\mathbf{k}$ with integer coefficients thanks to Theorem~\ref{cor:permehrhart}. The quantity $\sum_{\mu \neq \lambda \in S} \#(s^{\mathrm{sym}}_{\mathbf{k}})^{-1}(\mu)$ is a polynomial  in~$\mathbf{k}$ with integer coefficients by induction on rank and on root order.  Since the above equality holds for all good $\mathbf{k} \in \mathbb{N}[\Phi]^W$, we conclude that$L^{\mathrm{sym}}_{\lambda}(\mathbf{k})= \#(s^{\mathrm{sym}}_{\mathbf{k}})^{-1}(\lambda)$ is indeed a polynomial in~$\mathbf{k}$ with integer coefficients.
\end{proof}

\begin{table}
\begin{center}
\begin{tabular}{c | c | r }
$\Phi$ & $\lambda$ & $L^{\mathrm{sym}}_{\lambda}(\mathbf{k})$ \\ \specialrule{2.5pt}{1pt}{1pt}
$A_2$ & $0$ & $3k^2 + 3k+1$ \\ \hline
$A_2$ & $\omega_1$ & $3k^2 + 6k + 3$ \\ \hline
$A_2$ & $\omega_2$ & $3k^2 + 6k + 3$ \\ \hline
$A_2$ & $\omega_1+\omega_2$ & $6k+6$ \\ \specialrule{2.5pt}{1pt}{1pt}
$B_2$ & $0$ & $2k_{l}^2+4k_lk_s + k_s^2 + 2k_l + 2k_s + 1$ \\ \hline
$B_2$ & $\omega_1$ & $4k_l+4k_s +4$ \\ \hline
$B_2$ & $\omega_2$ & $2k_l^2 + 4k_lk_s + k_s^2 + 6k_l + 4k_s + 4$ \\ \hline
$B_2$ & $\omega_1+\omega_2$ & $4k_l+4k_s+8$ \\ \specialrule{2.5pt}{1pt}{1pt}
$G_2$ & $0$ & $9k_l^2+12k_lk_s+3k_s^2+3k_l+3k_s+1$ \\ \hline
$G_2$ & $\omega_1$ & $12k_l+6k_s+6$ \\ \hline
$G_2$ & $\omega_2$ & $6k_l+6k_s+6$ \\ \hline
$G_2$ & $\omega_1+\omega_2$ & $6k_l+6k_s+12$ \\ \specialrule{2.5pt}{1pt}{1pt}
\end{tabular}
\end{center}
\caption[Symmetric Ehrhart-like polynomials for rank~$2$ root systems]{The polynomials $L^{\mathrm{sym}}_{\lambda}(\mathbf{k})$ for the irreducible rank $2$ root systems.} \label{tab:sympolys}
\end{table}

Table~\ref{tab:sympolys} records the polynomials $L^{\mathrm{sym}}_{\lambda}(\mathbf{k})$ for the irreducible rank $2$ root systems, for all $\lambda \in P_{\geq 0}$ with $I^{0,1}_{\lambda}=[n]$. Compare these polynomials to the graphs of the corresponding symmetric interval-firing processes in Example~\ref{ex:rank2graphs}.

\begin{remark}
The evaluation of the polynomial $L^{\mathrm{sym}}_{\lambda}(\mathbf{k})$ for $\mathbf{k} \in \mathbb{N}[\Phi]^W$ not good may not count the number of weights in the connected component of $\Gamma_{\mathrm{sym},\mathbf{k}}$ containing~$\eta_{\mathbf{k}}(\lambda)$. For example, take $\Phi = B_2$ and $\mathbf{k}$ defined by $k_s \coloneqq  0$ and $k_l \coloneqq  1$, as in Example~\ref{ex:notgood}. Then, with $\lambda\coloneqq 0$, looking at Table~\ref{tab:sympolys} we see
\[L^{\mathrm{sym}}_{\lambda}(\mathbf{k})=2k_{l}^2+4k_lk_s + k_s^2 + 2k_l + 2k_s + 1=5,\]  
while there are only four weights in the connected component of $\Gamma_{\mathrm{sym},\mathbf{k}}$ containing the sink~$\eta_{\mathbf{k}}(\lambda)$. (Here the ``missing'' weight is of course the origin.)
\end{remark}

\begin{conj} \label{conj:symehrhart}
The polynomials $L^{\mathrm{sym}}_{\lambda}(\mathbf{k})$ have nonnegative integer coefficients.
\end{conj}

When $\lambda \in \Omega_m^{0}$, we know thanks to Proposition~\ref{prop:saturatedccs} that $(s^{\mathrm{sym}}_{\mathbf{k}})^{-1}(\lambda) = \Pi^Q(\lambda+\rho_{\mathbf{k}})$, so Corollary~\ref{cor:permehrhart} implies that Conjecture~\ref{conj:symehrhart} is true in this case. Very recently, the second and fourth authors have proved Conjecture~\ref{conj:symehrhart} in general~\cite{hopkins2018positive}. The first step in their proof of positivity is to give a more refined version of Theorem~\ref{thm:polypluszoneehrhart} that gives an explicit formula for the number of lattice points in a polytope plus dilating zonotope.

\section{Cubical subcomplexes}

In order to proceed further in our investigation of the stabilization maps $s^{\mathrm{sym}}_{\mathbf{k}}$ and~$s^{\mathrm{tr}}_{\mathbf{k}}$, and the relation between them, we need to understand a bit more about the connected components of~$\Gamma_{\mathrm{sym},\mathbf{k}}$. We know that the connected component of $\Gamma_{\mathrm{sym},\mathbf{k}}$ containing the sink~$\eta_{\mathbf{k}}(\lambda)$ is contained in the discrete permutohedron $w_{\lambda}\Pi^Q_{I^{0,1}_{\lambda}}(\lambda_{\mathrm{dom}})$ (Theorem~\ref{thm:permcc}); but it can sometimes contain all of this permutohedron (see Proposition~\ref{prop:saturatedccs}) and can sometimes contain relatively little of it. In this section we will show that there is a small amount of $w_{\lambda}\Pi^Q_{I^{0,1}_{\lambda}}(\lambda_{\mathrm{dom}})$ that this connected component must always contain.

The permutohedron $\Pi_I(\lambda)$ has the structure of a polyhedral complex. The  \emph{cubical subcomplex} of $\Pi_I(\lambda)$ is the union of all faces of $\Pi_I(\lambda)$ that are cubes; here a \emph{cube} means a product of pairwise orthogonal intervals. We denote the cubical subcomplex by~$^\square\Pi_I(\lambda)$. Note that every edge is a cube, and hence $^\square\Pi_I(\lambda)$ contains at least the $1$-skeleton of $\Pi_I(\lambda)$, but it may contain more. We use $^\square\Pi^Q_I(\lambda) \coloneqq  \, ^\square\Pi_I(\lambda) \cap (Q+\lambda)$. 

\begin{prop} \label{prop:symccscubes}
Let $\lambda \in P$ with $\<\lambda,\alpha^\vee\> \neq -1$ for all $\alpha \in \Phi^{+}$ and let $\mathbf{k} \in \mathbb{N}[\Phi]^W$.  Let $Y_{\lambda}\coloneqq \{\mu\in P\colon \mu \baAstU{\mathrm{sym},\mathbf{k}} \eta_{\mathbf{k}}(\lambda) \}$ be the connected component of $\Gamma_{\mathrm{sym},\mathbf{k}}$ containing the sink~$\eta_{\mathbf{k}}(\lambda)$. Then $Y_{\lambda}$ contains the discrete cubical subcomplex~$w_{\lambda} \, ^\square\Pi^Q_{I^{0,1}_{\lambda}}(\eta_\mathbf{k}(\lambda_{\mathrm{dom}}))$.
\end{prop}
\begin{proof}
By the usual projection argument that we have by now carried out many times, we can assume that $I^{0,1}_{\lambda}=[n]$ and consequently that $\lambda$ is dominant. 

For any simple root~$\alpha_i$ we have that $\<\eta_{\mathbf{k}}(\lambda),\alpha_i^\vee\> \in \{\mathbf{k}(\alpha),\mathbf{k}(\alpha)+1\}$. This means that we can ``unfire'' $\alpha_i$ from $\eta_{\mathbf{k}}(\lambda)$; that is, $\<\eta_{\mathbf{k}}(\lambda)-\alpha_i,\alpha_i^\vee\> \leq \mathbf{k}(\alpha)-1$, so that there will be an edge $\eta_{\mathbf{k}}(\lambda)-\alpha_i \raU{\mathrm{sym},\mathbf{k}} \eta_{\mathbf{k}}(\lambda)$ of $\Gamma_{\mathrm{sym},\mathbf{k}}$. In fact, we can keep ``unfiring'' the simple root~$\alpha_i$ until we reach~$s_{\alpha_i}(\eta_{\mathbf{k}}(\lambda))$; i.e., in $\Gamma_{\mathrm{sym},\mathbf{k}}$ there are sequence of edges
\[s_{\alpha_i}(\eta_{\mathbf{k}}(\lambda)) \raU{\mathrm{sym},\mathbf{k}} s_{\alpha_i}(\eta_{\mathbf{k}}(\lambda)) + \alpha_i  \raU{\mathrm{sym},\mathbf{k}} \cdots  \raU{\mathrm{sym},\mathbf{k}} \eta_{\mathbf{k}}(\lambda)-\alpha_i \raU{\mathrm{sym},\mathbf{k}} \eta_{\mathbf{k}}(\lambda).\]
(Note that it is possible that $s_{\alpha_i}(\eta_{\mathbf{k}}(\lambda))=\eta_{\mathbf{k}}(\lambda)$, in which case we would not actually be able to unfire $\alpha_i$ at all). This means that all the $(Q+\eta_{\mathbf{k}}(\lambda))$-points of the entire edge of $\Pi(\eta_{\mathbf{k}}(\lambda))$ between $\eta_{\mathbf{k}}(\lambda)$ and $s_{\alpha_i}(\eta_{\mathbf{k}}(\lambda))$ are reachable via unfirings from $\eta_{\mathbf{k}}(\lambda)$. Moreover, if $\alpha_i$ and $\alpha_j$ are orthogonal, then unfiring one of these does not affect our ability to unfire the other, and hence in this way we can reach any $(Q+\eta_{\mathbf{k}}(\lambda))$-point on a face of $\Pi(\eta_{\mathbf{k}}(\lambda))$ that is the orthogonal product of edges coming out of the vertex~$\eta_{\mathbf{k}}(\lambda)$ in the direction of a negative simple root. Since in particular $s_i(\eta_{\mathbf{k}}(\lambda))$ is reachable via firings and unfirings from $\eta_{\mathbf{k}}(\lambda)$, by applying the $W$-symmetry of~$\Gamma^{\mathrm{un}}_{\mathrm{sym},\mathbf{k}}$ (Theorem~\ref{thm:symmetry}) we see that all vertices of $\Pi(\eta_{\mathbf{k}}(\lambda))$ are so reachable. But note that any face of $\Pi(\eta_{\mathbf{k}}(\lambda))$ can be transported via $W$ to a face containing $\eta_{\mathbf{k}}(\lambda)$, such that the edges of this face which contain $\eta_{\mathbf{k}}(\lambda)$ are in the direction of a negative simple root (see the proof of Theorem~\ref{thm:traverseformula}). We thus conclude that we can reach any $(Q+\eta_{\mathbf{k}}(\lambda))$-point on any cubical face of $\Pi(\eta_{\mathbf{k}}(\lambda))$ via firings and unfirings from~$\eta_{\mathbf{k}}(\lambda)$. 
\end{proof}

\begin{cor} \label{cor:symccsweylorbit}
Let $\lambda \in P$ with $\<\lambda,\alpha^\vee\> \neq -1$ for all $\alpha \in \Phi^{+}$ and let $\mathbf{k} \in \mathbb{N}[\Phi]^W$.  Let $Y_{\lambda}\coloneqq \{\mu\in P\colon \mu \baAstU{\mathrm{sym},\mathbf{k}} \eta_{\mathbf{k}}(\lambda) \}$ be the connected component of $\Gamma_{\mathrm{sym},\mathbf{k}}$ containing the sink~$\eta_{\mathbf{k}}(\lambda)$. Then $Y_{\lambda}$ contains $w_{\lambda}W_{I^{0,1}_{\lambda}}(\eta_\mathbf{k}(\lambda_{\mathrm{dom}}))$. In the special case $\mathbf{k}=0$, $Y_{\lambda}$ is in fact  equal to $w_{\lambda}W_{I^{0,1}_{\lambda}}(\lambda_{\mathrm{dom}})$. 
\end{cor}
\begin{proof}
Note that $^\square\Pi_I(\mu)$ contains at least the $1$-skeleton of $\Pi_I(\mu)$. Thus $Y_{\lambda}$ contains $w_{\lambda}W_{I^{0,1}_{\lambda}} (\eta_\mathbf{k}(\lambda_{\mathrm{dom}}))$ by Proposition~\ref{prop:symccscubes}. Now suppose $\mathbf{k}=0$. If $\mu \raU{\mathrm{sym},0} \mu'$ then $\mu' = \mu + \alpha$ for some $\alpha \in \Phi^{+}$ with~$\<\mu,\alpha^\vee\> = -1$, which means that $\mu' = s_{\alpha}(\mu)$. Hence any two elements in a connected component of $\Gamma_{\mathrm{sym},0}$ must be related by a Weyl group element. By Corollary~\ref{cor:symconfluence}, each connected component of $\Gamma_{\mathrm{sym},0}$ contains only a single sink, and thus the component $Y_{\lambda}$ must be exactly $w_{\lambda}W_{I^{0,1}_{\lambda}}(\lambda_{\mathrm{dom}})$.
\end{proof}

\section{How interval-firing components decompose} \label{sec:decompose}

In this section, we study how symmetric and truncated interval-firing components ``decompose'' into smaller components. Let us explain what we mean by ``decompose'' more precisely. For any $\mathbf{k} \in \mathbb{N}[\Phi]^W$, $\Gamma_{\mathrm{tr},{\mathbf{k}}}$ is a subgraph of $\Gamma_{\mathrm{sym},{\mathbf{k}}}$, so the connected components of $\Gamma_{\mathrm{sym},{\mathbf{k}}}$ are unions of connected components of $\Gamma_{\mathrm{tr},{\mathbf{k}}}$. Similarly, $\Gamma_{\mathrm{sym},{\mathbf{k}}}$ is a subgraph of $\Gamma_{\mathrm{tr},{\mathbf{k}+1}}$ and so the connected components of $\Gamma_{\mathrm{tr},{\mathbf{k}+1}}$ are unions of connected components of $\Gamma_{\mathrm{sym},{\mathbf{k}}}$. What we want to show, in both cases, is that the way these components decompose into smaller components is consistent with the way we label the components by their sinks $\eta_{\mathbf{k}}(\lambda)$.

That the connected components of $\Gamma_{\mathrm{sym},{\mathbf{k}}}$ break into connected components of $\Gamma_{\mathrm{tr},{\mathbf{k}}}$ in a way consistent with the map $\eta_{\mathbf{k}}$ turns out to be a simple consequence of the fact that these connected components contain parabolic coset orbits (i.e., a consequence of Corollary~\ref{cor:symccsweylorbit} from the previous section). This is established in the next lemma and corollary.

\begin{lemma} \label{lem:symdecompose}
For $\lambda,\mu \in P$, if $\lambda$ and $\mu$ belong to the same connected component of~$\Gamma_{\mathrm{sym},0}$, then $\eta_{\mathbf{k}}(\lambda)$ and $\eta_{\mathbf{k}}(\mu)$ belong to the same connected component of $\Gamma_{\mathrm{sym},{\mathbf{k}}}$ for all~$\mathbf{k} \in \mathbb{N}[\Phi]^W$.
\end{lemma}
\begin{proof}
Let $\lambda,\mu \in P$ belong to the same connected component of $\Gamma_{\mathrm{sym},0}$. From Corollary~\ref{cor:symccsweylorbit}, we get that $\mu_{\mathrm{dom}} = \lambda_{\mathrm{dom}}$ and also that there is some $w\in w_{\mu}W_{I^{0,1}_{\lambda}}$ such that $w^{-1}(\lambda)$ is dominant. But by Corollary~\ref{cor:cosets} this means $w \in w_{\lambda}W_{I^{0}_{\lambda{\mathrm{dom}}}}$, and since the cosets of $W_{I^{0,1}_{\lambda}}$ are unions of cosets of $W_{I^{0}_{\lambda}}$, this means $w_{\mu}W_{I^{0,1}_{\lambda}} = w_{\lambda}W_{I^{0}_{\lambda{\mathrm{dom}}}}$. Thus, Corollary~\ref{cor:symccsweylorbit} tells us that indeed $\eta_{\mathbf{k}}(\lambda)$ and $\eta_{\mathbf{k}}(\mu)$ belong to the same connected component of $\Gamma_{\mathrm{sym},{\mathbf{k}}}$ for all $\mathbf{k} \in \mathbb{N}[\Phi]^W$.
\end{proof}

\begin{cor} \label{cor:symdecompose}
For all $\mu \in P$ and all good $\mathbf{k} \in \mathbb{N}[\Phi]^W$, we have
\[ s^{\mathrm{sym}}_{\mathbf{k}}(\mu) = s^{\mathrm{sym}}_{0}(s^{\mathrm{tr}}_{\mathbf{k}}(\mu)).\]
\end{cor}
\begin{proof}
Since $\Gamma_{\mathrm{tr},{\mathbf{k}}}$ is a subgraph of $\Gamma_{\mathrm{sym},{\mathbf{k}}}$, the $\raU{\mathrm{sym},{\mathbf{k}}}$-stabilization of $\mu$ is the same as the $\raU{\mathrm{sym},{\mathbf{k}}}$-stabilization of the $\raU{\mathrm{tr},{\mathbf{k}}}$-stabilization of $\mu$. But the $\raU{\mathrm{tr},{\mathbf{k}}}$-stabilization of $\mu$ is by definition $\eta_{\mathbf{k}}(\lambda)$ where $\lambda \coloneqq  s^{\mathrm{tr}}_{\mathbf{k}}(\mu)$. Let $\lambda'$ be the sink of the connected component of $\Gamma_{\mathrm{sym},0}$ containing $\lambda$; hence, $\lambda' =  s^{\mathrm{sym}}_{0}(\lambda)$. Then Lemma~\ref{lem:symdecompose} says that $\eta_{\mathbf{k}}(\lambda')$ is the sink of the connected component of $\Gamma_{\mathrm{sym},{\mathbf{k}}}$ containing $\eta_{\mathbf{k}}(\lambda)$. In other words, the $\raU{\mathrm{sym},{\mathbf{k}}}$-stabilization of $\lambda$ is $\eta_{\mathbf{k}}(\lambda')$, i.e., $s^{\mathrm{sym}}_{\mathbf{k}}(\mu)=\lambda'= s^{\mathrm{sym}}_{0}(s^{\mathrm{tr}}_{\mathbf{k}}(\mu))$.
\end{proof}

We want an analog of Lemma~\ref{lem:symdecompose} and Corollary~\ref{cor:symdecompose} for truncated interval-firing. But to show that the connected components of $\Gamma_{\mathrm{tr},{\mathbf{k}+1}}$ break into connected components of $\Gamma_{\mathrm{sym},{\mathbf{k}}}$ in a way consistent with the map $\eta_{\mathbf{k}}$ turns out to be much more involved. In fact, for technical reasons, we are able to achieve this only assuming that~$\Phi$ is simply laced. Nevertheless, the first few steps towards giving truncated analogs of Lemma~\ref{lem:symdecompose} and Corollary~\ref{cor:symdecompose} do not require the assumption that $\Phi$ be simply laced, so we state them for general~$\Phi$.

\begin{prop} \label{prop:trmove0}
Let $\lambda \in P$ be such that $\<\lambda,\alpha^\vee\> \neq -1$ for all $\alpha \in \Phi^{+}$. Suppose that $\lambda \raU{\mathrm{tr},1} \lambda + \beta$ for some $\beta \in \Phi^{+}$. Then $\lambda \raU{\mathrm{tr},1} \lambda + w_{\lambda}(\alpha_i)$ for some simple root $\alpha_i$. Moreover, in this case we have $\eta_{\mathbf{k}}(\lambda) \raU{\mathrm{tr},{\mathbf{k}+1}} \eta_{\mathbf{k}}(\lambda) + w_{\lambda}(\alpha_i)$ for all $\mathbf{k} \in \mathbb{N}[\Phi]^W$.
\end{prop}
\begin{proof}
If $\<\lambda,\alpha^\vee\> \neq -1$ for all $\alpha \in \Phi^{+}$, but $\lambda \raU{\mathrm{tr},1} \lambda + \beta$ for some $\beta \in \Phi^{+}$, this must mean that $\<\lambda,\beta^\vee\> = 0$. Applying $w^{-1}_{\lambda}$, we get $\<w_{\lambda}^{-1}(\lambda),w_{\lambda}^{-1}(\beta)^\vee\> = 0$. Since $w_{\lambda}^{-1}(\beta)^\vee$ is either a positive sum or a negative sum of simple coroots, and because $w_{\lambda}^{-1}(\lambda)=\lambda_{\mathrm{dom}}$ is dominant, this means there is some simple root $\alpha_i$ such that $\<w_{\lambda}^{-1}(\lambda),\alpha_i^\vee\>=0$. But then~$\<\lambda,w_{\lambda}(\alpha_i)^\vee\>=0$. And note by Proposition~\ref{prop:posparaboliccosets} that indeed $w_{\lambda}(\alpha)$ is positive.

To prove the last sentence of the proposition: note that 
\[\<\eta_{\mathbf{k}}(\lambda),w_{\lambda}(\alpha_i)^\vee\> = \<\lambda+w_{\lambda}(\rho_{\mathbf{k}}),w_{\lambda}(\alpha_i)^\vee\> = \<w_{\lambda}^{-1}(\lambda),\alpha_i^\vee\> + \<\rho_{\mathbf{k}},\alpha_i^\vee\> = 0+\mathbf{k}(\alpha)=\mathbf{k}(\alpha);\] 
so indeed, $\eta_{\mathbf{k}}(\lambda) \raU{\mathrm{tr},{\mathbf{k}+1}} \eta_{\mathbf{k}}(\lambda) + w_{\lambda}(\alpha_i)$.
\end{proof}

\begin{prop} \label{prop:trtraps}
Let $\lambda \in P$ be a weight such that $\<\lambda,\alpha^\vee\> \neq -1$ for all $\alpha \in \Phi^{+}$. Let~$\mathbf{k} \in \mathbb{N}[\Phi]^W$ be good, with $\mathbf{k} \geq 1$. Let $\mu \in w_{\lambda} \Pi^Q_{I^{0,1}_{\lambda}}(\eta_{\mathbf{k}}(\lambda_{\mathrm{dom}}))$. Then $\mu$ and $\eta_{\mathbf{k}}(\lambda)$ belong to the same connected component of $\Gamma_{\mathrm{tr},{\mathbf{k}+1}}$.
\end{prop}
\begin{proof}
First let us prove this proposition when $\lambda$ is dominant and $I^{0}_{\lambda} = [n]$. In this case, $\eta_{\mathbf{k}}(\lambda) \in \Pi(\rho_{\mathbf{k}+1})$.  Let $\omega\in \Omega_m^0$ be such that $\lambda\in Q+\rho+\omega$. Note that, since $\mathbf{k}\geq 1$, $\eta_{\mathbf{k}}(\lambda) - \omega$ is still dominant; hence, because $P^{\mathbb{R}}_{\geq 0}\subseteq Q^{\mathbb{R}}_{\geq 0}$, we get that $\eta_{\mathbf{k}}(\lambda) - \omega \in \Pi(\rho_{\mathbf{k}+1})$ by Proposition~\ref{prop:perm_containment}. But then by definition of $\omega$ we have that~$\eta_{\mathbf{k}}(\lambda) \in \Pi^Q(\rho_{\mathbf{k}+1})+\omega$. Thus by Lemma~\ref{lem:trccs} the connected component of~$\Gamma_{\mathrm{tr},{\mathbf{k}+1}}$ that $\eta_{\mathbf{k}}(\lambda)$ belongs to is $\Pi^Q(\rho_{\mathbf{k}+1})+\omega$. By Corollary~\ref{cor:symccsweylorbit}, the connected component of $\Gamma_{\mathrm{sym},{\mathbf{k}}}$ that $\eta_{\mathbf{k}}(\lambda)$ belongs to contains the Weyl orbit~$W(\eta_{\mathbf{k}}(\lambda))$. Hence also the the connected component of $\Gamma_{\mathrm{tr},{\mathbf{k}+1}}$ that~$\eta_{\mathbf{k}}(\lambda)$ belongs to contains $W(\eta_{\mathbf{k}}(\lambda))$. But this connected component is, as mentioned,~$\Pi^Q(\rho_{\mathbf{k}+1})+\omega$; in particular, it is a convex set intersected with $Q+\eta_{\mathbf{k}}(\lambda)$. Since~$\mu$ belongs to the convex hull of $W(\eta_{\mathbf{k}}(\lambda))$ and belongs to the coset~$Q+\eta_{\mathbf{k}}(\lambda)$, this means that~$\mu \in \Pi^Q(\rho_{\mathbf{k}+1})+\omega$. So indeed $\mu$ and~$\eta_{\mathbf{k}}(\lambda)$ belong to the same connected component of $\Gamma_{\mathrm{tr},{\mathbf{k}+1}}$ in this case.

Now let us address general $\lambda$. Note that $w_{\lambda}\Phi^{+}_{I^{0,1}_{\lambda}}$ is a choice of positive roots for the sub-root system $w_{\lambda}\Phi_{I^{0,1}_{\lambda}}$. Moreover, by Proposition~\ref{prop:posparaboliccosets}, $w_{\lambda}\Phi^{+}_{I^{0,1}_{\lambda}}$ is a subset of positive roots. Hence any truncated interval-firing move (with parameter $\mathbf{k}+1$) we can carry out in $w_{\lambda}\Phi_{I^{0,1}_{\lambda}}$ with choice of positive roots $w_{\lambda}\Phi^{+}_{I^{0,1}_{\lambda}}$, we can actually carry out in the original root system $\Phi$. But then note that $\<\lambda,w_{\lambda}(\alpha_i)^\vee\> \in \{0,1\}$ for all~$i \in I^{0,1}_{\lambda}$; hence the result follows from the previous paragraph by orthogonally projecting $\lambda$ and~$\mu$ onto~$\Span_{\mathbb{R}}(w_{\lambda}\Phi_{I^{0,1}_{\lambda}})$.
\end{proof}

The strategy will be to use Proposition~\ref{prop:trmove0} to say that whenever we have a $\raU{\mathrm{tr},{1}}$-move from a sink of $\Gamma_{\mathrm{sym},{0}}$, we have a corresponding $\raU{\mathrm{tr},{\mathbf{k}+1}}$-move from the corresponding sink of $\Gamma_{\mathrm{sym},{\mathbf{k}}}$; then we will apply Proposition~\ref{prop:trtraps} to say that that move actually gets us ``trapped'' in the correct connected component of $\Gamma_{\mathrm{tr},{\mathbf{k}+1}}$. But we have reached the point where to carry out this strategy we must assume that~$\Phi$ is simply laced.

\begin{prop} \label{prop:trmove1}
Suppose that $\Phi$ is simply laced. Let $\mu \in P_{\geq 0}$ be dominant. Suppose~$\mu \raU{\mathrm{tr},{1}} \lambda$ where $\lambda = \mu+\alpha_i$ for a simple root $\alpha_i$. Then $\lambda \in W_{I^{0,1}_{\lambda}}(\lambda_{\mathrm{dom}})$.
\end{prop}
\begin{proof}
If $\mu$ is dominant but $\mu \raU{\mathrm{tr},{1}} \lambda$, this must mean that $\<\mu,\alpha_i^\vee\> = 0$. Let $\Phi'$ be the irreducible sub-root system of~$\Phi_{I^{0}_{\mu}}$ that contains $\alpha_i$. Let $\theta'$ be the highest root of $\Phi'$. We claim that $\lambda_{\mathrm{dom}} = \mu +\theta'$. First of all, because $\Phi'$ is also simply laced, the Weyl group~$W'$ of $\Phi'$ acts transitively on $\Phi'$ so that there is some $w \in W'$ with~$w(\theta') = \alpha_i$. But~$W'\subseteq W_{I^{0}_{\mu}}$, the stabilizer of $\mu$, so we indeed have $w(\mu +\theta') = \mu+\alpha_i=\lambda$. Why is  $\mu +\theta'$ dominant? Let $\D$ be the Dynkin diagram of $\Phi$ (which is just an undirected graph since $\Phi$ is simply laced). For $I\subseteq[n]$ use $\D[I]$ to denote the restriction of the Dynkin diagram to the vertices in~$I$. Note that $\Phi'= \Phi_{I}$ where $I$ is (the set of vertices of) the connected component of $\D[I^{0}_{\mu}]$ containing $\alpha_i$. Hence $\theta' = \sum_{j\in I} c_j\alpha_j$ for some coefficients $c_j$. First of all, $\theta'$ is dominant in $\Phi'$, so if $j \in I$ then $\<\theta',\alpha_j^\vee\> \geq 0$ and hence certainly $\<\mu+\theta',\alpha_j^\vee\> \geq 0$. Now suppose $j \notin I$ and $j$ is not adjacent in $\D$ to any vertex in $I$; then clearly $\<\theta',\alpha_j^\vee\> = 0$ and so again $\<\mu+\theta',\alpha_j^\vee\> \geq 0$. Finally, suppose $j \notin I$ but $j$ is adjacent in $\D$ to some vertex in $I$; then, since $\Phi$ is simply laced and $\theta'$ is a positive root of $\Phi$, we certainly have $\<\theta',\alpha_j^\vee\> \geq -1$; but $\<\lambda,\alpha_j^\vee\> \geq 1$ since $j\notin I^{0}_{\mu}$, and thus $\<\lambda+\theta',\alpha_j^\vee\> \geq 0$. So indeed $\mu+\theta'$ is dominant and so~$\lambda_{\mathrm{dom}} = \mu +\theta'$, as claimed.

Suppose for a moment that $\Phi'\neq A_1$. Then, writing $\theta'=\sum_{j=1}^{n}c_j\omega_j$, we will have that $c_j \in \{0,1\}$ for all $j \in I$; this can be seen for instance by noting that these coefficients $c_j$ are precisely the number of edges between $j$ and the ``affine node'' in the affine Dynkin diagram extending~$\D[I]$ (see~\cite[VI,\S3]{bourbaki2002lie}). This means that we have~$W' \subseteq W_{I^{0,1}_{\lambda}}$, and so $w(\lambda_{\mathrm{dom}}) = \lambda$ for some $w \in W_{I^{0,1}_{\lambda}}$; or in other words, we have~$\lambda \in W_{I^{0,1}_{\lambda}}(\lambda_{\mathrm{dom}})$. On the other hand, if $\Phi'=A_1$, then actually $\theta'=\alpha_i$ and so~$\lambda = \lambda_{\mathrm{dom}}$ and the claim is clear.
\end{proof}

\begin{remark} \label{rem:decomposeproblem}
Note that Proposition~\ref{prop:trmove1} is in general false when $\Phi$ is not simply laced. For example, take $\Phi=B_2$. Then, with $\mu \coloneqq  0$ and $\lambda \coloneqq  \alpha_1$ (the long simple root, with numbering as in Figure~\ref{fig:dynkinclassification}), we have $\mu \raU{\mathrm{tr},{1}} \lambda$ but $\lambda \notin W_{I^{0,1}_{\lambda}}(\lambda_{\mathrm{dom}})$.\end{remark}

\begin{prop} \label{prop:trmove2}
Suppose that $\Phi$ is simply laced. Let $\mu \in P_{\geq 0}$ be dominant. Suppose that~$\mu \raU{\mathrm{tr},{1}} \lambda$ where $\lambda = \mu+\alpha_i$ for a simple root $\alpha_i$. Then $\eta_k(\mu) + \alpha \in \Pi^{Q}_{I^{0,1}_{\lambda}}(\eta_k(\lambda_{\mathrm{dom}}))$ for all $k \geq 0$.
\end{prop}
\begin{proof}
The statement in the case $k=0$ follows immediately from Proposition~\ref{prop:trmove1}; so assume $k \geq 1$. By Proposition~\ref{prop:trmove1} we have that $\lambda \in W_{I^{0,1}_{\lambda}}(\lambda_{\mathrm{dom}})$, which means, by Proposition~\ref{prop:perm_containment}, that $\lambda_{\mathrm{dom}}-\lambda$ is a nonnegative sum of simple roots in $I^{0,1}_{\lambda}$.  Since $\mu$ is dominant we have $\eta_k(\mu) = \mu + k\rho$. Then note that $\eta_k(\mu)+\alpha = \lambda +k\rho = \mu+k\rho+\alpha$ is actually dominant as well, because $\mu$ is dominant, and $k\rho+\alpha$ is dominant since $\Phi$ is simply laced. Further, observe that $\eta_k(\lambda_{\mathrm{dom}})-(\eta_k(\mu)+\alpha)=\lambda_{\mathrm{dom}}-\lambda$. But then the fact that $\eta_k(\lambda_{\mathrm{dom}})-(\eta_k(\mu)+\alpha)$ is a nonnegative sum of simple roots in $I^{0,1}_{\lambda}$, together with the fact that $\eta_k(\mu)+\alpha$ is dominant, implies, via Proposition~\ref{prop:perm_containment}, that we have~$\eta_k(\mu) + \alpha \in \Pi^{Q}_{I^{0,1}_{\lambda}}(\eta_k(\lambda_{\mathrm{dom}}))$.
\end{proof}

\begin{prop} \label{prop:trmove3}
Suppose that $\Phi$ is simply laced. Let $\mu \in P$ satisfy $\<\mu,\alpha^\vee\> \neq -1$ for all $\alpha\in \Phi^{+}$. Suppose that $\mu \raU{\mathrm{tr},1} \lambda$ where $\lambda = \mu + w_{\mu}(\alpha_i)$ for some simple root~$\alpha_i$. Then for all $k \geq 0$, $\eta_k(\mu)$ and $\eta_k(\lambda)$ belong to the same connected component of~$\Gamma_{\mathrm{tr},{k+1}}$.
\end{prop}
\begin{proof}
If $k=0$ the claim is obvious. So assume $k \geq 1$.

Let $\lambda'$ be the sink of the connected component of $\Gamma_{\mathrm{sym},0}$ containing $\lambda$; hence by Corollary~\ref{cor:symccsweylorbit}, we have that $\lambda' \in w_{\lambda} W_{I^{0,1}_{\lambda}}(\lambda_{\mathrm{dom}})$, so in particular $\lambda'_{\mathrm{dom}}= \lambda_{\mathrm{dom}}$. Now, if $\<\mu,\alpha^\vee\> \neq -1$ for all $\alpha\in \Phi^{+}$ and $\mu \raU{\mathrm{tr},1} \lambda$ this means that $\<\mu,w_{\mu}(\alpha_i)^\vee\> = 0$. Hence we also have $\mu_{\mathrm{dom}} \raU{\mathrm{tr},1} \mu_{\mathrm{dom}} +\alpha_i$. Then $\mu_{\mathrm{dom}} + \alpha_i \in W_{I^{0,1}_{\lambda}}(\lambda_{\mathrm{dom}})$ by Proposition~\ref{prop:trmove1}; and so by applying $w_{\mu}$ we get $\lambda \in w_{\mu}W_{I^{0,1}_{\lambda}}(\lambda_{\mathrm{dom}})$. This implies $\lambda' \in w_{\mu}W_{I^{0,1}_{\lambda}}(\lambda_{\mathrm{dom}})$, so that $(w_{\mu}w)^{-1}(\lambda')$ is dominant for some $w \in W_{I^{0,1}_{\lambda}}$. But because of Corollary~\ref{cor:cosets} that means that $w_{\mu}w=w_{\lambda'}w'$ for some $w' \in W_{I^{0}_{\lambda}}$.

By Proposition~\ref{prop:trmove2} we get that $\eta_k(\mu_{\mathrm{dom}}) + \alpha_i \in \Pi^{Q}_{I^{0,1}_{\lambda}}(\eta_k(\lambda_{\mathrm{dom}}))$. By applying~$w_{\mu}$ we get $\eta_k(\mu)+w_{\mu}(\alpha_i) \in w_{\mu} \Pi^{Q}_{I^{0,1}_{\lambda}}(\eta_k(\lambda_{\mathrm{dom}}))$. Note that since $w \in W_{I^{0,1}_{\lambda}}$, we have that~$w_{\mu} \Pi^{Q}_{I^{0,1}_{\lambda}}(\eta_k(\lambda_{\mathrm{dom}})) = w_{\mu}w \Pi^{Q}_{I^{0,1}_{\lambda}}(\eta_k(\lambda_{\mathrm{dom}}))$. Similarly,~$w' \in W^{I^{0}_{\lambda}} \subseteq W^{I^{0,1}_{\lambda}}$ implies that~$w_{\lambda'}w'\Pi^{Q}_{I^{0,1}_{\lambda}}(\eta_k(\lambda_{\mathrm{dom}}))=w_{\lambda'}\Pi^{Q}_{I^{0,1}_{\lambda}}(\eta_k(\lambda_{\mathrm{dom}}))$. Hence, we can conclude that~$\eta_k(\mu)+w_{\mu}(\alpha_i) \in w_{\lambda'}\Pi^{Q}_{I^{0,1}_{\lambda}}(\eta_k(\lambda_{\mathrm{dom}}))$. Since $\lambda'$ is a sink of $\Gamma_{\mathrm{sym},0}$ (and thus, by Lemma~\ref{lem:symsinks}, satisfies $\<\lambda',\alpha^\vee\> \neq -1$ for all~$\alpha \in \Phi^{+}$), we can apply Proposition~\ref{prop:trtraps} to conclude that~$\eta_k(\lambda')$ and~$\eta_k(\mu)+w_{\mu}(\alpha_i)$ belong to the same connected component of~$\Gamma_{\mathrm{tr},{k+1}}$.

But since $\lambda$ and $\lambda'$ belong to the same connected component of $\Gamma_{\mathrm{sym},0}$, Lemma~\ref{lem:symdecompose} tells us that $\eta_k(\lambda)$ and $\eta_k(\lambda')$ belong to the same connected component of $\Gamma_{\mathrm{sym},k}$, and hence also belong to the same connected component of $\Gamma_{\mathrm{tr},{k+1}}$. Then note by Proposition~\ref{prop:trmove0} that we have $\eta_k(\mu) \raU{\mathrm{tr},{k+1}} \eta_k(\mu)+w_{\mu}(\alpha_i)$, so $\eta_k(\mu)$ and $ \eta_k(\mu)+w_{\mu}(\alpha_i)$ belong to the same connected component of $\Gamma_{\mathrm{tr},{k+1}}$. Putting it all together, $\eta_k(\mu)$ and~$ \eta_k(\lambda)$ belong to the same connected component of $\Gamma_{\mathrm{tr},{k+1}}$, as claimed.
\end{proof}

Finally, we are able to prove the desired analogs of Lemma~\ref{lem:symdecompose} and Corollary~\ref{cor:symdecompose} in the simply laced case.

\begin{lemma} \label{lem:trdecompose}
Suppose that $\Phi$ is simply laced.  For $\lambda,\mu \in P$, if $\lambda$ and $\mu$ belong to the same connected component of~$\Gamma_{\mathrm{tr},1}$, then $\eta_k(\lambda)$ and $\eta_k(\mu)$ belong to the same connected component of $\Gamma_{\mathrm{tr},{k+1}}$ for all~$k \geq 0$.
\end{lemma}
\begin{proof}
Clearly it suffices to prove this when $\lambda$ is a sink of $\Gamma^{\mathrm{tr},1}$. So let us describe one way to compute the $\raU{\mathrm{tr},{k+1}}$-stabilization of $\eta_k(\mu)$. If $\mu$ is not a sink of $\Gamma_{\mathrm{sym},0}$, then by Lemma~\ref{lem:symdecompose} we know that $\eta_k(\mu)$ is in the same connected component of~$\Gamma_{\mathrm{sym},k}$ as $\eta_k(\mu')$, where $\mu'$ is the sink of the component of $\Gamma_{\mathrm{sym},{0}}$ containing $\mu$; so then to compute the $\raU{\mathrm{tr},{k+1}}$-stabilization of $\eta_k(\mu)$ we instead compute the $\raU{\mathrm{tr},{k+1}}$-stabilization of~$\eta_k(\mu')$. So now assume that $\mu$ is a sink of $\Gamma_{\mathrm{sym},0}$. Then, if $\mu$ is not a sink of~$\Gamma_{\mathrm{tr},1}$, by Proposition~\ref{prop:trmove0} there is a simple root~$\alpha_i$ with $\mu \raU{\mathrm{tr},1} \mu'$ where $\mu' = \mu+w_{\mu}(\alpha_i)$. By Proposition~\ref{prop:trmove3} we get that $\eta_k(\mu)$ and $\eta_k(\mu')$ are in the same connected component of $\Gamma_{\mathrm{tr},{k+1}}$; so again to compute the $\raU{\mathrm{tr},{k+1}}$-stabilization of $\eta_k(\mu)$ we instead compute the $\raU{\mathrm{tr},{k+1}}$-stabilization of $\eta_k(\mu')$. Because $\raU{\mathrm{tr},{1}}$ is terminating, this procedure will eventually terminate; in fact, it must terminate at computing the $\raU{\mathrm{tr},{k+1}}$-stabilization of $\eta_k(\mu)$ where~$\mu$ is a sink of $\Gamma_{\mathrm{tr},1}$. But there is only one sink of the connected component of~$\Gamma_{\mathrm{tr},1}$ containing $\mu$, namely,~$\lambda$; so the lemma is proved.
\end{proof}

\begin{cor} \label{cor:trdecompose}
Suppose that $\Phi$ is simply laced. Then for all $\mu \in P$ and all $k \geq 0$, we have
\[ s^{\mathrm{tr}}_{k+1}(\mu) = s^{\mathrm{tr}}_{1}(s^{\mathrm{sym}}_{k}(\mu)).\]
\end{cor}
\begin{proof}
This follows from Lemma~\ref{lem:trdecompose} in the same way that Corollary~\ref{cor:symdecompose} follows from Lemma~\ref{lem:symdecompose}. Since $\Gamma_{\mathrm{sym},{k}}$ is a subgraph of $\Gamma_{\mathrm{sym},{k+1}}$, the $\raU{\mathrm{tr},{k+1}}$-stabilization of $\mu$ is the same as the $\raU{\mathrm{tr},{k+1}}$-stabilization of the $\raU{\mathrm{sym},{k}}$-stabilization of $\mu$. But the $\raU{\mathrm{sym},{k}}$-stabilization of $\mu$ is by definition $\eta_k(\lambda)$ where $\lambda \coloneqq  s^{\mathrm{sym}}_{k}(\mu)$. Let $\eta_1(\lambda')$ be the sink of the connected component of $\Gamma_{\mathrm{tr},1}$ containing $\lambda$; hence, $\lambda' =  s^{\mathrm{tr}}_{1}(\lambda)$. Then Lemma~\ref{lem:symdecompose} says that $\eta_{k}(\eta_{1}(\lambda')) = \eta_{k+1}(\lambda')$ (this equality follows from Proposition~\ref{prop:etafacts}) is the sink of the connected component of $\Gamma_{\mathrm{tr},{k+1}}$ containing $\eta_{k}(\lambda)$. In other words, the $\raU{\mathrm{tr},{k+1}}$-stabilization of $\lambda$ is $\eta_{k+1}(\lambda')$, i.e., $s^{\mathrm{tr}}_{k+1}(\mu)=\lambda'= s^{\mathrm{tr}}_{1}(s^{\mathrm{sym}}_{k}(\mu))$.
\end{proof}

We expect that (with the appropriate care regarding the goodness of $\mathbf{k}\in \mathbb{N}[\Phi]^W$) Lemma~\ref{lem:symdecompose} and Corollary~\ref{cor:symdecompose} should hold in the non-simply laced case as well, but, as we mentioned in Remark~\ref{rem:decomposeproblem}, our method of proof does not work there.

\section{Truncated Ehrhart-like polynomials}\label{sec:tr_Ehrhart}

The existence of the truncated Ehrhart-like polynomials, in the simply laced case, follows easily from the fact that truncated components decompose into symmetric ones in a consistent way (together with the existence of the symmetric Ehrhart-like polynomials).

\begin{thm} \label{thm:trehrhart}
Suppose that $\Phi$ is simply laced. Then, for any $\lambda \in P$, for all~$k \geq 1$ the quantity $L^{\mathrm{tr}}_{\lambda}(k)$ is given by a polynomial in~$k$ with integer coefficients.
\end{thm}
\begin{proof}
By Corollary~\ref{cor:trdecompose}, for any $k\geq 1$ and any $\lambda \in P$ we have
\begin{align*}
\#(s^{\mathrm{tr}}_{k})^{-1}(\lambda) &= \#(s^{\mathrm{sym}}_{k-1})^{-1}((s^{\mathrm{tr}}_1)^{-1}(\lambda))\\
&= \sum_{\mu \in (s^{\mathrm{tr}}_1)^{-1}(\lambda)} L^{\mathrm{sym}}_{\mu}(k-1).
\end{align*}
The right-hand side of this expression is an evaluation of a polynomial (with integer coefficients) because of Theorem~\ref{thm:symehrhart}. Since this identity holds for all $k\geq 1$, we conclude that the desired polynomial $L^{\mathrm{tr}}_{\lambda}(k)$ does exist.
\end{proof}
This finishes the proof of Theorem~\ref{thm:Ehrhart_intro}.

\begin{table}
\begin{center}
\definecolor{Gray}{gray}{0.9}
\begin{tabular}{c | c | r }
$\lambda$ & $L^{\mathrm{tr}}_{\lambda}(k)$ \\  \specialrule{2.5pt}{1pt}{1pt}
\rowcolor{Gray} $0$ & $3k^2+3k+1$ \\ \hline
\rowcolor{Gray} $\omega_1$ & $3k^2+3k+1$ \\ \hline
$-\omega_1+\omega_2$ & $2k+1$ \\ \hline
$-\omega_2$ & $k+1$ \\ \hline
\rowcolor{Gray} $\omega_2$ & $3k^2+3k+1$ \\ \hline
$\omega_1-\omega_2$ & $2k+1$ \\ \hline
$-\omega_1$ & $k+1$ \\ \hline
\rowcolor{Gray} $\omega_1 + \omega_2$ & $2k+1$ \\ \hline
$-\omega_1+2\omega_2$ & $k+1$ \\ \hline
$2\omega_1 - \omega_2$ & $k+1$ \\ \hline
$-2\omega_1 + \omega_2$ & $k+1$ \\ \hline
$\omega_1-2\omega_2$ & $k+1$ \\ \hline
$-\omega_1-\omega_2$ & $1$ \\ \specialrule{2.5pt}{1pt}{1pt}
\end{tabular}
\end{center}
\caption[Truncated Ehrahrt-like polynomials for~$A_2$]{The polynomials $L^{\mathrm{tr}}_{\lambda}(k)$ for $\Phi=A_2$.} \label{tab:trpolys}
\end{table}

\begin{conj} \label{conj:trehrhart}
For any~$\Phi$ and~$\lambda \in P$, for all good $\mathbf{k} \in \mathbb{N}[\Phi]^W$ the quantity $L^{\mathrm{tr}}_{\lambda}(\mathbf{k})$ is given by a polynomial with nonnegative integer coefficients in $\mathbf{k}$.
\end{conj}

Note that the fact we can take $\mathbf{k}=0$ in Conjecture~\ref{conj:trehrhart} means that the constant term of the $L^{\mathrm{tr}}_{\lambda}(\mathbf{k})$ polynomials should be~$1$ (which, compared to the symmetric polynomials, makes them even more like Ehrhart polynomials of zonotopes). Strictly speaking, our Theorem~\ref{thm:trehrhart} does not establish that these polynomials have constant term~$1$ even in the simply laced case.

\begin{remark}
Table~\ref{tab:trpolys} records the polynomials $L^{\mathrm{tr}}_{\lambda}(k)$ for $\Phi=A_2$, for all $\lambda \in P$ with $I_{\lambda_{\mathrm{dom}}}^{0,1}=[n]$. Compare these polynomials to the graphs of the $A_2$ truncated interval-firing processes in Example~\ref{ex:rank2graphs}. In agreement with Conjecture~\ref{conj:trehrhart}, all these polynomials have constant coefficient $1$. Note that, for $\lambda \in P$ with $L^{\mathrm{sym}}_{\lambda}(k) \neq 0$, the constant term of $L^{\mathrm{sym}}_{\lambda}(k)$ is by definition equal to the number of vertices in the connected component of $\Gamma_{\mathrm{sym},{0}}$ containing $\lambda$, which by Lemma~\ref{lem:symdecompose} is also equal to the number of connected components of~$\Gamma_{\mathrm{tr},{k}}$ contained in the connected component of~$\Gamma_{\mathrm{sym},{k}}$ with sink $\eta_{k}(\lambda)$ for all $k\geq 0$.
\end{remark}

We know that Conjecture~\ref{conj:trehrhart} holds for $\lambda \in \Omega_m^0$. That is because, for $\lambda \in \Omega_m^0$, Lemma~\ref{lem:trccs} tells us that $(s^{\mathrm{tr}}_{\mathbf{k}})^{-1}(\lambda)=\Pi(\rho_{\mathbf{k}})+\lambda$, and hence $\#(s^{\mathrm{tr}}_{\mathbf{k}})^{-1}(\lambda)$ is literally the Ehrhart polynomial of a zonotope.

Polynomials with nonnegative integer coefficients occupy a special place in algebraic combinatorics. Of course it would be great, in the course of positively resolving Conjectures~\ref{conj:symehrhart} and~\ref{conj:trehrhart}, to also give a combinatorial interpretation of the coefficients of the coefficients of these polynomials. (In fact, for the symmetric polynomials, this is precisely what is done in~\cite{hopkins2018positive}.) It would also be extremely interesting to relate these polynomials to the representation theory or algebraic geometry attached to the root system~$\Phi$, and establish positivity in that way. These polynomials arose for us in the course of a purely combinatorial investigation, but it is hard to imagine that they do not have some deeper significance if they indeed have nonnegative integer coefficients.

\begin{remark} \label{rem:ehrhart_symmetry}
It is also worth considering how the stabilization maps $s^{\mathrm{sym}}_{\mathbf{k}}$ and $s^{\mathrm{tr}}_{\mathbf{k}}$ interact with the symmetries of $\Gamma^{\mathrm{un}}_{\mathrm{sym},{\mathbf{k}}}$ and $\Gamma^{\mathrm{un}}_{\mathrm{tr},{\mathbf{k}}}$ coming from Theorem~\ref{thm:symmetry}. For the symmetric stabilization maps: if $\lambda \in P$ and $w \in W^{I^{0,1}_{\lambda}}$, then it is not hard to deduce from Lemma~\ref{lem:symdecompose} that
\[ (s^{\mathrm{sym}}_{\mathbf{k}})^{-1}(w(\lambda)) = w((s^{\mathrm{sym}}_{\mathbf{k}})^{-1}(\lambda))\]
for all good $\mathbf{k} \in \mathbb{N}[\Phi]^W$. Of course this implies that
\[ L^{\mathrm{sym}}_{w(\lambda)}(\mathbf{k}) = L^{\mathrm{sym}}_{\lambda}(\mathbf{k}),\]
in this case. Meanwhile, it appears that if  $w \in C\subseteq W$ and $\varphi\colon P\to P$ is the affine map $\varphi\colon v\mapsto w(v-\rho/h)+\rho/h$, then
\[ (s^{\mathrm{tr}}_{\mathbf{k}})^{-1}(\varphi(\lambda)) = \varphi((s^{\mathrm{tr}}_{\mathbf{k}})^{-1}(\lambda))\]
for all $\lambda \in P$ and all good $\mathbf{k} \in \mathbb{N}[\Phi]^W$. But even in the simply laced case, where we have Lemma~\ref{lem:trdecompose} at our disposal, in order to conclude that $s^{\mathrm{tr}}_{k}$ indeed respects the symmetry $\varphi$ in this way, we would need to know that this is the case for $\mathbf{k}=1$; and, as we mention in the next section, we do not currently have a great understanding of $\Gamma_{\mathrm{tr},{1}}$. So to show that the truncated stabilization maps and polynomials have the expected symmetries coming from the subgroup $C$ would require some more work.
\end{remark}

\section{Iterative descriptions of the stabilization} \label{sec:iterate}

Finally, let us focus a little more on what our decomposition results tell us about the relationship between the polynomials $L^{\mathrm{sym}}_{\lambda}(\mathbf{k})$ and $L^{\mathrm{tr}}_{\lambda}(\mathbf{k})$, and between the stabilization map $s^{\mathrm{sym}}_{\mathbf{k}}$ and $s^{\mathrm{tr}}_{\mathbf{k}}$. So, let us assume that $\Phi$ is simply laced for the remainder of this section. It is clear that Corollaries~\ref{cor:symdecompose} and~\ref{cor:trdecompose} imply the following identities relating these polynomials for all $\lambda \in P$ and all $k\geq 1$:
\begin{align*}
L^{\mathrm{sym}}_{\lambda}(k) &= \sum_{\mu \in (s^{\mathrm{sym}}_0)^{-1}(\lambda)} L^{\mathrm{tr}}_{\mu}(k); \\
L^{\mathrm{tr}}_{\lambda}(k) &= \sum_{\mu \in (s^{\mathrm{tr}}_1)^{-1}(\lambda)} L^{\mathrm{sym}}_{\mu}(k-1).
\end{align*}
What is more, these corollaries also immediately imply some striking, iterative descriptions of the stabilization functions:

\begin{cor} \label{cor:iterative}
Suppose that $\Phi$ is simply laced. Then for all $\mu \in P$ and all $k \geq 1$:
\begin{itemize}
\item $s^{\mathrm{sym}}_1(\mu) = s^{\mathrm{sym}}_0(s^{\mathrm{tr}}_1(\mu))$;
\item $s^{\mathrm{sym}}_k(\mu) = (s^{\mathrm{sym}}_1)^{k}(\mu)$;
\item $s^{\mathrm{tr}}_k(\mu) = s^{\mathrm{tr}}_1((s^{\mathrm{sym}}_1)^{k-1}(\mu))$.
\end{itemize}
\end{cor}

Corollary~\ref{cor:iterative} says that the information of all of the stabilization maps is contained just in $s^{\mathrm{sym}}_0$ and $s^{\mathrm{tr}}_1$. Now,  $s^{\mathrm{sym}}_0$ is pretty simple to understand: for example, its fibers are just parabolic Weyl coset orbits (see Corollary~\ref{cor:symccsweylorbit}). So somehow all of the complexity of all truncated and symmetric interval-firing processes (or, at least all the complexity related to \emph{stabilization} for these interval-firing processes) is contained just in~$\Gamma_{\mathrm{tr},{1}}$. Admittedly, we do not understand $\Gamma_{\mathrm{tr},{1}}$ very well. It would be very interesting, for example, to try to find an explicit description of the connected components of $\Gamma_{\mathrm{tr},{1}}$.

Finally, we end the paper by discussing another surprising consequence of Corollary~\ref{cor:iterative}: for all $\lambda \in P$ and all~$k\geq 1$,
\[ \#((s^{\mathrm{sym}}_1)^{k})^{-1}(\lambda) = L^{\mathrm{sym}}_{\lambda}(k).\]
In other words, we have a map $f\colon X \to X$ from some discrete set to itself, such that the sizes $\#(f^k)^{-1}(x)$ of fibers of iterates of this map are given by polynomials (in $k$) for every point $x \in X$. In fact, we have many such maps, one for each simply laced root system. This is a very special property for a self-map of a discrete set to have. In the next two examples we show what this looks like in the simplest cases.

\begin{figure}
\begin{center}
\begin{tikzpicture}[>={Stealth[width=2mm,length=2mm,bend]},decoration={markings,mark=at position 0.6 with {\arrow{>}}}]
\def\scl{0.3}
\node[scale=1,anchor=east] (-4) at (-4,0) {$\dots$};
\node[scale=\scl,draw,circle,fill=green] (-3) at (-3,0) {};
\node[scale=\scl,draw,circle,fill=red] (-2) at (-2,0) {};
\node[scale=\scl,draw,circle,fill=green] (-1) at (-1,0) {};
\node[scale=\scl,draw,circle,fill=red] (0) at (0,0) {};
\node[scale=\scl,draw,circle,fill=green] (1) at (1,0) {};
\node[scale=\scl,draw,circle,fill=red] (2) at (2,0) {};
\node[scale=\scl,draw,circle,fill=green] (3) at (3,0) {};
\node[scale=1,anchor=west] (4) at (4,0) {$\dots$};
\draw[->,thick] (0) to [bend left] (1);
\draw[->,thick] (1) to [bend left] (0);
\draw[->,thick] (-1) -- (0);
\draw[postaction={decorate},thick] (-2) to [bend right,looseness=1.5] (1);
\draw[->,thick] (-3) -- (-2);
\draw[->,thick] (2) -- (1);
\draw[->,thick] (3) -- (2);
\draw[->,thick] (-4) -- (-3); 
\draw[->,thick] (4) -- (3); 
\node[anchor=south] at (-3.north) {\small \textcolor{black}{$1$}};
\node[anchor=south] at (-2.north) {\small  \textcolor{black}{$1$}};
\node[anchor=south] at (-1.north) {\small  \textcolor{black}{$0$}};
\node[anchor=south] at (0.north) {\small  \textcolor{black}{$k+1$}};
\node[anchor=south] at (1.north) {\small  \textcolor{black}{$k+2$}};
\node[anchor=south] at (2.north) {\small  \textcolor{black}{$1$}};
\node[anchor=south] at (3.north) {\small  \textcolor{black}{$1$}};
\node[anchor=north] at (0.south) {\textcolor{red}{$0$}};
\node[anchor=north] at (1.south) {\textcolor{green}{$\omega_1$}};
\node[anchor=north] at (2.south) {\parbox{0.5in}{\begin{center}\textcolor{red}{$\alpha_1$}\\\textcolor{red}{$=2\omega_1$}\end{center}}};
\end{tikzpicture}
\end{center}
\caption[The stabilization map for~$A_1$]{The map $s^{\mathrm{sym}}_1\colon P \to P$ for $\Phi=A_1$. We write $L^{\mathrm{sym}}_{\lambda}(k)$ above each weight~$\lambda \in P$.} \label{fig:symtabmapa1}
\end{figure}
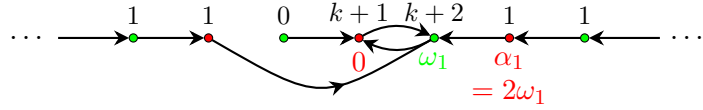

\begin{example}
Although we have so far been eschewing one-dimensional examples, in fact $s^{\mathrm{sym}}_1$ is interesting even for $A_1$. Figure~\ref{fig:symtabmapa1} depicts $s^{\mathrm{sym}}_1$ for $\Phi=A_1$. Of course in this picture we draw an arrow from $\mu$ to $\lambda$ to mean that $s^{\mathrm{sym}}_1(\mu)=\lambda$. The colors of the vertices correspond to classes of weights modulo the root lattice. We write the polynomials $L^{\mathrm{sym}}_{\lambda}(k)$ above the weights in this figure. One can verify by hand that in this case~$\#((s^{\mathrm{sym}}_1)^{k})^{-1}(\lambda) = L^{\mathrm{sym}}_{\lambda}(k)$ for all $\lambda \in P$ and all~$k \geq 0$.
\end{example}

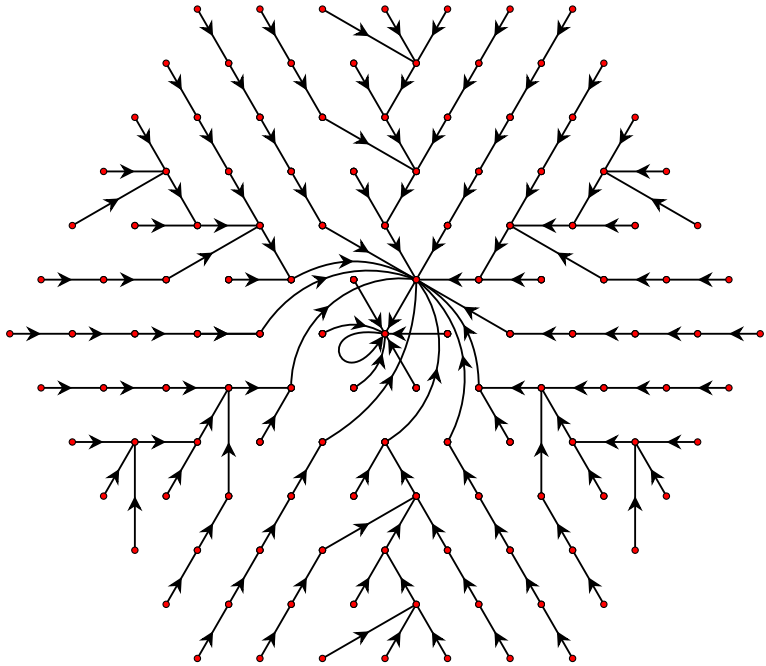
\begin{figure}
\begin{center}
\begin{tikzpicture}[scale=0.12,>={Stealth[width=2mm,length=2mm,bend]}]
\coordinate (L1) at (90:12.00);
\coordinate (L2) at (330:12.00);
\coordinate (L3) at (210:12.00);

\tikzset{
kb1/.style={postaction={decorate,decoration={markings,mark=at position 0.5 with {\arrow{>}}}},line width=0.25mm},
kb2/.style={postaction={decorate,decoration={markings,mark=at position 0.9 with {\arrow{>}}}},line width=0.25mm},   
kb3/.style={postaction={decorate,decoration={markings,mark=at position 0.7 with {\arrow{>}}}},line width=0.25mm},   
}

\draw[kb2] (barycentric cs:L1=2.0,L2=0.0,L3=1.0) -- (barycentric cs:L1=1.0,L2=1.0,L3=1.0);
\draw[kb2] (barycentric cs:L1=0.0,L2=2.0,L3=1.0) -- (barycentric cs:L1=1.0,L2=1.0,L3=1.0);
\draw[kb2] (barycentric cs:L1=2.0,L2=1.0,L3=0.0) -- (barycentric cs:L1=1.0,L2=1.0,L3=1.0);
\draw[kb2] (barycentric cs:L1=1.0,L2=2.0,L3=0.0) -- (barycentric cs:L1=1.0,L2=1.0,L3=1.0);
\draw[kb3] (barycentric cs:L1=0.0,L2=1.0,L3=2.0) to [bend right] (barycentric cs:L1=1.0,L2=1.0,L3=1.0);
\draw[kb3] (barycentric cs:L1=1.0,L2=0.0,L3=2.0) to [bend left] (barycentric cs:L1=1.0,L2=1.0,L3=1.0);

\draw[->,line width=0.25mm] (-0.1,0) to[out=170,in=240,looseness=100] (-0.2,-0.2); 

\draw[kb1] (barycentric cs:L1=3.0,L2=1.0,L3=-1.0) -- (barycentric cs:L1=2.0,L2=1.0,L3=0.0);
\draw[kb1] (barycentric cs:L1=2.0,L2=2.0,L3=-1.0) -- (barycentric cs:L1=2.0,L2=1.0,L3=0.0);
\draw[kb1] (barycentric cs:L1=1.0,L2=3.0,L3=-1.0) -- (barycentric cs:L1=2.0,L2=1.0,L3=0.0);
\draw[kb1] (barycentric cs:L1=0.0,L2=3.0,L3=0.0) to [bend right] (barycentric cs:L1=2.0,L2=1.0,L3=0.0);
\draw[kb1] (barycentric cs:L1=-1.0,L2=3.0,L3=1.0) to [bend right=40] (barycentric cs:L1=2.0,L2=1.0,L3=0.0);
\draw[kb1] (barycentric cs:L1=-1.0,L2=2.0,L3=2.0) to [bend right=50] (barycentric cs:L1=2.0,L2=1.0,L3=0.0);
\draw[kb1] (barycentric cs:L1=-1.0,L2=1.0,L3=3.0) to [bend right] (barycentric cs:L1=2.0,L2=1.0,L3=0.0);
\draw[kb1] (barycentric cs:L1=3.0,L2=0.0,L3=0.0) -- (barycentric cs:L1=2.0,L2=1.0,L3=0.0);
\draw[kb1] (barycentric cs:L1=3.0,L2=-1.0,L3=1.0) -- (barycentric cs:L1=2.0,L2=1.0,L3=0.0);
\draw[kb1] (barycentric cs:L1=2.0,L2=-1.0,L3=2.0) to [bend left] (barycentric cs:L1=2.0,L2=1.0,L3=0.0);
\draw[kb1] (barycentric cs:L1=1.0,L2=-1.0,L3=3.0) to [bend left=40] (barycentric cs:L1=2.0,L2=1.0,L3=0.0);
\draw[kb1] (barycentric cs:L1=0.0,L2=0.0,L3=2.0) to [bend left=50] (barycentric cs:L1=2.0,L2=1.0,L3=0.0);

\draw[kb1] (barycentric cs:L1=3.0,L2=2.0,L3=-2.0) -- (barycentric cs:L1=2.0,L2=2.0,L3=-1.0);
\draw[kb1] (barycentric cs:L1=4.0,L2=2.0,L3=-3.0) -- (barycentric cs:L1=3.0,L2=2.0,L3=-2.0);
\draw[kb1] (barycentric cs:L1=5.0,L2=2.0,L3=-4.0) -- (barycentric cs:L1=4.0,L2=2.0,L3=-3.0);
\draw[kb1] (barycentric cs:L1=6.0,L2=2.0,L3=-5.0) -- (barycentric cs:L1=5.0,L2=2.0,L3=-4.0);

\draw[kb1] (barycentric cs:L1=4.0,L2=1.0,L3=-2.0) -- (barycentric cs:L1=3.0,L2=1.0,L3=-1.0);
\draw[kb1] (barycentric cs:L1=5.0,L2=1.0,L3=-3.0) -- (barycentric cs:L1=4.0,L2=1.0,L3=-2.0);
\draw[kb1] (barycentric cs:L1=6.0,L2=1.0,L3=-4.0) -- (barycentric cs:L1=5.0,L2=1.0,L3=-3.0);
\draw[kb1] (barycentric cs:L1=7.0,L2=1.0,L3=-5.0) -- (barycentric cs:L1=6.0,L2=1.0,L3=-4.0);

\draw[kb1] (barycentric cs:L1=4.0,L2=0.0,L3=-1.0) -- (barycentric cs:L1=3.0,L2=0.0,L3=-0.0);
\draw[kb1] (barycentric cs:L1=5.0,L2=0.0,L3=-2.0) -- (barycentric cs:L1=4.0,L2=0.0,L3=-1.0);
\draw[kb1] (barycentric cs:L1=6.0,L2=0.0,L3=-3.0) -- (barycentric cs:L1=5.0,L2=0.0,L3=-2.0);
\draw[kb1] (barycentric cs:L1=7.0,L2=0.0,L3=-4.0) -- (barycentric cs:L1=6.0,L2=0.0,L3=-3.0);

\draw[kb1] (barycentric cs:L1=4.0,L2=-1.0,L3=0.0) -- (barycentric cs:L1=3.0,L2=0.0,L3=0.0);
\draw[kb1] (barycentric cs:L1=5.0,L2=-1.0,L3=-1.0) -- (barycentric cs:L1=4.0,L2=0.0,L3=-1.0);
\draw[kb1] (barycentric cs:L1=6.0,L2=-1.0,L3=-2.0) -- (barycentric cs:L1=5.0,L2=-1.0,L3=-1.0);
\draw[kb1] (barycentric cs:L1=7.0,L2=-1.0,L3=-3.0) -- (barycentric cs:L1=6.0,L2=-1.0,L3=-2.0);

\draw[kb1] (barycentric cs:L1=4.0,L2=-2.0,L3=1.0) -- (barycentric cs:L1=3.0,L2=-1.0,L3=1.0);
\draw[kb1] (barycentric cs:L1=5.0,L2=-2.0,L3=0.0) -- (barycentric cs:L1=4.0,L2=0.0,L3=-1.0);
\draw[kb1] (barycentric cs:L1=6.0,L2=-2.0,L3=-1.0) -- (barycentric cs:L1=5.0,L2=-1.0,L3=-1.0);
\draw[kb1] (barycentric cs:L1=7.0,L2=-2.0,L3=-2.0) -- (barycentric cs:L1=6.0,L2=-1.0,L3=-2.0);

\draw[kb1] (barycentric cs:L1=7.0,L2=-3.0,L3=-1.0) -- (barycentric cs:L1=6.0,L2=-1.0,L3=-2.0);

\draw[kb1] (barycentric cs:L1=6.0,L2=-3.0,L3=0.0) -- (barycentric cs:L1=5.0,L2=-2.0,L3=0.0);
\draw[kb1] (barycentric cs:L1=7.0,L2=-4.0,L3=0.0) -- (barycentric cs:L1=6.0,L2=-3.0,L3=0.0);

\draw[kb1] (barycentric cs:L1=5.0,L2=-3.0,L3=1.0) -- (barycentric cs:L1=4.0,L2=-2.0,L3=1.0);
\draw[kb1] (barycentric cs:L1=6.0,L2=-4.0,L3=1.0) -- (barycentric cs:L1=5.0,L2=-3.0,L3=1.0);
\draw[kb1] (barycentric cs:L1=7.0,L2=-5.0,L3=1.0) -- (barycentric cs:L1=6.0,L2=-4.0,L3=1.0);

\draw[kb1] (barycentric cs:L1=3.0,L2=-2.0,L3=2.0) -- (barycentric cs:L1=2.0,L2=-1.0,L3=2.0);
\draw[kb1] (barycentric cs:L1=4.0,L2=-3.0,L3=2.0) -- (barycentric cs:L1=3.0,L2=-2.0,L3=2.0);
\draw[kb1] (barycentric cs:L1=5.0,L2=-4.0,L3=2.0) -- (barycentric cs:L1=4.0,L2=-3.0,L3=2.0);
\draw[kb1] (barycentric cs:L1=6.0,L2=-5.0,L3=2.0) -- (barycentric cs:L1=5.0,L2=-4.0,L3=2.0);

\draw[kb1] (barycentric cs:L1=2.0,L2=-2.0,L3=3.0) -- (barycentric cs:L1=2.0,L2=-1.0,L3=2.0);
\draw[kb1] (barycentric cs:L1=3.0,L2=-3.0,L3=3.0) -- (barycentric cs:L1=3.0,L2=-2.0,L3=2.0);
\draw[kb1] (barycentric cs:L1=4.0,L2=-4.0,L3=3.0) -- (barycentric cs:L1=3.0,L2=-3.0,L3=3.0);
\draw[kb1] (barycentric cs:L1=5.0,L2=-5.0,L3=3.0) -- (barycentric cs:L1=4.0,L2=-4.0,L3=3.0);

\draw[kb1] (barycentric cs:L1=1.0,L2=-2.0,L3=4.0) -- (barycentric cs:L1=1.0,L2=-1.0,L3=3.0);
\draw[kb1] (barycentric cs:L1=2.0,L2=-3.0,L3=4.0) -- (barycentric cs:L1=3.0,L2=-2.0,L3=2.0);
\draw[kb1] (barycentric cs:L1=3.0,L2=-4.0,L3=4.0) -- (barycentric cs:L1=3.0,L2=-3.0,L3=3.0);
\draw[kb1] (barycentric cs:L1=4.0,L2=-5.0,L3=4.0) -- (barycentric cs:L1=4.0,L2=-4.0,L3=3.0);

\draw[kb1] (barycentric cs:L1=3.0,L2=-5.0,L3=5.0) -- (barycentric cs:L1=4.0,L2=-4.0,L3=3.0);

\draw[kb1] (barycentric cs:L1=2.0,L2=-5.0,L3=6.0) -- (barycentric cs:L1=2.0,L2=-4.0,L3=5.0);
\draw[kb1] (barycentric cs:L1=2.0,L2=-4.0,L3=5.0) -- (barycentric cs:L1=2.0,L2=-3.0,L3=4.0);

\draw[kb1] (barycentric cs:L1=1.0,L2=-5.0,L3=7.0) -- (barycentric cs:L1=1.0,L2=-4.0,L3=6.0);
\draw[kb1] (barycentric cs:L1=1.0,L2=-4.0,L3=6.0) -- (barycentric cs:L1=1.0,L2=-3.0,L3=5.0);
\draw[kb1] (barycentric cs:L1=1.0,L2=-3.0,L3=5.0) -- (barycentric cs:L1=1.0,L2=-2.0,L3=4.0);
\draw[kb1] (barycentric cs:L1=1.0,L2=-2.0,L3=4.0) -- (barycentric cs:L1=1.0,L2=-1.0,L3=3.0);

\draw[kb1] (barycentric cs:L1=0.0,L2=-4.0,L3=7.0) -- (barycentric cs:L1=0.0,L2=-3.0,L3=6.0);
\draw[kb1] (barycentric cs:L1=0.0,L2=-3.0,L3=6.0) -- (barycentric cs:L1=0.0,L2=-2.0,L3=5.0);
\draw[kb1] (barycentric cs:L1=0.0,L2=-2.0,L3=5.0) -- (barycentric cs:L1=0.0,L2=-1.0,L3=4.0);
\draw[kb1] (barycentric cs:L1=0.0,L2=-1.0,L3=4.0) -- (barycentric cs:L1=0.0,L2=0.0,L3=3.0);

\draw[kb1] (barycentric cs:L1=-1.0,L2=-3.0,L3=7.0) -- (barycentric cs:L1=-1.0,L2=-2.0,L3=6.0);
\draw[kb1] (barycentric cs:L1=-1.0,L2=-2.0,L3=6.0) -- (barycentric cs:L1=-1.0,L2=-1.0,L3=5.0);
\draw[kb1] (barycentric cs:L1=-1.0,L2=-1.0,L3=5.0) -- (barycentric cs:L1=0.0,L2=-1.0,L3=4.0);
\draw[kb1] (barycentric cs:L1=-1.0,L2=0.0,L3=4.0) -- (barycentric cs:L1=0.0,L2=0.0,L3=3.0);

\draw[kb1] (barycentric cs:L1=-2.0,L2=-2.0,L3=7.0) -- (barycentric cs:L1=-1.0,L2=-2.0,L3=6.0);
\draw[kb1] (barycentric cs:L1=-2.0,L2=-1.0,L3=6.0) -- (barycentric cs:L1=-1.0,L2=-1.0,L3=5.0);
\draw[kb1] (barycentric cs:L1=-2.0,L2=0.0,L3=5.0) -- (barycentric cs:L1=0.0,L2=-1.0,L3=4.0);
\draw[kb1] (barycentric cs:L1=-2.0,L2=1.0,L3=4.0) -- (barycentric cs:L1=-1.0,L2=1.0,L3=3.0);

\draw[kb1] (barycentric cs:L1=-3.0,L2=-1.0,L3=7.0) -- (barycentric cs:L1=-1.0,L2=-2.0,L3=6.0);

\draw[kb1] (barycentric cs:L1=-4.0,L2=0.0,L3=7.0) -- (barycentric cs:L1=-3.0,L2=0.0,L3=6.0);
\draw[kb1] (barycentric cs:L1=-3.0,L2=0.0,L3=6.0) -- (barycentric cs:L1=-2.0,L2=0.0,L3=5.0);

\draw[kb1] (barycentric cs:L1=-5.0,L2=1.0,L3=7.0) -- (barycentric cs:L1=-4.0,L2=1.0,L3=6.0);
\draw[kb1] (barycentric cs:L1=-4.0,L2=1.0,L3=6.0) -- (barycentric cs:L1=-3.0,L2=1.0,L3=5.0);
\draw[kb1] (barycentric cs:L1=-3.0,L2=1.0,L3=5.0) -- (barycentric cs:L1=-2.0,L2=1.0,L3=4.0);

\draw[kb1] (barycentric cs:L1=-5.0,L2=2.0,L3=6.0) -- (barycentric cs:L1=-4.0,L2=2.0,L3=5.0);
\draw[kb1] (barycentric cs:L1=-4.0,L2=2.0,L3=5.0) -- (barycentric cs:L1=-3.0,L2=2.0,L3=4.0);
\draw[kb1] (barycentric cs:L1=-3.0,L2=2.0,L3=4.0) -- (barycentric cs:L1=-2.0,L2=3.0,L3=2.0);
\draw[kb1] (barycentric cs:L1=-2.0,L2=2.0,L3=3.0) -- (barycentric cs:L1=-1.0,L2=2.0,L3=2.0);

\draw[kb1] (barycentric cs:L1=-5.0,L2=3.0,L3=5.0) -- (barycentric cs:L1=-4.0,L2=4.0,L3=3.0);
\draw[kb1] (barycentric cs:L1=-4.0,L2=3.0,L3=4.0) -- (barycentric cs:L1=-3.0,L2=3.0,L3=3.0);
\draw[kb1] (barycentric cs:L1=-3.0,L2=3.0,L3=3.0) -- (barycentric cs:L1=-2.0,L2=3.0,L3=2.0);
\draw[kb1] (barycentric cs:L1=-2.0,L2=3.0,L3=2.0) -- (barycentric cs:L1=-1.0,L2=2.0,L3=2.0);

\draw[kb1] (barycentric cs:L1=-5.0,L2=4.0,L3=4.0) -- (barycentric cs:L1=-4.0,L2=4.0,L3=3.0);
\draw[kb1] (barycentric cs:L1=-4.0,L2=4.0,L3=3.0) -- (barycentric cs:L1=-3.0,L2=3.0,L3=3.0);
\draw[kb1] (barycentric cs:L1=-3.0,L2=4.0,L3=2.0) -- (barycentric cs:L1=-2.0,L2=3.0,L3=2.0);

\draw[kb1] (barycentric cs:L1=-5.0,L2=5.0,L3=3.0) -- (barycentric cs:L1=-4.0,L2=4.0,L3=3.0);

\draw[kb1] (barycentric cs:L1=-5.0,L2=6.0,L3=2.0) -- (barycentric cs:L1=-4.0,L2=5.0,L3=2.0);
\draw[kb1] (barycentric cs:L1=-4.0,L2=5.0,L3=2.0) -- (barycentric cs:L1=-3.0,L2=4.0,L3=2.0);

\draw[kb1] (barycentric cs:L1=-5.0,L2=7.0,L3=1.0) -- (barycentric cs:L1=-4.0,L2=6.0,L3=1.0);
\draw[kb1] (barycentric cs:L1=-4.0,L2=6.0,L3=1.0) -- (barycentric cs:L1=-3.0,L2=5.0,L3=1.0);
\draw[kb1] (barycentric cs:L1=-3.0,L2=5.0,L3=1.0) -- (barycentric cs:L1=-2.0,L2=4.0,L3=1.0);
\draw[kb1] (barycentric cs:L1=-2.0,L2=4.0,L3=1.0) -- (barycentric cs:L1=-1.0,L2=3.0,L3=1.0);

\draw[kb1] (barycentric cs:L1=-4.0,L2=7.0,L3=0.0) -- (barycentric cs:L1=-3.0,L2=6.0,L3=0.0);
\draw[kb1] (barycentric cs:L1=-3.0,L2=6.0,L3=0.0) -- (barycentric cs:L1=-2.0,L2=5.0,L3=0.0);
\draw[kb1] (barycentric cs:L1=-2.0,L2=5.0,L3=0.0) -- (barycentric cs:L1=0.0,L2=4.0,L3=-1.0);
\draw[kb1] (barycentric cs:L1=-1.0,L2=4.0,L3=0.0) -- (barycentric cs:L1=0.0,L2=3.0,L3=0.0);

\draw[kb1] (barycentric cs:L1=-3.0,L2=7.0,L3=-1.0) -- (barycentric cs:L1=-1.0,L2=6.0,L3=-2.0);
\draw[kb1] (barycentric cs:L1=-2.0,L2=6.0,L3=-1.0) -- (barycentric cs:L1=-1.0,L2=5.0,L3=-1.0);
\draw[kb1] (barycentric cs:L1=-1.0,L2=5.0,L3=-1.0) -- (barycentric cs:L1=0.0,L2=4.0,L3=-1.0);
\draw[kb1] (barycentric cs:L1=0.0,L2=4.0,L3=-1.0) -- (barycentric cs:L1=0.0,L2=3.0,L3=0.0);

\draw[kb1] (barycentric cs:L1=-2.0,L2=7.0,L3=-2.0) -- (barycentric cs:L1=-1.0,L2=6.0,L3=-2.0);
\draw[kb1] (barycentric cs:L1=-1.0,L2=6.0,L3=-2.0) -- (barycentric cs:L1=-1.0,L2=5.0,L3=-1.0);
\draw[kb1] (barycentric cs:L1=0.0,L2=5.0,L3=-2.0) -- (barycentric cs:L1=0.0,L2=4.0,L3=-1.0);

\draw[kb1] (barycentric cs:L1=-1.0,L2=7.0,L3=-3.0) -- (barycentric cs:L1=-1.0,L2=6.0,L3=-2.0);

\draw[kb1] (barycentric cs:L1=0.0,L2=7.0,L3=-4.0) -- (barycentric cs:L1=0.0,L2=6.0,L3=-3.0);
\draw[kb1] (barycentric cs:L1=0.0,L2=6.0,L3=-3.0) -- (barycentric cs:L1=0.0,L2=5.0,L3=-2.0);

\draw[kb1] (barycentric cs:L1=1.0,L2=7.0,L3=-5.0) -- (barycentric cs:L1=1.0,L2=6.0,L3=-4.0);
\draw[kb1] (barycentric cs:L1=1.0,L2=6.0,L3=-4.0) -- (barycentric cs:L1=1.0,L2=5.0,L3=-3.0);
\draw[kb1] (barycentric cs:L1=1.0,L2=5.0,L3=-3.0) -- (barycentric cs:L1=1.0,L2=4.0,L3=-2.0);
\draw[kb1] (barycentric cs:L1=1.0,L2=4.0,L3=-2.0) -- (barycentric cs:L1=1.0,L2=3.0,L3=-1.0);

\draw[kb1] (barycentric cs:L1=2.0,L2=6.0,L3=-5.0) -- (barycentric cs:L1=2.0,L2=5.0,L3=-4.0);
\draw[kb1] (barycentric cs:L1=2.0,L2=5.0,L3=-4.0) -- (barycentric cs:L1=2.0,L2=4.0,L3=-3.0);
\draw[kb1] (barycentric cs:L1=2.0,L2=4.0,L3=-3.0) -- (barycentric cs:L1=3.0,L2=2.0,L3=-2.0);
\draw[kb1] (barycentric cs:L1=2.0,L2=3.0,L3=-2.0) -- (barycentric cs:L1=2.0,L2=2.0,L3=-1.0);

\draw[kb1] (barycentric cs:L1=3.0,L2=5.0,L3=-5.0) -- (barycentric cs:L1=4.0,L2=3.0,L3=-4.0);
\draw[kb1] (barycentric cs:L1=3.0,L2=4.0,L3=-4.0) -- (barycentric cs:L1=3.0,L2=3.0,L3=-3.0);
\draw[kb1] (barycentric cs:L1=3.0,L2=3.0,L3=-3.0) -- (barycentric cs:L1=3.0,L2=2.0,L3=-2.0);

\draw[kb1] (barycentric cs:L1=4.0,L2=4.0,L3=-5.0) -- (barycentric cs:L1=4.0,L2=3.0,L3=-4.0);
\draw[kb1] (barycentric cs:L1=4.0,L2=3.0,L3=-4.0) -- (barycentric cs:L1=3.0,L2=3.0,L3=-3.0);

\draw[kb1] (barycentric cs:L1=5.0,L2=3.0,L3=-5.0) -- (barycentric cs:L1=4.0,L2=3.0,L3=-4.0);

\foreach \x in {-3,-2,...,3}{
	\foreach \y in {-3,-2,...,3}{
		\foreach \z in {-3,-2,...,3}{
			\draw[fill=red] (barycentric cs:L1=1.0+\x+\z,L2=1.0-\x+\y,L3=1.0-\y-\z) circle (0.35);
		}
	}
}

\end{tikzpicture}
\end{center}
\caption[The stabilization map for~$A_2$]{The map $s^{\mathrm{sym}}_1\colon Q\to Q$ for $\Phi=A_2$. The origin is the central point (i.e., the one with a loop). The root $\alpha_1+\alpha_2$ is the point immediately north-east of the origin.} \label{fig:symtabmapa2}
\end{figure}

\begin{example}
Note that when $\Phi=A_2$, we have $\rho \in Q$ and hence $s^{\mathrm{sym}}_1$ preserves the root lattice and so descends to a map $s^{\mathrm{sym}}_1\colon Q\to Q$. Figure~\ref{fig:symtabmapa2} depicts $s^{\mathrm{sym}}_1\colon Q \to Q$ for $\Phi=A_2$. (As with our previous drawings for rank~$2$ interval-firing processes, we of course only depict the ``interesting,'' finite portion of this function near the origin.) Compare this figure to the symmetric interval-firing graphs for $A_2$ in Example~\ref{ex:rank2graphs} and the polynomials $L^{\mathrm{sym}}_{\lambda}(k)$ for $A_2$ recorded in Table~\ref{tab:sympolys}. Observe that indeed $((s^{\mathrm{sym}}_1)^k)^{-1}(0) = \Pi^Q(k\rho)$ for all $k\geq 1$. Also observe that $((s^{\mathrm{sym}}_1)^k)^{-1}(\alpha_1+\alpha_2)$ is the set of $Q$-lattice points on the boundary of $\Pi((k+1)\rho)$.
\end{example}

\bibliography{interval_firing}{}
\bibliographystyle{alpha}

\end{document}